\documentclass[11pt]{article}

\title{High-frequency analysis of parabolic stochastic PDEs with multiplicative noise}

\author{Carsten Chong\thanks{Institut de mathématiques, École Polytechnique Fédérale de Lausanne, Station 8, CH-1015 Lausanne, e-mail: carsten.chong@epfl.ch}}

\date{}
\usepackage{latexsym}
\usepackage{amssymb}
\usepackage{amsmath}
\usepackage{amsthm}
\usepackage{mathrsfs}
\usepackage[numbers]{natbib}
\usepackage{stmaryrd}
\usepackage{url}
\usepackage{longtable}
\usepackage{color}
\usepackage{dsfont}
\usepackage[T1]{fontenc}
\usepackage{lmodern}
\usepackage[linktocpage]{hyperref}
\usepackage{mathtools}

\usepackage{graphicx}
\usepackage[latin1]{inputenc}
\usepackage{enumitem}
\usepackage{amsopn}
\usepackage{comment}
\usepackage{breqn}
\newcommand{\bfi}{\begin{fig}}
\newcommand{\efi}{\end{fig}}
\newcommand{\btab}{\begin{tab}}
\newcommand{\etab}{\end{tab}}
\newcommand{\barr}{\begin{array}}
\newcommand{\earr}{\end{array}}
\newcommand{\beq}{\begin{equation}}
\newcommand{\eeq}{\end{equation}}
\newcommand{\bdis}{\begin{displaymath}}
\newcommand{\edis}{\end{displaymath}\noindent}

\newcommand{\bbn}{\mathbb{N}}

\newcommand{\bbe}{\mathbb{E}}

\newcommand{\bone}{\mathds 1}

\newcommand{\R}{\mathbb{R}}
\newcommand{\B}{\mathbb{B}}
\newcommand{\N}{\mathbb{N}}
\newcommand{\E}{\mathbb{E}}
\renewcommand{\P}{\mathbb{P}}

\newcommand{\limd}{\stackrel{\rm d}{\longrightarrow}}

\newcommand{\limst}{\stackrel{\mathrm{st}}{\Longrightarrow}}
\newcommand{\lims}{\stackrel{\mathrm{st}}{\longrightarrow}}
\newcommand{\limL}{\stackrel{L^1}{\Longrightarrow}}
\newcommand{\limp}{\stackrel{\mathbb{P}}{\longrightarrow}}

\newcommand{\lec}{\lesssim}

\newcommand{\calc}{{\cal C}}

\newcommand{\cali}{{\cal I}}

\newcommand{\calf}{{\cal F}}

\newcommand{\calb}{{\cal B}}

\newcommand{\calz}{{\cal Z}}

\newcommand{\al}{{\alpha}}
\newcommand{\la}{{\lambda}}
\newcommand{\La}{{\Lambda}}
\newcommand{\eps}{{\epsilon}}

\newcommand{\ga}{{\gamma}}
\newcommand{\Ga}{{\Gamma}}

\newcommand{\si}{{\sigma}}
\newcommand{\om}{{\omega}}
\newcommand{\Om}{{\Omega}}
\renewcommand{\phi}{\varphi}

\newcommand{\cov}{{\mathrm{Cov}}}

\newcommand{\ov}{\overline}
\newcommand{\un}{\underline}
\newcommand{\wh}{\widehat}
\newcommand{\wt}{\widetilde}

\newcommand{\dd}{\mathrm{d}}
\newcommand{\ee}{\mathrm{e}}
\newcommand{\ii}{\mathrm{i}}

\newcommand{\bb}{\mathrm{b}}
\newcommand{\cc}{\mathrm{c}}

\newcommand{\sumt}{\sum_{i=1}^{t(n)}}

\newcommand{\sumtla}{\sum_{i=\lambda_n+1}^{t(n)}}
\newcommand{\del}{\Delta_n}
\newcommand{\Del}{\Delta_n}
\newcommand{\Delh}{\sqrt{\Delta_n}}
\newcommand{\delh}{\sqrt{\Delta_n}}

\newcommand{\qtimes}{\mathrel{\hphantom{:=} \mathrel{\times}}}
\newcommand{\qminus}{\mathrel{\hphantom{:=} \mathrel{-}}}

\newcommand{\cquad}{\mathrel{\hphantom{:=}}\quad}
\newcommand{\cqquad}{\mathrel{\hphantom{:=}}\qquad}

\theoremstyle{plain}
\newtheorem{Theorem}{Theorem}[section]
\newtheorem{Corollary}[Theorem]{Corollary}
\newtheorem{Lemma}[Theorem]{Lemma}
\newtheorem{Proposition}[Theorem]{Proposition}
\theoremstyle{definition}
\newtheorem{Definition}[Theorem]{Definition}
\newtheorem{Example}[Theorem]{Example}
\newtheorem{Assumption}{Assumption}

\theoremstyle{remark}
\newtheorem{Remark}[Theorem]{Remark}

\newcommand{\bthm}{\begin{Theorem}}
\newcommand{\ethm}{\end{Theorem}}

\newcommand{\bcor}{\begin{Corollary}}
\newcommand{\ecor}{\end{Corollary}}

\newcommand{\blem}{\begin{Lemma}}
\newcommand{\elem}{\end{Lemma}}

\newcommand{\bprop}{\begin{Proposition}}
\newcommand{\eprop}{\end{Proposition}}

\newcommand{\bdf}{\begin{Definition}}
\newcommand{\edf}{\end{Definition}}

\newcommand{\bex}{\begin{Example}}
\newcommand{\eex}{\end{Example}}

\newcommand{\brem}{\begin{Remark}}
\newcommand{\erem}{\end{Remark}}

\newcommand{\bass}{\begin{Assumption}}
\newcommand{\eass}{\end{Assumption}}

\newcommand{\bpr}{\begin{proof}}
\newcommand{\epr}{\end{proof}}

\newcommand{\benu}{\begin{enumerate}}
\newcommand{\eenu}{\end{enumerate}}

\newcommand{\bit}{\begin{itemize}}
\newcommand{\eit}{\end{itemize}}

\numberwithin{equation}{section}

\allowdisplaybreaks[3]

\begin{document}
	
	\maketitle
	
	\begin{abstract} We consider the stochastic heat equation driven by a multiplicative Gaussian noise that is white in time and spatially homogeneous in space. Assuming that the spatial correlation function is given by a Riesz kernel of order $\alpha \in (0,1)$, we prove a central limit theorem for power variations and other related functionals of the solution. 
	To our surprise, there is no asymptotic bias despite the low regularity of the noise coefficient in the multiplicative case. We trace this circumstance back to cancellation effects between error terms arising naturally in second-order limit theorems for power variations.
	\end{abstract}
	
	\noindent
	\begin{tabbing}
	\end{tabbing}

	
	\vspace{0.5cm}

\section{Introduction}\label{sect:intro}

Consider the stochastic heat equation
\beq\label{eq:SHE-si}\begin{aligned}
	\partial_t u(t,x) &= \frac12 \Delta u(t,x)+\si(u(t,x))\dot W(t,x),&&(t,x)\in(0,\infty)\times\R^d,\\
	u(0,x)&=u_0(x),&&x\in\R^d,\end{aligned}\eeq
subject to a sufficiently regular initial condition $u_0$. Here, $\si\colon\R\to\R$ is a globally Lipschitz function and $\dot W$ is a centered Gaussian noise that is white in time and spatially homogeneous with spectral measure $\mu(\dd \xi):=|\xi|^{\al-d}\,\dd\xi$ for some $\al\in(0,1]$. In particular, if $d=\al=1$, $\dot W$ is a space-time white noise. We will recall the precise definition of the noise as well as the solution theory for \eqref{eq:SHE-si} in Section~\ref{sect:main} below. 

In this paper, we are interested in the \emph{normalized power variations} (and other related functionals) of the solution $u(t,x)$, defined as
\beq\label{eq:Vnp} V^n_p(u,t) := \del \sum_{i=1}^{[t/\del]} \biggl|\frac{ u(i\Del,x)-u((i-1)\Del,x)}{\tau_n}\biggr|^p,\qquad t\in(0,\infty), \eeq
where $p>0$, $x\in\R^d$ is an arbitrary but fixed spatial point, $[\cdot]$ denotes the integer part, $\Del$ is a small step size that converges to $0$ as $n\to\infty$ (e.g., $\Del=\frac1n$), and $\tau_n=\tau_n(\al)$ is a normalizing sequence depending on $\al$, chosen in such a way that $V^n_p(u,t)$ converges as $n\to\infty$. Indeed, if $\tau_n$ is given as in \eqref{eq:taun}, then $V^n_p(u,t)$ satisfies the following \emph{law of large numbers (LLN)} as $\Del\to0$ (see \cite{Foondun15, Pospisil07} if $\al=d=1$ and $p\in\{2,4\}$ and  \cite[Theorem~2.1 and Remark~2.7]{Chong19} for the general case):
\beq\label{eq:LLN-0} V^n_p(u,t) \limL V_p(u,t):=m_p \int_0^t |\si(u(s,x))|^p \,\dd s, \eeq 
where $m_p$ is the $p$th absolute moment of the standard normal distribution and $\limL$ denotes uniform $L^1$-convergence on any compact interval $[0,T]$ with $T\in(0,\infty)$. 

The main goal of this paper is to study the second-order behavior of $V^n_p(u,t)$, that is, we want to show a \emph{central limit theorem (CLT)} of the type
\beq\label{eq:CLT-0} \frac1{\Delh}\Bigl(V^n_p(u,t) - V_p(u,t)\Bigr)\limst \calz, \eeq
where $\calz$ is some mixed Gaussian process and $\limst$ stands for stable convergence in law in the space of càdlàg functions with the local uniform topology (again, we will recall this in Section~\ref{sect:main}).

A CLT like \eqref{eq:CLT-0} plays a fundamental role in statistical estimation problems related to \eqref{eq:SHE-si}. Consider,  
for example, the \emph{parabolic Anderson model} [i.e., $\si(x)=\si_0 x$ for some $\si_0>0$] and suppose that $\si_0$ is an unknown parameter that we would like to estimate based on the \emph{high-frequency} observations $u(i\Del,x)$, $i=1,\dots,[T/\Del]$, where $T>0$ is arbitrary but fixed. Then \eqref{eq:LLN-0}, together with classical Riemann sum approximation, readily shows that 
\[ (\si^n_0)^p:=\frac{V^n_p(u,T)}{m_p\Del \sum_{i=1}^{[T/\Del]} {|u(i\Del,x)|}^p} \limp \si_0^p;\]
see also \cite[Section~3]{Cialenco18} in the case $d=\al=1$.
In other words, the left-hand side is a \emph{weakly consistent estimator} of $\si_0^p$ (if the value of $\al$ is known; otherwise, in order to be able to compute $\tau_n=\tau_n(\al)$ and hence $V^n_p(u,T)$ from our observations, we have to replace $\tau_n$ by an appropriate estimator first; see \cite[Section~2.3]{Chong19}). But if we want to construct \emph{asymptotic confidence regions} and \emph{statistical tests} for $\si_0$, then we need to study the second-order fluctuations of $(\si^n_0)^p$, and the CLT for $V^n_p(u,t)$ is a crucial step for this purpose. We will pursue this direction, which, in addition to the CLT \eqref{eq:CLT-0}, also involves analyzing the discretization error of the Riemann sum approximation in the denominator of $(\si_0^n)^p$, in a separate work.

Let us briefly recall the existing literature in this framework. The idea of using quadratic variation [i.e., $p=2$ in \eqref{eq:Vnp}] to estimate the variance of a Gaussian process is, of course, as old as stochastic analysis itself and originates from the observation that the quadratic variation of $\si_0 B$  on $[0,1]$, where $B$ is a standard Brownian motion, is precisely $\si_0^2$. In the last decade, this basic idea was developed into a systematic theory of power variations and related  functionals for a large class of semimartingales.
The resulting statistical applications have attracted a lot of interest, especially in the context of \emph{volatility estimation} in financial econometrics. The interested reader is referred to the two monographs \cite{AitSahalia14, Jacod12} for more information on this aspect.

As it turns out, power variations are very robust in the sense that their first- and second-order behavior remains (essentially) the same for many processes that are not semimartingales. Examples include fractional Brownian motion with Hurst parameter less than $\frac34$ (see \cite{Corcuera06,Istas97}) or, more generally, moving average processes with fractional kernels (see \cite{BN11,Corcuera13, Podolskij15}). 

When it comes to the stochastic heat equation \eqref{eq:SHE-si}, clearly the process $t\mapsto u(t,x)$ (with $x$ fixed) is not a semimartingale. Indeed, if $\al=d=1$, it has a finite and nontrivial quartic variation by \cite{Swanson07,Walsh81}; in general, it has a finite and nontrivial $(\frac12-\frac\al4)^{-1}$-variation, which follows from \eqref{eq:Vnp} and the definition of $\tau_n$ in \eqref{eq:taun}. Nevertheless, as we have indicated in \eqref{eq:LLN-0}, the LLN for $V^n_p(u,t)$ holds true in a form that is similar to the semimartingale or moving average setting. More generally, if $Y$ denotes the solution to 
\beq\label{eq:SHE-2}\begin{aligned}
	\partial_t Y(t,x) &= \frac12 \Delta Y(t,x)+\sigma (t,x)\dot W(t,x),&&(t,x)\in(0,\infty)\times\R^d,\\
	Y(0,x)&=u_0(x),&&x\in\R^d,\end{aligned}\eeq
where $\sigma(t,x)$ is a predictable, $L^2$-continuous and $L^{p+\eps}$-bounded random field, then 
\beq\label{eq:LLN-1}  V^n_p(Y,t) \limL m_p\int_0^t |\sigma(s,x)|^p\,\dd s \eeq
by \cite[Theorem~2.1 and Remark~2.7]{Chong19}.
In particular, if we set $\sigma(t,x):=\si (u(t,x))$, where $u$ is the solution to \eqref{eq:SHE-si}, then $Y=u$ and \eqref{eq:LLN-0} follows from \eqref{eq:LLN-1}. 

The picture is much less complete (and, as we will see, quite different) when it comes to the second-order limit theorem. Until now, a CLT of the form \eqref{eq:CLT-0}
has only been shown if, among other things, the random field $\sigma(t,x)$ in \eqref{eq:SHE-2} is  
\bit
\item essentially $\frac12$-Hölder continuous in time and 
\item pathwise twice differentiable in space; 
\eit
see \cite[Theorem~2.3]{Chong19}. This covers the case of additive noise (previously considered also in   \cite{Bibinger17,Bibinger19, Cialenco17} when  $\al=d=1$ and $p\in\{2,4\}$) but certainly \emph{not} the case of multiplicative noise as in \eqref{eq:SHE-si}. Indeed, it is well known from \cite[Theorem~2.1]{Sanzsole02} that the solution $u$ to \eqref{eq:SHE-si} is only $(\frac12-\frac\al 4-\eps)$-Hölder continuous in time and $(1-\frac\al2-\eps)$-Hölder continuous in space. So if we let $\sigma(t,x)=\si (u(t,x))$ as above, then $\sigma(t,x)$ would \emph{not} meet the assumptions needed for the CLT in \cite[Theorem~2.3]{Chong19}.

Among the various reasons why the strong regularity conditions on $\sigma(t,x)$ above were needed for  the CLT in \cite{Chong19}, this is the most important one: at many places within the proof, we had to \emph{discretize} the random field $\sigma(t,x)$ along the grid points $t=i\Del$ [e.g., replace $\sigma(t,x)$ by $\sigma([t/\Del]\Del,x)$], so that the increments of the noise become locally (conditionally) Gaussian. Of course, if $\sigma(t,x)=\si (u(t,x))$, this discretization induces an error of order $\Del^{1/2-\al/ 4}$ by the Hölder properties of $u$, which blows up once we divide by $\sqrt{\Del}$ as in \eqref{eq:CLT-0}. Thus, we can only justify this discretization procedure if $\sigma(t,x)$ is essentially $\frac12$-Hölder continuous in time.

 The obvious question is now whether this problem is an artifact of our proof. The answer is yes and no:  if $0<\al<1$, we will show that the CLT as in \eqref{eq:CLT-0} is valid. However, in the case $\al=1$ (in dimension $1$, this is the case when $\dot W$ is a space-time white noise), it turns out that an additional asymptotic bias emerges in the CLT  for $p\geq4$, while it holds true without bias for $p=2$. In fact, the bias term for $p\geq4$ originates exactly from the discretization of $\sigma(t,x)=\si (u(t,x))$. 

In this paper, we focus on the case $0<\al<1$, where we have no asymptotic bias; the critical case $\al=1$ will be discussed in a second paper \cite{Chong19c} separately. In principle, if $d\geq2$, one could also consider \eqref{eq:SHE-si} with $\al\in(1,2)$. Whether the CLT holds in that case is open, although we strongly conjecture that it will \emph{not} hold and that the rate will no longer be $1/\Delh$. Let us briefly mention the papers \cite{Huang18, Huang19}, which, despite their titles, are not related to the CLT we consider in the present work.

The remaining article is organized as follows: In Section~\ref{sect:main}, we properly introduce the equations and variation functionals under investigation and announce our main result (Theorem~\ref{thm:CLT}). 
In Section~\ref{sect:proof}, we outline the main strategy of the proof, while the details are given in Section~\ref{sect:details}.

In what follows, we often write $\iint_a^b :=\int_a^b \int_{\R^d}$, $\iint:=\iint_0^\infty$, $\iiint_a^b := \int_a^b\int_{\R^d}\int_{\R^d}$, and $\iiint := \iiint_0^\infty$. Moreover, $\bbn:=\{1,2,\ldots\}$ and $\bbn_0:=\{0,1,2,\ldots\}$. We also write $A\lec B$ if there is a constant $C\in(0,\infty)$ that is independent of all quantities of interest such that $A\leq CB$.

\section{Set-up and main result}\label{sect:main} Let $\B:=(\Om,\calf,(\calf_t)_{t\geq0},\P)$ be a stochastic basis satisfying the usual conditions that supports an $L^2$-valued centered Gaussian measure $W(A)$, indexed by bounded Borel sets $A\in\calb_\bb([0,\infty)\times \R^d)$, such that $W(A)$ is $\calf_t$-measurable if $A\subseteq [0,t]\times\R^d$, $W(A)$ and $\calf_t$ are independent if $A\cap [0,t]\times\R^d = \emptyset$, and 
\beq\label{eq:cov} \bbe[ W(A_1)W(A_2) ] = \iiint\bone_{A_1}(s,y)\bone_{A_2}(s,y')\,\La(\dd y,\dd y')\,\dd s.   \eeq
In the last line, the spatial covariance measure is given by $\La(\dd y,\dd y') := F(y-y')\,\dd y\,\dd y'$, where $F$ is the \emph{Riesz kernel} $F(y):=c_\al |y|^{-\al}$ for some $\al\in(0,1)$ and $c_\al := \pi^{d/ 2-\al} \Ga(\frac \al 2)/\Ga(\frac{d-\al}{2})$ ($|\cdot|$ denotes the Euclidean norm on $\R^d$). 
By \cite{Dalang99}, if 
\beq\label{eq:FT} \calf \phi(\xi):=\int_{\R^d} \ee^{-2\pi \ii \xi\cdot x}\phi(x)\,\dd x,\quad \xi\in\R^d,\quad\phi\in C^\infty_\cc(\R^d),\eeq
is the \emph{Fourier transform} on the space $C^\infty_\cc(\R^d)$ of compactly supported smooth functions, then there is a measure $\mu$, called the \emph{spectral measure} of $W$, such that for all $\psi_1,\psi_2\in C^\infty_\cc([0,\infty)\times \R^d)$,
\beq\label{eq:isometry} \begin{split} &\iiint \psi_1(s,y)\psi_2(s,y')\,\La(\dd y,\dd y')\,\dd s= \iint \calf \psi_1(s,\cdot)(\xi)\ov{\calf\psi_2(s,\cdot)(\xi)}\,\mu(\dd \xi)\,\dd s, \end{split} \eeq
where $\ov{(\cdot)}$ denotes complex conjugation. With our choice of $F$, it is well known (see Chapter~V, \S 1, Lemma~1(a) in \cite{Stein70}) that the spectral measure is given by $\mu(\dd \xi) =|\xi|^{\al-d}\,\dd \xi$.

It is a classical result from \cite[Theorem~13]{Dalang99} that if $u_0$ is bounded and continuous, then \eqref{eq:SHE-si} admits a unique \emph{mild solution} $u$ that satisfies
\beq\label{eq:mom} \sup_{(t,x)\in[0,T]\times\R^d} \E[|u(t,x)|^p]<\infty\eeq
for all $T\geq0$ and $p\geq1$. Recall that a predictable process $u$ is a mild solution to \eqref{eq:SHE-si} if for every $(t,x)\in(0,\infty)\times\R^d$, we have
\beq\label{eq:mild} u(t,x) = u^{(0)}(t,x)+\iint_0^t G(t-s,x-y)\si(u(s,y)) \,W(\dd s,\dd y) \eeq
almost surely, where 
\beq\label{eq:initial} u^{(0)}(t,x):=\int_{\R^d} G(t,x-y)u_0(y)\,\dd y\eeq is the solution to the homogeneous heat equation and
\beq\label{heatkernel} G(t,x):=G_x(t):=(2\pi t)^{-\frac d2}\ee^{-\frac{|x|^2}{2 t}} \bone_{t>0},\qquad (t,x)\in[0,\infty)\times\R^d, \eeq
is the \emph{heat kernel}.

Now suppose that we have measurements of the solution $u(t,x)$ at time points $t=i\Del$, $i=1,\dots,[T/\Del]$, (where $T>0$ and $\Del\to 0$ as $n\to\infty$) and space points $x=x_1,\dots,x_K$ (where $K\in\N$ and $x_1,\dots,x_K$ are different points in $\R^d$). Given two further numbers $L,M\in\N$ and an evaluation function $f\colon \R^{K\times L} \to \R^M$, we can now consider the \emph{normalized variation functional} $V^n_f(u,t):=(V^n_f(u,t)_1,\dots,V^n_f(u,t)_M)'$ defined as follows:
\beq\label{eq:Vnf} V^n_f(u,t)_m := \Delta_n \sumt f_m\biggl( \frac{\Delta^n_i u}{\tau_n},\ldots,\frac{\Delta^n_{i+L-1} u}{\tau_n} \biggr),\quad t\in[0,T],\quad m=1,\dots,M, \eeq
where $\Delta^n_i u:= (u(i\Del,x_1)-u((i-1)\Del,x_1),\dots,u(i\Del,x_K)-u((i-1)\Del,x_K))'$ is the column vector containing the $i$th \emph{increment} of $u$ at the different points $x_1,\dots,x_K$ and 
\beq\label{eq:taun}  \tau_n:=\tau_n(\al):=\sqrt{C_\al}\Del^{\frac12-\frac\al4},\quad C_\al:=\frac{\pi^{\frac d2 - \al}\Ga(\frac\al2)}{2^{\frac\al2}(1-\frac\al2)\Ga(\frac d2)},\quad t(n):=[t/\Del]-L+1. \eeq
Typical examples include 
\bit
\item $K=L=M=1$ and $f(z):=|z|^p$: then $V^n_f(u,t)=V^n_p(u,t)$ if $x=x_1$;
\item $K=M$ and $f_k((z_{kl})_{k,l=1}^{K,L})=\prod_{l=1}^L |z_{kl}|^{p_{kl}}$ [resp., $f_k((z_{kl})_{k,l=1}^{K,L})=\prod_{l=1}^L z_{kl}^{p_{kl}}$] with $p_{kl}\geq0$ (resp., $p_{kl}\in\N_0$): then $V^n_f(u,t)$ is a \emph{multipower variation} (resp., \emph{signed multipower variation}); cf.\ \cite[Section~2.2]{Chong19}; 
\item $K=2$, $L=M=1$ and $f((\begin{smallmatrix} z_1\\ z_2 \end{smallmatrix}))= z_1z_2$: then $V^n_f(u,t)$ is a \emph{quadratic covariation}.
\eit

The following LLN was proved in \cite[Theorem~2.1 and Remark~2.7]{Chong19}:
\bthm[Law of large numbers]\label{thm:LLN} Suppose that $f$ is continuous with at most polynomial growth [i.e., there is $p\in(0,\infty)$ such that ${|f(z)|}\lec 1+|z|^p$ for all $z\in\R^{K\times L}$, where $|z|$ denotes the Euclidean norm of the matrix $z$ when viewed as a vector in $\R^{KL}$]. Further assume that the initial condition $u_0$ is bounded and Hölder continuous on $\R^d$ with some exponent larger than $1-\frac\al2$. Then
\beq\label{eq:LLN} V^n_f(u,t)\limL V_f(u,t):=\int_0^t \mu_f\Bigl(\si^2(u(s,x_1)),\ldots,\si^2(u(s,x_K)) \Bigr)\,\dd s.\eeq
Here, $\mu_f\colon [0,\infty)^K \to \R^M$, $(w_1,\dots,w_K)\mapsto \bbe[f(Z)]$, where $Z=(Z_{kl})_{k,l=1}^{K,L}$ is multivariate Gaussian with mean $0$ and 
\beq\label{eq:cov-Z}\cov(Z_{k_1l_1},Z_{k_2 l_2}) = \Ga_{|l_1-l_2|}w_k \bone_{k_1=k_2=k},\eeq
and
\beq\label{eq:Ga-formula} \Ga_0:=1,\quad \Ga_r:=\frac12 \Big( (r+1)^{1-\frac\al 2}-2r^{1-\frac\al 2}+(r-1)^{1-\frac\al 2}\Big),\quad r\geq1.  \eeq 
\ethm

In order to describe our main result, we need some technical terminology. A \emph{product filtered extension} of $\B$ is an enlarged stochastic basis $\ov\B:=(\ov\Om,\ov\calf,(\ov\calf_t)_{t\geq0},\ov\P)$ obtained by taking the product
\[ \ov\Om:=\Om\times\Om',\quad\ov\calf:=\calf\otimes \calf',\quad\ov\calf_t:=\bigcap_{s>t} \calf_s\otimes\calf'_s,\quad\ov\P:=\P\otimes\P' \]
with an additional stochastic basis $\B':=(\Om',\calf',(\calf'_t)_{t\geq0},\P')$. In what follows, $\ov\E$ denotes expectation on $\ov \B$ and random elements $X$ on $\B$ are canonically extended to $\ov\B$ by setting $X(\om,\om'):=X(\om)$ (without changing notation). Now a sequence $(X_n)_{n\in\N}$ of $E$-valued random variables defined on the original basis $\B$ (where $E$ is a metric space) is said to \emph{converge stably in law} to a variable $X$ defined on the product filtered extension $\ov\B$ (and denoted $X_n\lims X$) if for all real-valued random variables $Y$ on $\B$ and bounded continuous functions $g\colon E\to\R$, we have $\E[Yg(X_n)] \to \ov\E[Yg(X)]$ as $n\to\infty$. This clearly implies $X_n\limd X$ in the ordinary sense. But what sets it apart, is the following property:
\beq\label{eq:stable} X_n\lims X,\ Y_n\limp Y \implies (X_n,Y_n)\lims (X,Y);  \eeq 
see \cite[Proposition~2(i)]{Podolskij10}. To have stable (and not just ordinary) convergence in law in our main result below, is essential for the statistical purposes mentioned in the Introduction.

Recall that $\limst$ denotes stable convergence in law on $E:=D([0,\infty),\R^M)$, the space of $\R^M$-valued càdlàg functions, equipped with the local uniform topology.
\bthm[Central limit theorem]\label{thm:CLT} Grant the following hypotheses:
\benu
\item[H1.] The function $f\colon \R^{K\times L}\to\R^M$ is even [i.e., we have $f(z)=f(-z)$ for all $z\in\R^{K\times L}$] and four times continuously differentiable with all partial derivatives up to order four being of polynomial growth at most. 
\item[H2.] For each $m=1,\ldots, M$, $f_m(z)$ only depends on a single row in the matrix $z$, that is, there is $k=k(m)$ such that $f_m(z)$ only depends on $z_{k1},\ldots,z_{kL}$.
\item[H3.] The function $\si$ is four times differentiable with all four derivatives being of polynomial growth at most.
\item[H4.] The initial condition $u_0$ is bounded and differentiable with a derivative that is Hölder continuous with some exponent larger than $\frac12-\frac\al4$.
\eenu
Then, as $\Delta_n\to0$,
\beq\label{eq:CLT-statement} \frac1{\Delh}\Bigl( V^n_f(u,t)-V_f(u,t)\Bigr) \limst \calz, \eeq
where $\calz= ((\calz_1(t),\ldots,\calz_M(t))')_{t\geq0}$ is a continuous process defined on a product filtered extension $\ov\B$ of the original basis $\B$ with the following properties: conditionally on the $\si$-field $\calf$, $\calz$ is a centered Gaussian process with independent increments such that the covariance function $\calc_{m_1m_2}(t) := \ov\bbe[\calz_{m_1}(t) \calz_{m_2}(t)\mid\calf]$, for $m_1,m_2=1,\ldots,M$, is given by 
\beq\label{eq:calc} \begin{split}\calc_{m_1m_2}(t) &= \int_0^t \rho_{f_{m_1},f_{m_2}}(0;\si^2(u(s,x_1)),\ldots,\si^2(u(s,x_K)))\,\dd s\\
	&\quad+\sum_{r=1}^\infty \int_0^t \rho_{f_{m_1},f_{m_2}}\Bigl(r;\si^2(u(s,x_1)),\ldots,\si^2(u(s,x_K))\Bigr)\,\dd s\\
	&\quad+\sum_{r=1}^\infty \int_0^t \rho_{f_{m_2},f_{m_1}}\Bigl(r;\si^2(u(s,x_1)),\ldots,\si^2(u(s,x_K))\Bigr)\,\dd s. \end{split}\eeq 
In the last line, for $r\in\N_0$, we define 
\beq\rho_{f_{m_1},f_{m_2}}(r;w_1,\ldots,w_K):=\cov(f_{m_1}(Z^{(1)}),f_{m_2}(Z^{(2)})),\eeq 
where $Z^{(1)}=(Z^{(1)}_{kl})_{k,l=1}^{K,L}$ and $Z^{(2)}=(Z^{(2)}_{kl})_{k,l=1}^{K,L}$ are jointly Gaussian, both with the same law as the matrix $Z$ in Theorem~\ref{thm:LLN} 
and cross-covariances  
\beq\label{eq:cross-cov} \cov(Z^{(1)}_{k_1l_1}, Z^{(2)}_{k_2l_2})=\Ga_{|l_1-l_2+r|}w_k\bone_{k_1=k_2=k}. \eeq
Part of the statement is that the series in \eqref{eq:calc} converge in the $L^1$-sense.
\ethm

\brem\label{rem:1} Let us comment on the assumptions as well as the conclusions of Theorem~\ref{thm:CLT}.
\benu
	\item The conclusion of Theorem~\ref{thm:CLT} is exactly the same as that of \cite[Theorem~2.3]{Chong19}. As indicated in the Introduction, the main novelty here is that the process $\si(u(t,x))$ is neither $\frac12$-Hölder continuous in time nor pathwise twice differentiable in space. This extension might seem to be somehow expected. But in view of the fact that \eqref{eq:CLT-statement} fails when $\al=1$ (as we will demonstrate separately in \cite{Chong19c}), we find the positive result for $\al\in(0,1)$ rather surprising. 
	\item Hypothesis H2 means that we are not allowed to mix observations at different spatial points for a given coordinate of $f$ (e.g., only the first two examples before Theorem~\ref{thm:LLN} satisfy this hypothesis). This condition is already needed in the case of additive noise (i.e., when $\si\equiv 1$) for the reasons explained in \cite[Remark~2.4 (2)]{Chong19}.
	\item Hypothesis H2 has the following consequences on $\mu_f$ and $\rho_{f_{m_1},f_{m_2}}$: for each $m=1,\dots,M$, there is a single value of $k$ (which, of course, may depend on $m$) such that the $m$th coordinate of $\mu_f$ only depends on $w_k$. Also, if $f_{m_1}(z)$ and $f_{m_2}(z)$ depend on the same row of $z$ (say, the $k$th one), then $\rho_{f_{m_1},f_{m_2}}(r;w_1,\dots,w_K)$ only depends on $w_k$. And if $f_{m_1}(z)$ and $f_{m_2}(z)$ depend on different rows of $z$, then $\rho_{f_{m_1},f_{m_2}}(r;\cdot)\equiv0$. 
\eenu
\erem

\section{Overview of the proof of Theorem~\ref{thm:CLT}}\label{sect:proof}
The strategy of proof for the CLT in \cite{Chong19} can be briefly described as follows: first, the variation functionals in \eqref{eq:Vnf} are, little by little, approximated by terms that get closer and closer to martingale sums; second, for the resulting martingale array, the actual CLT is shown using  Jacod's martingale CLT from \cite{Jacod97}; finally, the centering appearing in the martingale array of the previous step is transformed back to the actual centering $V_f(u,t)$ from \eqref{eq:CLT-statement}. A close inspection of the proof reveals that the second step does not rely on the regularity properties of $\si(u(t,x))$ at all, so the difficulties due to the low regularity of $\si(u(t,x))$ all appear in the first and the third step described above.

Now here is the surprising thing: while certain error terms that arise during these two steps do not vanish on a $\Delh$-rate when considered individually, their sum does so due to cancellation effects! This is the crucial observation that allows us to overcome the lack of regularity of $\si(u(t,x))$ and ultimately prove Theorem~\ref{thm:CLT}. Clearly, the special structure that $\si(u(t,x))$ is a function of the solution itself (rather than any random field with a similar regularity) is of utmost importance for obtaining the mentioned cancellations.

In order to simplify notations in the subsequent exposition, let us introduce the following abbreviations for $0\leq a\leq b\leq c<\infty$ and $(t,x)\in (0,\infty)\times\R^d$:
\beq\label{eq:abbr} \begin{split} 
	u(t,x)_b^c&:=u^{(0)}(t,x)+\iint_{b}^{c} G_{x-y}(t-s)\si(u(s,y))\,W(\dd s,\dd y), \\
	 u(t,x)_{b,a}^c&:=u^{(0)}(t,x)+\iint_{b}^{c} G_{x-y}(t-s)\si(u(s,y)_0^a)\,W(\dd s,\dd y), \\
	 \wt u(t,x)_{b,a}^c&:=u^{(0)}(t,x)+\iint_{b}^{c} G_{x-y}(t-s)\si(u(a,y))\,W(\dd s,\dd y), \\
	\Delta^n_i u_b^c &:= \iint_{b}^{c} \Delta^n_i G_y(s)\si(u(s,y))\,W(\dd s,\dd y),\\
		\Delta^n_i  u_{b,a}^c &:= \iint_{b}^{c} \Delta^n_i G_y(s)\si(u(s,y)_0^a)\,W(\dd s,\dd y),\\
		\Delta^n_i \wt u_{b,a}^c &:= \iint_{b}^{c} \Delta^n_i G_y(s)\si(u(a,y))\,W(\dd s,\dd y), \end{split}\eeq
where $\Delta^n_i G_y(s)$ is the matrix in $\R^{K\times L}$ with 
$$\Delta^n_i G_y(s)_{kl}:= G_{x_k-y}((i+l-1)\Del-s)-G_{x_k-y}((i+l-2)\Del-s)$$  
and
$G_y(s):=0$ if $s<0$. 
For brevity, we often use the abbreviations $\nu:=(k,l)$ [and $\nu_1:=(k_1,l_1)$, $\nu_2:=(k_2,l_2)$, etc.] and
$\cali:=\{1,\dots,K\}\times\{1,\dots,L\}$.

Notice that in the last three definitions of \eqref{eq:abbr}, the contribution of the initial condition to the increments is neglected. Also, observe that $u(t,x)_b^c$ and $\Delta^n_i u^c_b$ are $\calf_c$-measurable and that $u(t,x)_{b,a}^c$, $\wt u(t,x)_{b,a}^c$, $\Delta^n_i u_{b,a}^c$, and $\Delta^n_i \wt u_{b,a}^c$ are $\calf_a$-conditionally Gaussian.

\brem\label{rem:2} 
The only place where the different coordinates of $f$ interact is in the CLT in Lemma~\ref{lem:CLT-core} below. All other steps can be carried out coordinate by coordinate, so for them, we may and will assume $M=1$ in the sequel. Moreover, as a consequence of hypothesis H2 and the third part of Remark~\ref{rem:1}, if we assume $M=1$, we could also assume $K=1$, which would simplify the notation a bit. But we will mostly refrain from doing so because in the second part \cite{Chong19c} of this paper, H2 will no longer be (always) in place. Thus, keeping $K$ general will spare us repeating many proofs in \cite{Chong19c}.
\erem

\subsection{Martingale approximations and the core CLT}\label{sect:steps12}

In order to avoid notational difficulties, we first remove terms from \eqref{eq:Vnf} that correspond to small values of $i$.
\blem\label{lem:remove-lan}
If $\la_n:=[\Del^{-a}]$ where $a\in(0,1)$ is larger than but close to $\frac{1}{2\Theta}$ and $\Theta:=1+\frac\al2$, then, as $n\to\infty$,
\beq\label{eq:remove-lan} 
 \frac{1}{\Delh}\Biggl(V^n_f(u,t)-\Del\sumtla  f\biggl( \frac{\Delta^n_i u}{\tau_n},\dots,\frac{\Delta^n_{i+L-1} u}{\tau_n}\biggr)\Biggr) \limL 0.\eeq
\elem

Next, we cut off the stochastic integrals defining $\Delta^n_i u$ at a reasonable temporal distance from $i\Del$. This procedure decreases the overlap and hence improves the conditional independence between different increments. The same step already appeared in \cite{Chong19} (see also \cite{BN11, Corcuera13}), but the proofs there were much simpler due to the  $\frac12$-Hölder regularity in time of the noise coefficient.
\blem\label{lem:trunc} 
For any $a>\frac{1}{2\Theta}$, we have, as $n\to\infty$,
\beq\label{eq:approx-ga} \Delh\sumtla \Biggl\{ f\biggl( \frac{\Delta^n_i u}{\tau_n},\dots,\frac{\Delta^n_{i+L-1} u}{\tau_n}\biggr) - f\Biggl( \frac{\Delta^n_i u_{(i-\la_n)\Del}^{(i+L-1)\Del}}{\tau_n}\Biggr) \Biggr\} \limL 0.\eeq
\elem

In the case of additive noise (i.e., $\si$ is a constant function), if we center the resulting terms by their respective conditional expectation, then 
\beq\label{eq:notconv} \ov V^n_{\la_n}(u,t):=\Delh\sumtla \Biggl\{ f\Biggl( \frac{\Delta^n_i u_{(i-\la_n)\Del}^{(i+L-1)\Del}}{\tau_n}\Biggr) -\E\Biggl[f\Biggl( \frac{\Delta^n_i u_{(i-\la_n)\Del}^{(i+L-1)\Del}}{\tau_n}\Biggr) \mathrel{\bigg|} \calf^n_{i-\la_n}  \Biggr] \Biggr\} \limst \calz, \eeq
where $\calz$ is the process described in Theorem~\ref{thm:CLT} and $\calf^n_i := \calf_{i\Del}$; see \cite[Equation~(3.15)]{Chong19}. This, however, is no longer true under multiplicative noise. We will now identify the term that prevents $\ov V^n_{\la_n}(u,t)$ from converging and show that the remainder does converge stably in law to $\calz$. 

To this end, let $\ov a \in(\frac{a}{\Theta},a)$, $\ov\la_n :=[\Del^{-\ov a}]$, and
\beq\label{eq:betaga}\begin{split} \beta^n_i &:= \frac{\Delta^n_i  u_{(i-\la_n)\Del}^{(i+L-1)\Del}}{\tau_n},\qquad
	\beta_{i,k}^n:=\frac{\Delta^n_i \wt u_{(i-\la_n)\Del,(i-k-\la_n)\Del}^{(i+L-1)\Del}}{\tau_n},\\
	\ov\beta^n_i &:= \frac{\Delta^n_i u_{(i-\ov\la_n)\Del}^{(i+L-1)\Del}}{\tau_n},\qquad
	\ov\beta_{i,k}^n:=\frac{\Delta^n_i \wt u_{(i-\ov\la_n)\Del,(i-k-\la_n)\Del}^{(i+L-1)\Del}}{\tau_n},\\
	\ga^n_{i,k}&:=f(\beta^n_{i,k})-\E\Bigl[f(\beta^n_{i,k})\mathrel{\big|} \calf^n_{i-\la_n} \Bigr].
\end{split}\eeq
For an arbitrary number $m\in\N$, we can then consider the decomposition 
\beq\label{eq:BC}\ov V^n_{\la_n}(u,t)=B^n_1(t)+B^{n,m}_2(t)+B^{n,m}_3(t)+C^{n,m}(t),\eeq 
where
\beq\label{eq:B}\begin{split}
B^n_1(t)&:=\Delh\sumtla \Bigl\{ f( \beta^n_i) -f( \ov\beta^n_i)-\E\Bigl[f(\beta^n_i) -f(\ov\beta^n_i) \mathrel{\big|} \calf^n_{i-\la_n}  \Bigr]\Bigr\},\\
 B^{n,m}_2(t)&:=\Delh\sumtla \Bigl\{ f(\ov\beta^n_i) -f(\ov\beta^n_{i,k^{n,m}_i})-\E\Bigl[ f(\ov\beta^n_i) -f(\ov\beta^n_{i,k^{n,m}_i}) \mathrel{\big|} \calf^n_{i-\la_n} \Bigr]\Bigr\},\\
  B^{n,m}_3(t)&:=\Delh\sumtla \Bigl\{ f(\ov\beta^n_{i,k^{n,m}_i}) -f(\beta^n_{i,k^{n,m}_i})-\E\Bigl[ f(\ov\beta^n_{i,k^{n,m}_i}) -f(\beta^n_{i,k^{n,m}_i}) \mathrel{\big|} \calf^n_{i-\la_n} \Bigr]\Bigr\},\\
  C^{n,m}(t)&:=\Delh\sumtla \ga^n_{i,k^{n,m}_i},
\end{split}\eeq
and $k^{n,m}_i$ is defined as follows: let $\ov k^{n,m}_i$ be equal to the remainder of dividing $i$ by $(m+1)\la_n+L-1$, except that the outcome $0$ is replaced by $(m+1)\la_n+L-1$ [e.g., if $m=2$, $\la_n=4$, and $L=1$, we have $(m+1)\la_n+L-1=12$, so $\ov k^{n,m}_i=\ov k^{n,m}_{12+i} = \ov k^{n,m}_{24+i} = \dots = i$ for $i=1,\dots,12$]. Now if $\ov k^{n,m}_i \in \{1,\dots,m\la_n\}$, we define $k^{n,m}_i:=\ov k^{n,m}_i$, while if $\ov k^{n,m}_i \in \{m\la_n+1,\dots,(m+1)\la_n+L-1\}$, we define $k^{n,m}_i:=\ov k^{n,m}_i-m\la_n$ (e.g., in the situation just considered, $k^{n,m}_i=k^{n,m}_{12+i}=k^{n,m}_{24+i} = \dots = i$ for $i=1,\dots,8$ and $k^{n,m}_{8+i}=k^{n,m}_{20+i}=k^{n,m}_{32+i} = \dots = i$ for $i=1,\dots,4$). This complicated definition is owed to the classical block splitting technique used in the proof of \cite[Theorem~2.3]{Chong19}, from which we can now deduce the following result:
\blem\label{lem:CLT-core} If $a>\frac{1}{2\Theta}$ is sufficiently small, then there is a decomposition $C^{n,m}=C^{n,m}_1+C^{n,m}_2$ such that for every $m\in\N$,
$$C^{n,m}_1\limst \sqrt{\frac{m}{m+1}}\calz$$ 
as $n\to\infty$ and such that for every $T>0$,
\beq\label{eq:C2} \lim_{m\to\infty} \limsup_{n\to\infty} \E\Bigl[(C^{n,m}_2)^\ast_T\Bigr] =0, \eeq
where  $X^\ast_T :=\sup_{t\in[0,T]} |X(t)|$ for a stochastic process $(X(t))_{t\geq0}$.
\elem

Apart from this lemma, the reader does not need to remember the detailed definition of $k^{n,m}_i$. In fact, the only information to keep in mind is that for fixed $m$, the shift $k^{n,m}_i$ is at most $O(\la_n)$. As a consequence, if a $\beta$-term carries the subscript $k^{n,m}_i$, then the time parameter $s$ of $\si(u(s,y))$ in the definition of the $u$-increment is frozen at the point $(i-\la_n-k^{n,m}_i)\Del$, which, again for $m$ fixed, is always $O(\la_n\Del)$ from the point $i\Del$ apart. 

Two $B$-terms can be handled with ease. 
\blem\label{lem:B13} We have $B^n_1\limL0$ and, for every $m\in\N$, $B^{n,m}_3\limL0$ as $n\to\infty$.
\elem

Lemmas~\ref{lem:remove-lan}--\ref{lem:B13} together with \eqref{eq:BC} show that for every $m\in\N$,
\beq\label{eq:step1} \begin{split}& \frac{1}{\Delh}\Biggl(V^n_f(u,t)-\Del\sumtla  \E\Biggl[f\Biggl( \frac{\Delta^n_i u_{(i-\la_n)\Del}^{(i+L-1)\Del}}{\tau_n}\Biggr) \mathrel{\bigg|} \calf^n_{i-\la_n}  \Biggr]\Biggr) -B^{n,m}_2(t)-C^{n,m}_2(t)\\ &\quad\limst \sqrt{\frac{m}{m+1}}\calz \end{split} \eeq
as $n\to\infty$.
The term $C^{n,m}_2$ is harmless as it vanishes when first $n$ and then $m$ tend to infinity. This, however, is not true for $B^{n,m}_2$: it is the first term that we encounter and is \emph{not} asymptotically negligible. We postpone its analysis to Section~\ref{sect:bad}. 

\subsection{Approximating the conditional expectation}\label{sect:step3}
We move on to investigating the conditional expectation in \eqref{eq:step1}. 
By Taylor's formula, we can develop this into
$\sum_{j=1}^4 \delta^{n,i}_j$,
where
\begin{align*}
\delta^{n,i}_1
 &:= \E\Biggl[f\Biggl( \frac{ \Delta^n_i  u_{(i-\la_n)\Del,(i-\la_n)\Del}^{(i+L-1)\Del}}{\tau_n}\Biggr) \mathrel{\bigg|} \calf^n_{i-\la_n}  \Biggr],\\
\delta^{n,i}_2&:= \sum_{\nu\in\cali}\E\Biggl[\partial_\nu f\Biggl( \frac{\Delta^n_i  u_{(i-\la_n)\Del,(i-\la_n)\Del}^{(i+L-1)\Del}}{\tau_n}\Biggr) \frac{(\Delta^n_i u_{(i-\la_n)\Del}^{(i+L-1)\Del}-\Delta^n_i  u_{(i-\la_n)\Del,(i-\la_n)\Del}^{(i+L-1)\Del})_\nu}{\tau_n}\mathrel{\bigg|} \calf^n_{i-\la_n}  \Biggr],\\
\delta^{n,i}_3&:=\frac12\sum_{\nu_1,\nu_2\in\cali} \E\Biggl[\partial^2_{\nu_1\nu_2} f\Biggl( \frac{\Delta^n_i  u_{(i-\la_n)\Del,(i-\la_n)\Del}^{(i+L-1)\Del}}{\tau_n}\Biggr) \frac{(\Delta^n_i u_{(i-\la_n)\Del}^{(i+L-1)\Del}-\Delta^n_i  u_{(i-\la_n)\Del,(i-\la_n)\Del}^{(i+L-1)\Del})_{\nu_1}}{\tau_n}\\
&\mathrel{\hphantom{:=} \mathrel{\times}}\frac{(\Delta^n_i u_{(i-\la_n)\Del}^{(i+L-1)\Del}-\Delta^n_i  u_{(i-\la_n)\Del,(i-\la_n)\Del}^{(i+L-1)\Del})_{\nu_2}}{\tau_n}\mathrel{\bigg|} \calf^n_{i-\la_n}  \Biggr],\\
\delta^{n,i}_4&:=\frac16\sum_{\nu_1,\nu_2,\nu_3\in\cali} \E\Biggl[\partial^3_{\nu_1\nu_2\nu_3} f(\eps^n_i) \frac{(\Delta^n_i u_{(i-\la_n)\Del}^{(i+L-1)\Del}-\Delta^n_i  u_{(i-\la_n)\Del,(i-\la_n)\Del}^{(i+L-1)\Del})_{\nu_1}}{\tau_n}\\
&\mathrel{\hphantom{:=} \mathrel{\times}}\frac{(\Delta^n_i u_{(i-\la_n)\Del}^{(i+L-1)\Del}-\Delta^n_i  u_{(i-\la_n)\Del,(i-\la_n)\Del}^{(i+L-1)\Del})_{\nu_2}}{\tau_n}\\
&\mathrel{\hphantom{:=} \mathrel{\times}}\frac{(\Delta^n_i u_{(i-\la_n)\Del}^{(i+L-1)\Del}-\Delta^n_i  u_{(i-\la_n)\Del,(i-\la_n)\Del}^{(i+L-1)\Del})_{\nu_3}}{\tau_n}\mathrel{\bigg|} \calf^n_{i-\la_n}  \Biggr]
\end{align*}
for some intermediate value $\eps^{n,i}_1$.  Accordingly, we have the decomposition
\beq\label{eq:D}
\Delh\sumtla \E\Biggl[f\Biggl( \frac{\Delta^n_i u_{(i-\la_n)\Del}^{(i+L-1)\Del}}{\tau_n}\Biggr) \mathrel{\bigg|} \calf^n_{i-\la_n}  \Biggr] = \sum_{j=1}^4 D^n_j(t),\qquad D^n_j(t):=\Delh\sumtla \delta^{n,i}_j.
\eeq
\blem\label{lem:D45} If $a>\frac1{2\Theta}$ is sufficiently small, then $D^n_4\limL0$ 
as $n\to\infty$.
\elem

Unfortunately, $D^n_2$ and $D^n_3$ are two further terms that do not vanish asymptotically; they will be further investigated in Section~\ref{sect:bad}. For the analysis of $D^n_1$, we have to introduce some more notations. Given $r\in\N_0$, $h\in\R^d$, and $n\in\N$, we consider measures on $[0,\infty)\times\R^d\times\R^d$ defined via
\begin{align*} \Pi^n_{r,h}(A)
	&:=\iiint_A\frac{ G_y(s)-G_y(s-\Delta_n)}{\tau_n} \frac{G_{y'+h}(s+r\del)-G_{y'+h}(s+(r-1)\del)}{\tau_n}\,\La(\dd y,\dd y')\,\dd s,\\
	|\Pi^n_{r,h}|(A)
	&:=\iiint_A\frac{ |G_y(s)-G_y(s-\Delta_n)||G_{y'+h}(s+r\del)-G_{y'+h}(s+(r-1)\del)|}{\tau^2_n}\,\La(\dd y,\dd y') \,\dd s
\end{align*}
for $A\in \calb([0,\infty)\times\R^d\times \R^d)$ [recall the convention $G_y(s):=0$ if $s<0$]. Moreover, let $\un\mu_f$ be the function that maps $(v_{\nu_1,\nu_2})_{\nu_1,\nu_2\in\cali} \in\R^{\cali\times\cali}$ to $\bbe[f(\un Z)]$ where $\un Z=(\un Z_\nu)_{\nu\in\cali}$ has a multivariate normal distribution with mean $0$ and $\cov(\un Z_{\nu_1},\un Z_{\nu_2})=v_{\nu_1,\nu_2}$. In particular, if $v(w)_{k_1l_1,k_2l_2}$ denotes the right-hand side of \eqref{eq:cov-Z} for some $w=(w_1,\dots,w_K)\in[0,\infty)^K$, then 
\beq\label{eq:muid}\un\mu_f(v(w))=\mu_f(w)\eeq
for the mapping $\mu_f$ defined in Theorem~\ref{thm:LLN}.

Since the argument of $f$ in the definition of $\delta^{n,i}_1$ is $\calf^n_{i-\la_n}$-conditionally Gaussian, we can rewrite
\beq\label{eq:Dn1} D^n_1(t)=\Delh\sumtla \un\mu_f\Bigl((v^{n,i}_{\nu_1,\nu_2})_{\nu_1,\nu_2\in\cali}\Bigr), \eeq
with 
\beq\label{eq:v}\begin{split}
v^{n,i}_{\nu_1,\nu_2}&:=\iiint_{(i-\la_n)\Del}^{(i+l_{12}-1)\Del} \frac{\Delta^n_i G_y(s)_{\nu_1}\Delta^n_i G_{y'}(s)_{\nu_2}}{\tau_n^2}\si(u(s,y)_0^{(i-\la_n)\Del})\si(u(s,y')_0^{(i-\la_n)\Del}) \,\La(\dd y,\dd y')\,\dd s\\
&\hphantom{:}= \iiint_0^{(\la_n+l_{12}-1)\Del} \si(u(t^{n,i}_{l_{12}-1}-s,x_{k_{12}}-y)_0^{(i-\la_n)\Del})\\
&\mathrel{\hphantom{:=} \mathrel{\times}}\si(u(t^{n,i}_{l_{12}-1}-s,x_{k_{12}}-y')_0^{(i-\la_n)\Del})\,\Pi^n_{|l_1-l_2|,x_{k_1}-x_{k_2}}(\dd s,\dd y,\dd y'),
\end{split}\raisetag{2\baselineskip}\eeq
where $t^{n,i}_j:=(i+j)\Del$,  $l_{12}:=l_1\wedge l_2$, and $k_{12}:=k_1$ if $l_{12}=l_1$ and $k_{12}:=k_2$ otherwise. 

Let us pause for a moment and take a careful look at \eqref{eq:v}: if we could simply replace the last two $\si$-terms by $\si(u(i\Del,x))$ and the remaining $\Pi^n$-integral by $\Ga_{|l_1-l_2|}\bone_{k_1=k_2}$, this would lead to an expression in \eqref{eq:Dn1} that is a Riemann sum approximation of $V_f(u,t)$, the limit in \eqref{eq:LLN} that we want to get to. The goal is therefore to try to justify these approximations; and this is exactly the path that was taken in \cite{Chong19} as well. But the proofs there require nice Hölder and differentiability properties of $\si(u(s,y))$, which we simply do not have in the multiplicative case. Thus, it is not immediately clear whether we will succeed or not.

There is, however, a small but important difference between \cite{Chong19} and our present situation that helps us out: we do not have $\si(u(t^{n,i}_{l_{12}-1}-s,x_{k_{12}}-y^{(\prime)}))$ but $\si(u(t^{n,i}_{l_{12}-1}-s,x_{k_{12}}-y^{(\prime)})_0^{(i-\la_n)\Del})$ in \eqref{eq:v}. Since the main weight of the $\Pi^n$-measures is on the origin [in fact, for any $\eps>0$, we have $|\Pi^n_{r,h}|([\eps,\infty)\times\R^d\times\R^d)\to0$ uniformly in $r$ and $h$ by \eqref{eq:Theta}] and since, at the same time, the regularity of $(s,y^{(\prime)})\mapsto\si(u(t^{n,i}_{l_{12}-1}-s,x_{k_{12}}-y^{(\prime)})_0^{(i-\la_n)\Del})$ improves when $s$ is small (because then $t^{n,i}_{l_{12}-1}-s\approx i\Del$, which is ``far'' from $(i-\la_n)\Del$, so we skip the most singular part of the heat kernel in the mild formulation \eqref{eq:mild} of $u$), we \emph{gain} regularity compared to the situation when we simply had $\si(u(t^{n,i}_{l_{12}-1}-s,x_{k_{12}}-y^{(\prime)}))$. 
Of course, this comes at a price: first, we have already produced two bad terms $D^n_2$ and $D^n_3$ to get this truncation $(\cdot)_0^{(i-\la_n)\Del}$. And second, unsurprisingly, we will not be able to turn $D^n_1$ into $V_f(u,t)$ but only into the right-hand side of \eqref{eq:LLN} with $u(s,x_k)$ replaced by $u(s,x_k)_0^{s-\la_n\Del}$. This will lead to yet another term that does not vanish on a $\Delh$-rate. The remaining part of this subsection specifies the details of this discussion. 

We define $v^{\prime n,i}_{\nu_1,\nu_2}$ to be the last integral in \eqref{eq:v} but with $t^{n,i}_{l_{12}-1}-s$  replaced by $i\Del$ and the domain of integration replaced by $[0,\la_n\Del]\times\R^d\times\R^d$. Moreover, we define $v^{\prime\prime n,i}_{\nu_1,\nu_2}$ in the same way but with $y$ and $y'$ further replaced by $0$ (in addition to the two previous modifications). 

\blem\label{lem:remove-s}
If $a>\frac1{2\Theta}$, we have, as $n\to\infty$,
\beq\label{eq:help12} D^n_1(t) - \Delh\sumtla \un\mu_f\Bigl(( v^{\prime n,i}_{\nu_1,\nu_2})_{\nu_1,\nu_2\in\cali}\Bigr)  \limL 0. \eeq
\elem

\blem\label{lem:remove-y}
If $a>\frac1{2\Theta}$, we have, as $n\to\infty$,
\[ \Delh\sumtla \biggl(\un\mu_f\Bigl(( v^{\prime n,i}_{\nu_1,\nu_2})_{\nu_1,\nu_2\in\cali}\Bigr)-\un\mu_f\Bigl(( v^{\prime\prime n,i}_{\nu_1,\nu_2})_{\nu_1,\nu_2\in\cali}\Bigr)\biggr)  \limL 0. \]
\elem

Notice that by definition,
\beq\label{eq:vpp}\begin{split} v^{\prime\prime n,i}_{\nu_1,\nu_2} 
&= \iiint_0^{\la_n\Del} \si^2(u(i\Del,x_{k_{12}})_0^{(i-\la_n)\Del}) \,\Pi^n_{|l_1-l_2|,x_{k_1}-x_{k_2}}(\dd s,\dd y,\dd y')\\
&=\si^2(u(i\Del,x_{k_{12}})_0^{(i-\la_n)\Del})\Pi^n_{|l_1-l_2|,x_{k_1}-x_{k_2}}([0,\la_n\Del]\times\R^d\times\R^d).
 \end{split}\eeq

\blem\label{lem:remove-Pi} If $a>\frac1{2\Theta}$, then, as $n\to\infty$,
\begin{align*} &\sqrt{\Del}\sumtla \biggl( \un\mu_f\Bigl((v^{\prime\prime n,i}_{\nu_1,\nu_2})_{\nu_1,\nu_2\in\cali}\Bigr)-\mu_f\Bigl(\si^2(u(i\Del,\un x)_0^{(i-\la_n)\Del})\Bigr) \biggr)\limL0, \end{align*}
where $\si^2(u(t,\un x)_b^c):=(\si^2(u(t,x_1)_b^c),\dots,\si^2(u(t,x_K)_b^c))$ [and similarly without the $(\cdot)_b^c$-part].
\elem

\blem\label{lem:iDels}  If $a>\frac1{2\Theta}$ is small, then, as $n\to\infty$,
\beq \label{eq:help2}\frac{1}{\sqrt{\Del}} \Biggl( \Del\sumtla \mu_f\Bigl(\si^2(u(i\Del,\un x)_0^{(i-\la_n)\Del})\Bigr)- \int_{\la_n\Del}^t \mu_f\Bigl(\si^2(u(s,\un x)_0^{s-\la_n\Del})\Bigr)\,\dd s  \Biggr)\limL 0. \eeq
\elem

The last remaining portion 
\beq\label{eq:last}\begin{split}
	H^n(t):=\frac{1}{\sqrt{\Del}}\Biggl(\int_{\la_n\Del}^t \mu_f\Bigl(\si^2(u(s,\un x)_0^{s-\la_n\Del})\Bigr)\,\dd s-V_f(u,t)\Biggr)
	\end{split}\eeq
 is another term that does not vanish as $n\to\infty$.
The following proposition summarizes what we have shown so far in Sections~\ref{sect:steps12} and \ref{sect:step3} [combine \eqref{eq:step1}, \eqref{eq:D}, and \eqref{eq:last} with Lemmas~\ref{lem:D45}--\ref{lem:iDels}]:
\bprop\label{prop:1} Under the hypotheses of Theorem~\ref{thm:CLT}, if $a$ is close to but larger than $\frac{1}{2\Theta}$ and $\ov a\in (\frac{a}{\Theta},a)$, then for every $m\in\N$, as $n\to\infty$,
\beq\label{eq:CLT-prelim} \frac1{\Delh}\Bigl( V^n_f(u,t)-V_f(u,t)\Bigr) -\Bigl(C^{n,m}_2(t)+B^{n,m}_2(t)+D^n_2(t)+D^n_3(t)+H^n(t)\Bigr) \limst \sqrt{\frac{m}{m+1}}\calz. \eeq
\eprop

\subsection{The bad terms}\label{sect:bad}

In the previous two subsections, we identified four terms that were ill-behaved: $B^{n,m}_2$, $D^n_2$, $D^n_3$, and $H^n$. Our goal is now to approximate them by simpler (but still nonvanishing) terms. This will later allow us to exploit cancellation effects among them. 

\subsubsection*{The term $\boldsymbol{B^{n,m}_2}$}

Looking back at the definition in \eqref{eq:B}, we realize that $B^{n,m}_2$ has a very special property. For every $i$, the expression in braces has a vanishing $\calf^n_{i-\la_n}$-conditional expectation. This gives the whole sum a certain kind of martingale structure, which is helpful for many estimates (to be more precise, $B^{n,m}_2$ is a typical term where \emph{martingale size estimates} as explained in Section~\ref{sect:details} apply). In what follows, we will ensure that all approximations of $B^{n,m}_2$ retain this important feature.
\blem\label{lem:Bn1}
Define
\beq\label{eq:Bn1} B^{n,m,1}_2(t):=\Delh\sumtla \sum_{\nu\in\cali} \Bigl\{(\theta^{n,m,i}_1)_\nu-\E[(\theta^{n,m,i}_1)_\nu \mid \calf^n_{i-\la_n}]\Bigr\}, \eeq
where
\beq\label{eq:theta1} (\theta^{n,m,i}_1)_\nu:=\partial_\nu f(\ov\beta^n_{i,k^{n,m}_i}) \Bigl( \ov\beta^n_i-\ov\beta^n_{i,k^{n,m}_i} \Bigr)_\nu. 
 \eeq
Then $B^{n,m}_2-B^{n,m,1}_2\limL0$ for every $m\in\N$ as $n\to\infty$ if $a$ is sufficiently close to $\frac{1}{2\Theta}$.
\elem

Writing out the $\beta$-difference in \eqref{eq:theta1} as a stochastic integral and using Taylor's theorem on the resulting $\si$-difference, we split $\theta^{n,m,i}_1$ into four parts $\theta^{n,m,i}_{1,1}+\theta^{n,m,i}_{1,2}+\theta^{n,m,i}_{1,3}+\theta^{n,m,i}_{1,4}$ given by 
\begin{align*}
(\theta^{n,m,i}_{1,1})_\nu &:= \partial_\nu f(\ov\beta^n_{i,k^{n,m}_i}) \iint_{(i-\ov\la_n)\Del}^{(i+l-1)\Del} \frac{\Delta^n_i G_y(s)_\nu}{\tau_n}\si'(u((i-k^{n,m}_i-\la_n)\Del,y))\\
&\qtimes\iint_0^{(i-\la_n)\Del} \Bigl(G_{y-z}(s-r)-G_{y-z}((i-k^{n,m}_i-\la_n)\Del-r)\Bigr) \\
&\cquad{\mathrel{\times}\si(u(r,z))\,W(\dd r,\dd z)\,W(\dd s,\dd y),}\\
(\theta^{n,m,i}_{1,2})_\nu &:= \partial_\nu f(\ov\beta^n_{i,k^{n,m}_i}) \iint_{(i-\ov\la_n)\Del}^{(i+l-1)\Del} \frac{\Delta^n_i G_y(s)_\nu}{\tau_n}\si'(u((i-k^{n,m}_i-\la_n)\Del,y))\\
&\qtimes\iint_{(i-\la_n)\Del}^{(i-\ov\la_n)\Del}  G_{y-z}(s-r)  \si(u(r,z))\,W(\dd r,\dd z)\,W(\dd s,\dd y),\\
(\theta^{n,m,i}_{1,3})_\nu &:= \partial_\nu f(\ov\beta^n_{i,k^{n,m}_i}) \iint_{(i-\ov\la_n)\Del}^{(i+l-1)\Del} \frac{\Delta^n_i G_y(s)_\nu}{\tau_n}\si'(u((i-k^{n,m}_i-\la_n)\Del,y))\\
&\qtimes\iint_{(i-\ov\la_n)\Del}^s  G_{y-z}(s-r)  \si(u(r,z))\,W(\dd r,\dd z)\,W(\dd s,\dd y),\\
(\theta^{n,m,i}_{1,4})_\nu &:= \frac12 \partial_\nu f(\ov\beta^n_{i,k^{n,m}_i}) \iint_{(i-\ov\la_n)\Del}^{(i+l-1)\Del} \frac{\Delta^n_i G_y(s)_\nu}{\tau_n}\si''(\xi^n_i(s,y))\Biggl(\iint_0^s  \Bigl(G_{y-z}(s-r)\\
&\qminus G_{y-z}((i-k^{n,m}_i-\la_n)\Del-r)\Bigr) \si(u(r,z))\,W(\dd r,\dd z)\Biggr)^2\,W(\dd s,\dd y),
\end{align*}
where $\xi^n_i(s,y)$ is some intermediate value.
 By \eqref{eq:Bn1}, we obtain $B^{n,m,1}_2=\sum_{j=1}^4 B^{n,m,1,j}_2$ where
 \beq\label{eq:Bnm1j} B^{n,m,1,j}_2(t):=\Delh\sumtla \sum_{\nu\in\cali} \Bigl\{(\theta^{n,m,i}_{1,j})_\nu-\E[(\theta^{n,m,i}_{1,j})_\nu \mid \calf^n_{i-\la_n}]\Bigr\}. \eeq
Only $B^{n,m,1,2}_2$ is asymptotically relevant. Furthermore, if $(\theta^{n,i}_2)_\nu$ [resp., $B^{n,2}_2(t)$] denotes the term resulting from replacing  $k^{n,m}_i$ [resp., $(\theta^{n,m,i}_{1,2})_\nu$] in the definition of $(\theta^{n,m,i}_{1,2})_\nu$ [resp., $B^{n,m,1,2}_2(t)$] by $0$ [resp., $(\theta^{n,i}_2)_\nu$], then $B^{n,m,1,2}_2$ can be approximated by $B^{n,2}_2$, a term that is independent of $m$.
 \blem\label{lem:Bn2}
 If $a$ and $\ov a$ are sufficiently close to $\frac1{2\Theta}$ and $\frac a \Theta$, respectively, then $B^{n,m,1}_2-B^{n,2}_2\limL0$  as $n\to\infty$ for every $m\in\N$.
 \elem

Next, define
\[  B^{n,3}_2(t):=\Delh\sumtla \sum_{\nu\in\cali} (\theta^{n,i}_3)_\nu, \]
where $(\theta^{n,i}_3)_{\nu}$
is the conditional expectation of 
\beq\label{eq:wttheta3}\begin{split}
(\wt \theta^{n,i}_3)_{\nu}&:=
\partial_\nu f\Biggl(\frac{\Delta^n_i  u^{(i+L-1)\Del}_{(i-\ov\la_n)\Del,(i-\la_n)\Del}}{\tau_n}\Biggr) \iint_{(i-\ov\la_n)\Del}^{(i+l-1)\Del} \frac{\Delta^n_i G_y(s)_\nu}{\tau_n}\si'(u(s,y)_0^{(i-\la_n)\Del})\\
&\qtimes\iint_{(i-\la_n)\Del}^{(i-\ov\la_n)\Del}  G_{y-z}(s-r)  \si(u(r,z)_0^{(i-\la_n)\Del})\,W(\dd r,\dd z)\,W(\dd s,\dd y)
\end{split}\eeq
with respect to $\calf^n_{i-\ov\la_n}$ (n.b., \emph{not} with respect to $\calf^n_{i-\la_n}$).
Compared to $(\theta^{n,i}_2)_\nu$, we have replaced $u((i-\la_n)\Del,\cdot)$ by $u(\cdot,\cdot)_0^{(i-\la_n)\Del}$ (because the latter is easier to work with later), but more importantly, we have also replaced $u(r,z)$ by $u(r,z)_0^{(i-\la_n)\Del}$ such that the integrand of the double stochastic integral in \eqref{eq:wttheta3} is now $\calf^n_{i-\la_n}$-measurable. 
 \blem\label{lem:Bn3}
  If $a$ and $\ov a$ are sufficiently close to $\frac1{2\Theta}$ and $\frac a \Theta$, respectively, then $B^{n,2}_2-B^{n,3}_2\limL0$  as $n\to\infty$. 
 \elem
 
As in Lemmas~\ref{lem:remove-s} and \ref{lem:remove-y} (see also the discussion before them), we now make the $\si$- and $\si'$-terms independent of $s$ and $y$. To this end, let
\begin{align*}  B^{n,4|5}_2(t)&:=\Delh\sumtla \sum_{\nu\in\cali} (\theta^{n,i}_{4|5})_\nu,\qquad (\theta^{n,i}_{4|5})_{\nu}:=\E[(\wt \theta^{n,i}_{4|5})_{\nu}\mid \calf^n_{i-\ov\la_n}],
  \end{align*}
where [recall that $\nu=(k,l)$ by convention]
\begin{align*}
(\wt \theta^{n,i}_4)_{\nu}&:=
\partial_\nu f\Biggl(\iint_{(i-\ov\la_n)\Del}^{(i+L-1)\Del} \frac{\Delta^n_i G_y(s)}{\tau_n}\si(u(i\Del,y)_0^{(i-\la_n)\Del})\,W(\dd s,\dd y) \Biggr)\\
&\qtimes \iint_{(i-\ov\la_n)\Del}^{(i+l-1)\Del} \frac{\Delta^n_i G_y(s)_\nu}{\tau_n}\si'(u(i\Del,y)_0^{(i-\la_n)\Del})\\
&\cquad{\mathrel{\times}\iint_{(i-\la_n)\Del}^{(i-\ov\la_n)\Del}  G_{y-z}(i\Del-r)  \si(u(r,z)_0^{(i-\la_n)\Del})\,W(\dd r,\dd z)\,W(\dd s,\dd y),}\\
(\wt \theta^{n,i}_5)_{\nu}&:=
\partial_\nu f\Biggl(\biggl(\iint_{(i-\ov\la_n)\Del}^{(i+l-1)\Del} \frac{\Delta^n_i G_y(s)_\nu}{\tau_n}\si(u(i\Del,x_k)_0^{(i-\la_n)\Del})\,W(\dd s,\dd y)\biggr)_{\nu\in\cali} \Biggr)\\
&\qtimes \iint_{(i-\ov\la_n)\Del}^{(i+l-1)\Del} \frac{\Delta^n_i G_y(s)_\nu}{\tau_n}\si'(u(i\Del,x_k)_0^{(i-\la_n)\Del})\\
&\cquad{\mathrel{\times}\iint_{(i-\la_n)\Del}^{(i-\ov\la_n)\Del}  G_{x_k-z}(i\Del-r)  \si(u(r,z)_0^{(i-\la_n)\Del})\,W(\dd r,\dd z)\,W(\dd s,\dd y).}
\end{align*}

\blem\label{lem:Bn345} If $a$ and $\ov a$ are sufficiently close to $\frac1{2\Theta}$ and $\frac a \Theta$, respectively, then $B^{n,3}_2-B^{n,4}_2\limL0$ and $B^{n,4}_2-B^{n,5}_2 \limL0$ as $n\to\infty$. 
\elem

 For $v=(v_{\nu_1,\nu_2})_{\nu_1,\nu_2\in\cali}, w=(w_{\nu_1,\nu_2})_{\nu_1,\nu_2\in\cali}, c=(c_{\nu_1,\nu_2})_{\nu_1,\nu_2\in\cali}\in \R^{(K\times L)\times(K\times L)}$, and a function $f\colon \R^{K\times L}\to\R^M$, define 
\beq\label{eq:ovmu}\ov\mu_f(v,w,c)_\nu:=\E[\partial_\nu f(\ov Z^{(1)})\ov Z^{(2)}_\nu],\eeq 
where $\ov Z^{(1)}$ and $\ov Z^{(2)}$ are $\R^{K\times L}$-valued jointly Gaussian random elements with mean zero and $\E[\ov Z^{(1)}_{\nu_1}\ov  Z^{(1)}_{\nu_2}]=v_{\nu_1,\nu_2}$, $\E[\ov Z^{(2)}_{\nu_1} \ov Z^{(2)}_{\nu_2}]=w_{\nu_1,\nu_2}$, and $\E[\ov Z^{(1)}_{\nu_1}\ov  Z^{(2)}_{\nu_2}]=c_{\nu_1,\nu_2}$ [whenever the expectation is well-defined and $\ov Z^{(1)}$ and $\ov Z^{(2)}$ can be found with the desired properties]. Then
\begin{equation}\label{eq:thetamu}
(\theta^{n,i}_4)_\nu = \ov\mu_f(\ov v^{\prime n,i},\ov w^{\prime n,i},\ov c^{\prime n,i})_\nu,\qquad (\theta^{n,i}_5)_\nu = \ov\mu_f(\ov v^{\prime\prime n,i},\ov w^{\prime\prime n,i},\ov c^{\prime\prime n,i})_\nu,
\end{equation}
where [recall the notations introduced after \eqref{eq:v}]
\begin{align*}
\ov v^{\prime n,i}_{\nu_1,\nu_2} &:= 
\iiint \si(u(i\Del,x_{k_{12}}-y)_0^{(i-\la_n)\Del})\si(u(i\Del,x_{k_{12}}-y')_0^{(i-\la_n)\Del}) \\
&\qtimes \bone_{s\in(0,(\ov\la_n+l_{12}-1)\Del]} \,\Pi^n_{|l_1-l_2|,x_{k_1}-x_{k_2}}(\dd s,\dd y,\dd y'),\\
\ov w^{\prime n,i}_{\nu_1,\nu_2}&:=
\iiint \si'(u(i\Del,x_{k_{12}}-y)_0^{(i-\la_n)\Del})\si'(u(i\Del,x_{k_{12}}-y')_0^{(i-\la_n)\Del})\\
&\qtimes\iint_{(i-\la_n)\Del}^{(i-\ov\la_n)\Del}  G_{x_{k_{12}}-y-z}(i\Del-r)  \si(u(r,z)_0^{(i-\la_n)\Del})\,W(\dd r,\dd z)\\
&\qtimes\iint_{(i-\la_n)\Del}^{(i-\ov\la_n)\Del}  G_{x_{k_{12}}-y'-z}(i\Del-r)  \si(u(r,z)_0^{(i-\la_n)\Del})\,W(\dd r,\dd z)\\
&\qtimes \bone_{s\in(0,(\ov\la_n+l_{12}-1)\Del]}\,\Pi^n_{|l_1-l_2|,x_{k_1}-x_{k_2}}(\dd s,\dd y,\dd y'),\\
\ov c^{\prime n,i}_{\nu_1,\nu_2}&:= 
\iiint \si(u(i\Del,x_{k_{12}}-y)_0^{(i-\la_n)\Del})\si'(u(i\Del,x_{k_{12}}-y')_0^{(i-\la_n)\Del})\\
&\qtimes\iint_{(i-\la_n)\Del}^{(i-\ov\la_n)\Del}  G_{x_{k_{12}}-y'-z}(i\Del-r)  \si(u(r,z)_0^{(i-\la_n)\Del})\,W(\dd r,\dd z)\\
&\qtimes \bone_{s\in(0,(\ov\la_n+l_{12}-1)\Del]}\,\Pi^n_{|l_1-l_2|,x_{k_1}-x_{k_2}}(\dd s,\dd y,\dd y'),
\end{align*}
and $\ov v^{\prime\prime n,i}_{\nu_1,\nu_2}$, $\ov w^{\prime\prime n,i}_{\nu_1,\nu_2}$, and $\ov c^{\prime\prime n,i}_{\nu_1,\nu_2}$ are defined by the same formulas but with $y$ and $y'$ set to $0$. In fact, 
\begin{align*}
\ov v^{\prime\prime n,i}_{\nu_1,\nu_2} &= 
 \si^2(u(i\Del,x_{k_{12}})_0^{i\Del-\la_n\Del})  \Pi^n_{|l_1-l_2|,x_{k_1}-x_{k_2}}((0,(\ov\la_n+l_{12}-1)\Del]\times\R^d\times\R^d),\\
\ov w^{\prime\prime n,i}_{\nu_1,\nu_2}&= 
 \Biggl(\iint_{i\Del-\la_n\Del}^{i\Del-\ov\la_n\Del}  G_{x_{k_{12}}-z}(i\Del-r)  \si(u(r,z)_0^{i\Del-\la_n\Del})\,W(\dd r,\dd z)\Biggr)^2\\
 &\quad\times(\si')^2(u(i\Del,x_{k_{12}})_0^{i\Del-\la_n\Del})\Pi^n_{|l_1-l_2|,x_{k_1}-x_{k_2}}((0,(\ov\la_n+l_{12}-1)\Del]\times\R^d\times\R^d),\\
\ov c^{\prime\prime n,i}_{\nu_1,\nu_2}&= 
  \iint_{i\Del-\la_n\Del}^{i\Del-\ov\la_n\Del}  G_{x_{k_{12}}-z}(i\Del-r)  \si(u(r,z)_0^{i\Del-\la_n\Del})\,W(\dd r,\dd z)\\
 &\quad\times(\si\si')(u(i\Del,x_{k_{12}})_0^{i\Del-\la_n\Del})\Pi^n_{|l_1-l_2|,x_{k_1}-x_{k_2}}((0,(\ov\la_n+l_{12}-1)\Del]\times\R^d\times\R^d).
\end{align*}
Let $\ov v''_n(s)_{\nu_1,\nu_2}$, $\ov w''_n(s)_{\nu_1,\nu_2}$, and $\ov c''_n(s)_{\nu_1,\nu_2}$ [resp., $\ov v_n(s)_{\nu_1,\nu_2}$, $\ov w_n(s)_{\nu_1,\nu_2}$, and $\ov c_n(s)_{\nu_1,\nu_2}$] be defined through the last three expressions above but with $i\Del$ replaced by $s$ and the $\Pi^n$-factors replaced by $\Ga_{|l_1-l_2|}\bone_{k_1=k_2}$ (resp., with these changes and further with $\ov\la_n$ replaced by $0$).

\blem\label{lem:easy} If $\ov a>\frac a \Theta$, then $B^{n,5}_2-B^{n,6}_2\limL0$ and $B^{n,6}_2-B^{n,7}_2\limL0$ as $n\to\infty$, where
\begin{align*} B^{n,6}_2(t)&:=\sqrt{\Del}\sumtla \sum_{\nu\in\cali} (\theta^{n,i}_6)_\nu,\qquad  (\theta^{n,i}_6)_\nu:=\ov\mu_f(\ov v''_n(i\Del),\ov w''_n(i\Del), \ov c''_n(i\Del))_\nu, \\
B^{n,7}_2(t)&:=\sqrt{\Del}\sumtla \sum_{\nu\in\cali} (\theta^{n,i}_7)_\nu,\qquad  (\theta^{n,i}_7)_\nu:=\ov\mu_f(\ov v_n(i\Del),\ov w_n(i\Del), \ov c_n(i\Del))_\nu. 
\end{align*}
\elem

\blem\label{lem:B56} If $a$ is close enough to $\frac1{2\Theta}$, then $B^{n,7}_2-B^{n,8}_2\limL0$ as $n\to\infty$, where
\beq\label{eq:Bn8} B^{n,8}_2(t):=\frac1{\sqrt{\Del}}\int_{\la_n\Del}^t \sum_{\nu\in\cali} \ov\mu_f(\ov v_n(s),\ov w_n(s), \ov c_n(s))_\nu\,\dd s. \eeq
\elem

\subsubsection*{The terms $\boldsymbol{D^n_2}$ and $\boldsymbol{D^n_3}$}

Again, we want to simplify the terms $D^n_2$ and $D^n_3$ from \eqref{eq:D}. While for $B^{n,m}_2$, it was the martingale structure that helped us in our estimates, here we will use the fact that certain terms arising in the approximations below are removed because they have zero $\calf^n_{i-\la_n}$-conditional expectation. 
\blem\label{lem:Dn21}
If $a>\frac1{2\Theta}$ is small enough, we have $D^n_2-D^{n,1}_2\limL0$, where
\[
D^{n,1}_2(t):=\Delh\sumtla\sum_{\nu\in\cali} (\rho^{n,i}_1)_\nu,\qquad (\rho^{n,i}_1)_\nu:= \E[(\wt\rho^{n,i}_1)_\nu \mid\calf^n_{i-\la_n}],
\]
and
\begin{equation}\label{eq:wtrhophi-2}\begin{split} 
(\wt \rho^{n,i}_1)_\nu  &:= \frac12\partial_\nu f\Biggl( \frac{\Delta^n_i  u_{(i-\la_n)\Del,(i-\la_n)\Del}^{(i+L-1)\Del}}{\tau_n}\Biggr)\iint_{(i-\la_n)\Del}^{(i+l-1)\Del} \frac{\Delta^n_i G_y(s)_\nu}{\tau_n} \\
&\qtimes\si''(u(s,y)_0^{(i-\la_n)\Del})\Bigl(u(s,y)-u(s,y)_0^{(i-\la_n)\Del}\Bigr)^2 \,W(\dd s,\dd y). 
\end{split}
\end{equation}
\elem

We expand the square in $(\wt \rho^{n,i}_1)_\nu$ using the integration by parts formula and identify the asymptotically relevant part.
\blem\label{lem:Dn22} Define $D^{n,2}_2$ in the same way as $D^{n,1}_2$ but with 
\beq\label{eq:rhoni2}\begin{split}
(\wt \rho^{n,i}_2)_\nu  &:= \frac12\partial_\nu f\Biggl( \frac{\Delta^n_i  u_{(i-\la_n)\Del,(i-\la_n)\Del}^{(i+L-1)\Del}}{\tau_n}\Biggr)\iint_{(i-\la_n)\Del}^{(i+l-1)\Del} \frac{\Delta^n_i G_y(s)_\nu}{\tau_n} \si''(u(s,y)_0^{(i-\la_n)\Del})\\
&\qtimes\iiint_{(i-\la_n)\Del}^s G_{y-z}(s-r)G_{y-z'}(s-r)\\
&\cquad{\mathrel{\times}\si(u(r,z)_0^{(i-\la_n)\Del})\si(u(r,z')_0^{(i-\la_n)\Del})\,\La(\dd z,\dd z')\,\dd r  \,W(\dd s,\dd y)}
\end{split}\eeq
instead of $(\wt \rho^{n,i}_1)_\nu $. If $a>\frac{1}{2\Theta}$ is small enough, then $D^{n,1}_2-D^{n,2}_2\limL0$ as $n\to\infty$.
\elem

Recalling the notation \eqref{eq:ovmu}, we see that $\rho^{n,i}_2 := \E[\wt \rho^{n,i}_2 \mid \calf^n_{i-\la_n}]$ can be written as $\rho^{n,i}_2 = \frac12\ov\mu_f(\wt v^{n,i},\wt w^{n,i},\wt c^{n,i})$, where $\wt v^{n,i}:=v^{n,i}$ from \eqref{eq:v} and
\beq\label{eq:wc}\begin{split}
\wt w^{n,i}_{\nu_1,\nu_2}&:=\iiint \si''(u(t^{n,i}_{l_{12}-1}-s,x_{k_{12}}-y)_0^{(i-\la_n)\Del})\si''(u(t^{n,i}_{l_{12}-1}-s,x_{k_{12}}-y')_0^{(i-\la_n)\Del})\\
&\qtimes\iiint_0^{(\la_n+l_{12}-1)\Del-s} G_{z}(r)G_{z'}(r)\si(u(t^{n,i}_{l_{12}-1}-s-r,x_{k_{12}}-y-z)_0^{(i-\la_n)\Del})\\
&\cquad{\mathrel{\times}\si(u(t^{n,i}_{l_{12}-1}-s-r,x_{k_{12}}-y-z')_0^{(i-\la_n)\Del})\,\La(\dd z,\dd z')\,\dd r}\\
&\qtimes\iiint_0^{(\la_n+l_{12}-1)\Del-s} G_{z}(r)G_{z'}(r)\si(u(t^{n,i}_{l_{12}-1}-s-r,x_{k_{12}}-y'-z)_0^{(i-\la_n)\Del})\\
&\cquad{\mathrel{\times}\si(u(t^{n,i}_{l_{12}-1}-s-r,x_{k_{12}}-y'-z')_0^{(i-\la_n)\Del})\,\La(\dd z,\dd z')\,\dd r}\\
&\qtimes\bone_{s\in[0,(\la_n+l_{12}-1)\Del]}\,\Pi^n_{|l_1-l_2|,x_{k_1}-x_{k_2}}(\dd s,\dd y,\dd y'),\\
\wt c^{n,i}_{\nu_1,\nu_2}&:=\iiint \si(u(t^{n,i}_{l_{12}-1}-s,x_{k_{12}}-y)_0^{(i-\la_n)\Del})\si''(u(t^{n,i}_{l_{12}-1}-s,x_{k_{12}}-y')_0^{(i-\la_n)\Del})\\
&\qtimes\iiint_0^{(\la_n+l_{12}-1)\Del-s} G_{z}(r)G_{z'}(r)\si(u(t^{n,i}_{l_{12}-1}-s-r,x_{k_{12}}-y'-z)_0^{(i-\la_n)\Del})\\
&\cquad{\mathrel{\times}\si(u(t^{n,i}_{l_{12}-1}-s-r,x_{k_{12}}-y'-z')_0^{(i-\la_n)\Del})\,\La(\dd z,\dd z')\,\dd r}\\
&\qtimes\bone_{s\in[0,(\la_n+l_{12}-1)\Del]}\,\Pi^n_{|l_1-l_2|,x_{k_1}-x_{k_2}}(\dd s,\dd y,\dd y').
\end{split}\eeq
Let $\wt w^{\prime n,i}_{\nu_1,\nu_2}$ and $\wt c^{\prime n,i}_{\nu_1,\nu_2}$ (resp., $\wt w^{\prime\prime n,i}_{\nu_1,\nu_2}$ and $\wt c^{\prime\prime n,i}_{\nu_1,\nu_2}$) be the same integrals as in \eqref{eq:wc} but with $t^{n,i}_{l_{12}-1}-s$  replaced by $i\Del$ and the domain of integration [of both the $\Pi^n$- and the $\La(\dd z,\dd z')\,\dd r$-integrals] replaced by $[0,\la_n\Del]\times\R^d\times\R^d$ (resp., with the same changes and additionally $y$ and $y'$ replaced by $0$). Similarly, let $\wt v^{\prime n,i}:= v^{\prime n,i}$ and  $\wt v^{\prime\prime n,i}:= v^{\prime\prime n,i}$ (as introduced before Lemma~\ref{lem:remove-s}). Then define 
\[ D^{n,3|4}_2(t):=\Delh\sumtla\sum_{\nu\in\cali} (\rho^{n,i}_{3|4})_\nu,\]
where $\rho^{n,i}_3 = \frac12\ov\mu_f(\wt v^{\prime n,i},\wt w^{\prime n,i},\wt c^{\prime n,i})$ and $\rho^{n,i}_4 =\frac12 \ov\mu_f(\wt v^{\prime \prime n,i},\wt w^{\prime\prime n,i},\wt c^{\prime\prime n,i})$.
\blem\label{lem:Dn234} If $a>\frac{1}{2\Theta}$ is sufficiently small, then $D^{n,2}_2-D^{n,3}_2\limL0$ and $D^{n,3}_2-D^{n,4}_2\limL0$.
\elem

By definition, we have
\begin{align*}
\wt v^{\prime\prime n,i}_{\nu_1,\nu_2}&=\si^2(u(i\Del,x_{k_{12}})_0^{i\Del-\la_n\Del})\Pi^n_{|l_1-l_2|,x_{k_1}-x_{k_2}}([0,\la_n\Del]\times\R^d\times\R^d),\\
\wt w^{\prime\prime n,i}_{\nu_1,\nu_2}&=(\si'')^2(u(i\Del,x_{k_{12}})_0^{i\Del-\la_n\Del})  \Biggl(\iiint_0^{\la_n\Del} G_z(r)G_{z'}(r)\si(u(i\Del-r,x_{k_{12}}-z)_0^{i\Del-\la_n\Del})\\
&\quad\times\si(u(i\Del-r,x_{k_{12}}-z')_0^{i\Del-\la_n\Del})\,\La(\dd z,\dd z')\,\dd r\Biggr)^2 \Pi^n_{|l_1-l_2|,x_{k_1}-x_{k_2}}([0,\la_n\Del]\times\R^d\times\R^d),\\
\wt c^{\prime\prime n,i}_{\nu_1,\nu_2} &=(\si\si'')(u(i\Del,x_{k_{12}})_0^{i\Del-\la_n\Del})\Biggl(\iiint_0^{\la_n\Del} G_z(r)G_{z'}(r)\si(u(i\Del-r,x_{k_{12}}-z)_0^{i\Del-\la_n\Del})\\
&\quad\times\si(u(i\Del-r,x_{k_{12}}-z')_0^{i\Del-\la_n\Del})\,\La(\dd z,\dd z')\,\dd r\Biggr) \Pi^n_{|l_1-l_2|,x_{k_1}-x_{k_2}}([0,\la_n\Del]\times\R^d\times\R^d).
\end{align*}
Let $\wt v_n(s)_{\nu_1,\nu_2}$, $\wt w_n(s)_{\nu_1,\nu_2}$, and $\wt c_n(s)_{\nu_1,\nu_2}$ be defined through the  three equations of the previous display but with $i\Del$ replaced by $s$ and the $\Pi^n$-factors replaced by $\Ga_{|l_1-l_2|}\bone_{k_1=k_2}$.

\blem\label{lem:Dn2456} If $a$ is close enough to $\frac1{2\Theta}$, then $D^{n,4}_2-D^{n,5}_2\limL0$ and $D^{n,5}_2-D^{n,6}_2\limL0$ as $n\to\infty$, where
\begin{align*} D^{n,5}_2(t)&:=\sqrt{\Del}\sumtla \sum_{\nu\in\cali} (\rho^{n,i}_5)_\nu,\qquad  (\rho^{n,i}_5)_\nu:=\frac12\ov\mu_f(\wt v_n(i\Del),\wt w_n(i\Del), \wt c_n(i\Del))_\nu, \\
D^{n,6}_2(t)&:=\frac1{2\sqrt{\Del}}\int_{\la_n\Del}^t \sum_{\nu\in\cali} \ov\mu_f(\wt v_n(s),\wt w_n(s), \wt c_n(s))_\nu\,\dd s. 
\end{align*}
\elem

The next problematic term is $D^n_3$ and its analysis largely parallels that of $D^n_2$. Therefore, we only state the final form of this term and only sketch the intermediate steps in the proof.
\blem\label{lem:Dn3} Define
\beq\label{eq:Dn6} \begin{split} D^{n,6}_3(t)&:=\frac{1}{2\Delh} \int_{\la_n\Del}^t \sum_{\nu_1,\nu_2\in\cali\colon k_1=k_2=k} (\si')^2(u(s,x_k)_0^{s-\la_n\Del})\wh\mu_f(\wh v_n(s),\wh w,\wh c_n(s))_{\nu_1,\nu_2}\\
	&\qtimes \iiint_0^{\la_n\Del} G_z(r)G_{z'}(r) \si(u(s-r,x_k-z)_0^{s-\la_n\Del})\\
	&\cquad{\mathrel{\times}\si(u(s-r,x_k-z')_0^{s-\la_n\Del})\,\La(\dd z,\dd z')\,\dd r\,\dd s,}\end{split} \eeq
where $\wh\mu_f(v,w,c)_{\nu_1,\nu_2}:=\E[\partial^2_{\nu_1\nu_2} f(\ov Z^{(1)})\ov Z^{(2)}_{\nu_1}\ov Z^{(2)}_{\nu_2}]$ [with $\ov Z^{(1)}$ and $\ov Z^{(2)}$ as described after \eqref{eq:ovmu}] and
\begin{align*}
\wh v_n(s)_{\nu_1,\nu_2}&:=\si^2(u(s,x_k)_0^{s-\la_n\Del})\Ga_{|l_1-l_2|}\bone_{k_1=k_2=k},\\
\wh w_{\nu_1,\nu_2}&:=\Ga_{|l_1-l_2|}\bone_{k_1=k_2=k},\\
\wh c_n(s)_{\nu_1,\nu_2}&:=\si(u(s,x_k)_0^{s-\la_n\Del})\Ga_{|l_1-l_2|}\bone_{k_1=k_2=k}.
\end{align*}
 If $a>\frac{1}{2\Theta}$ is small enough, then $D^n_{3}-D^{n,6}_3\limL0$.
\elem

\subsection{Putting pieces together: the final step of the proof}

The following proposition shows how the bad terms encountered in the previous subsections cancel each other:
\bprop\label{prop:cancel} If $a>\frac{1}{2\Theta}$ is small enough, then, as $n\to\infty$,
\beq\label{eq:cancel} B^{n,8}_2 + D^{n,6}_2+D^{n,6}_3 + H^n \limL 0. \eeq
\eprop

From this we can readily deduce Theorem~\ref{thm:CLT}.
\bpr[Proof of Theorem~\ref{thm:CLT}] The idea is to write the left-hand side of \eqref{eq:CLT-statement} as
\begin{align*}
& \frac1{\Delh}\Bigl( V^n_f(u,t)-V_f(u,t)\Bigr) -\Bigl(C^{n,m}_2(t)+B^{n,m}_2(t)+D^n_2(t)+D^n_3(t)+H^n(t)\Bigr)\\
&\quad +C^{n,m}_2(t) + (B^{n,m}_2(t)-B^{n,8}_2(t)) + (D^n_2(t)-D^{n,6}_2(t))+(D^n_2(t)-D^{n,6}_3(t))\\
&\quad+(B^{n,8}_2(t) + D^{n,6}_2(t)+D^{n,6}_3(t) + H^n(t)).
\end{align*}
For fixed $m$, the first line converges stably in law to $(\frac {m}{m+1})^{1/2}\calz$ by \eqref{eq:CLT-prelim}. At the same time, we have $B^{n,m}_2-B^{n,8}_2\limL0$,  $D^n_2(t)-D^{n,6}_2(t)\limL0$, $D^n_3(t)-D^{n,6}_3(t)\limL0$, and $B^{n,8}_2 + D^{n,6}_2+D^{n,6}_3 + H^n\limL0$ by Lemmas~\ref{lem:Bn1}--\ref{lem:B56}, Lemmas~\ref{lem:Dn21}--\ref{lem:Dn2456}, Lemma~\ref{lem:Dn3}, and Proposition~\ref{prop:cancel}, respectively. Thus, \eqref{eq:CLT-statement} follows from these observations in conjunction with \eqref{eq:C2} and \cite[Proposition~2.2.4]{Jacod12}.
\epr

\section{Details of the proof of Theorem~\ref{thm:CLT}}\label{sect:details}

Two types of estimates are frequently used in the following to determine the asymptotic behavior of complicated expressions: \emph{standard size estimates} and \emph{martingale size estimates}. The idea behind the former, which were already used in \cite{Chong19}, is best explained by considering, for example, the term
\beq\label{u-expr} \begin{split} U_n(t) &:= \delh\sum_{i=\ell_n}^{t(n)} h(z^n_i) \iint_{(i-\ell_n)\Del}^{(i-\ell'_n)\Del} \frac{\Delta^n_i G_{y}(s)}{\tau_n} (\si(u(s,y))-\si(u(s,y)_0^{(i-\ell_n)\Del}))\,W(\dd s,\dd y).\end{split} \eeq
Here, $\ell_n := [\del^{-\ell}]$ and $\ell^{\prime}_n := [\del^{-\ell^{\prime}}]$ with some $0< \ell'<\ell<1$, $h$ is a function with at most polynomial growth, and $z^n_i$ are random variables such that for every $p\in[2,\infty)$,
\beq\label{size} \sup_{n\in \N} \sup_{i=1,\ldots,T(n)} \bbe[|z^n_i|^p]<\infty. \eeq
In most cases, $z^n_i$ is a normalized increment, possibly truncated or with modified $u$
[e.g., $z^n_i=\iint_0^{(i+L-1)\Del} \frac{\Delta^n_i G_{y}(s)}{\tau_n} \si(u(s,y))\,W(\dd s,\dd y)$ or $z^n_i=\iint_{(i-\ell_n)\Del}^{(i+L-1)\Del} \frac{\Delta^n_i G_{y}(s)}{\tau_n} \si(u(s,y)_0^{(i-\ell_n)\Del})\,W(\dd s,\dd y)$]. 

Using  Hölder's inequality with exponents $\frac{p}{p-1}$ and $p$, we obtain 
\begin{align*}
\bbe[(U_n)^\ast_T]& \leq \delh\sum_{i=\ell_n}^{T(n)} \bbe[|h(z^n_i)|^{\frac p {p-1}}]^{\frac{p-1}{p}} \\
&\quad\times\bbe\Biggl[ \Biggl|\iint_{(i-\ell_n)\Del}^{(i-\ell'_n)\Del} \frac{\Delta^n_i G_{y}(s)}{\tau_n} (\si(u(s,y))-\si(u(s,y)_0^{(i-\ell_n)\Del}))\,W(\dd s,\dd y)\Biggr|^p\Biggr]^{\frac1p}.
\end{align*}
By \eqref{size} and the growth assumptions on $h$, we have $\bbe[|h(z^n_i)|^{p/(p-1)}]\lec 1$, uniformly in $i$ and $n$. Thus, by the Burkholder--Davis--Gundy inequality, Minkowski's integral inequality, and the Cauchy--Schwarz inequality, we get
\begin{align*}
\bbe[(U_n)^\ast_T]&\lec \delh\sum_{i=\ell_n}^{T(n)} \bbe\Biggl[\Biggl( \iiint_{(i-\ell_n)\Del}^{(i-\ell'_n)\Del} \frac{|\Delta^n_i G_{y}(s)||\Delta^n_i G_{z}(s)|}{\tau_n^2} |\si(u(s,y))-\si(u(s,y)_0^{(i-\ell_n)\Del})|\\
&\quad\times |\si(u(s,y'))-\si(u(s,y')_0^{(i-\ell_n)\Del})|\,\La(\dd y,\dd y')\,\dd s\Biggr)^{\frac p2}\Biggr]^{\frac1p}\\
&\leq \delh\sum_{i=\ell_n}^{T(n)} \Biggl( \iiint_{(i-\ell_n)\Del}^{(i-\ell'_n)\Del} \frac{|\Delta^n_i G_{y}(s) \Delta^n_i G_{z}(s)|}{\tau_n^2} \bbe\Bigl[|\si(u(s,y))-\si(u(s,y)_0^{(i-\ell_n)\Del})|^{\frac p2}\\
&\quad\times |\si(u(s,y'))-\si(u(s,y')_0^{(i-\ell_n)\Del})|^{\frac p 2}\Bigr]^{\frac 2p}\,\La(\dd y,\dd y')\,\dd s\Biggr)^{\frac12}\\
&\leq \delh\sum_{i=\ell_n}^{T(n)} \Biggl( \iiint_{(i-\ell_n)\Del}^{(i-\ell'_n)\Del} \frac{|\Delta^n_i G_{y}(s) \Delta^n_i G_{z}(s)|}{\tau_n^2}\bbe[|\si(u(s,y))-\si(u(s,y)_0^{(i-\ell_n)\Del})|^p]^{\frac1p}\\
&\quad\times \bbe[|\si(u(s,y'))-\si(u(s,y')_0^{(i-\ell_n)\Del})|^p]^{\frac1p}\,\La(\dd y,\dd y')\,\dd s\Biggr)^{\frac12}.
\end{align*}
Because $\si$ is differentiable, there is some $\om^{n,i}(s,y)$ between $u(s,y)$ and $u(s,y)_0^{(i-\ell_n)\Del}$ such that $\si(u(s,y))-\si(u(s,y)_0^{(i-\ell_n)\Del}) = \si'(\om^{n,i}(s,y))(u(s,y)-u(s,y)_0^{(i-\ell_n)\Del})$. Thus, using the Cauchy--Schwarz inequality once more, we obtain
\begin{align*}
\bbe[|\si(u(s,y))-\si(u(s,y)_0^{(i-\ell_n)\Del})|^p]^{\frac1p}\leq \E[|\si'(\om^{n,i}(s,y))|^{2p}]^{\frac1{2p}}\E[|u(s,y)-u(s,y)_0^{(i-\ell_n)\Del}|^{2p}]^{\frac1{2p}}.
\end{align*}
The first term is bounded uniformly in $(s,y)\in [0,T]\times\R^d$ by \eqref{eq:mom} and the assumption that $\si'$ is of polynomial growth. With similar estimates as above, we can bound the second term via
\beq\label{eq:uest}\begin{split}
&\E[|u(s,y)-u(s,y)_0^{(i-\ell_n)\Del}|^{2p}]^{\frac1{2p}}= \E\Biggl[\Biggl|\iint_{(i-\ell_n)\Del}^s G_{y-z}(s-r)\si(u(r,z))\,W(\dd r,\dd z)\Biggr|^{2p}\Biggr]^{\frac{1}{2p}}\\
&\quad\lec  \E\Biggl[\Biggl|\iiint_{(i-\ell_n)\Del}^s G_{y-z}(s-r)G_{y-z'}(s-r)\si(u(r,z))\si(u(r,z'))\,\La(\dd z,\dd z')\,\dd r\Biggr|^{p}\Biggr]^{\frac{1}{2p}}\\
&\quad\leq\Biggl(\iiint_{(i-\ell_n)\Del}^s G_{y-z}(s-r)G_{y-z'}(s-r) \E[|\si(u(r,z))|^{2p}]^{\frac1{2p}}\E[|\si(u(r,z'))|^{2p}]^{\frac{1}{2p}}\,\La(\dd z,\dd z')\,\dd r\Biggr)^{\frac1{2}} \\
&\quad\lec (\ell_n\Del)^{\frac12-\frac\al4},
\end{split}\raisetag{5\baselineskip}\eeq
which implies
\beq\label{eq:siest} \E[|\si(u(s,y))-\si(u(s,y)_{0}^{(i-\ell_n)\Del})|^p]^{\frac1p} \lec (\ell_n\Del)^{\frac12-\frac\al4}. \eeq
The last step in \eqref{eq:uest} uses \eqref{eq:mom}, the linear growth of $\si$, and the fact that 
\beq\label{eq:G-int} \int_0^{\ell_n\Del} G_z(r)G_{z'}(r)\,\La(\dd z,\dd z')\,\dd r \lec (\ell_n\Del)^{1-\frac\al2}, \eeq
which follows from \eqref{eq:Pin1} by setting $r=h=0$ and replacing $\Del$ by $\ell_n\Del$.
We conclude from \eqref{eq:Theta} that
\begin{align*}
\bbe[(U_n)^\ast_T]&\lec  \delh[T/\del] (\ell_n\del)^{\frac12-\frac\al4}\sup_{(r,h)\in\N_0\times\R^d} |\Pi^n_{r,h}|((\ell'_n\del,(\ell_n+L-1)\del)\times\R^d\times\R^d)^{\frac12}\\
&\lec \del^{-\frac12} (\ell_n\del)^{\frac12-\frac\al4}\del^{\frac\Theta 2 \ell'}.
\end{align*}
The three factors in this final bound can be attributed to the components in \eqref{u-expr}. The factor $\Del^{-1/2}$ comes from $\Del^{1/2}\sum_{i=\la_n}^{T(n)}$, the factor $(\ell_n\del)^{1/2-\al/4}$ comes from the difference $\si(u(s,y))-\si(u(s,y)_0^{(i-\ell_n)\Del})$ (essentially as a consequence of the Hölder properties of $u$), and the factor $\del^{\Theta  \ell'/2}$ comes from the stochastic integral of $\Delta^n_i G_y(s)/\tau_n$ on the interval $((i-\ell_n)\Del, (i-\ell'_n)\Del]$ (which is $\ell'_n\Del$ away from $i\Del$).  The variables $h(z^n_i)$  have contribution $1$ by \eqref{size}.  

Therefore, in the following proofs, if we encounter a term like $U_n(t)$, we will simply say that the $\si$-difference in \eqref{u-expr} is of \emph{size} (or \emph{magnitude} or \emph{order}) $(\ell_n\Del)^{1/2-\al/4}$, while the $W$-integral is of size $\Del^{\Theta\ell'/2}$ and the $h$-variables are of size $1$, and then
\emph{directly} conclude that $$\bbe[(U_n)^\ast_T]\lec \del^{-\frac12} (\ell_n\del)^{\frac12-\frac\al4}\del^{\frac\Theta 2 \ell'},$$ without going through similar arguments again. 

We will also apply standard size estimates to expressions that are more complicated than \eqref{u-expr}, for example, if the stochastic integral in \eqref{u-expr} is squared or replaced by the product of two stochastic integrals. Then the generalized Hölder's inequality with exponent $\frac{p}{p-2}$ for $h(z^n_i)$ and exponent $p$ for each appearing stochastic integral can be used to factorize $\bbe[(U_n)^\ast_T]$ in the same manner as before. The key observation is that the total size of such expressions can  always be  determined  component by component. 

Standard size estimates (as explained above) are surprisingly sharp and cannot be further  improved in general. However, there are situations where we can do better. Consider, for instance, 
\beq\label{eq:Uprime} U'_n(t):= \Delh\sum_{i=1}^{t(n)} \om^n_i \eeq
with  random variables $\om^n_i$   that are $\calf^n_i$-measurable and have zero $\calf^n_{i-\ell_n}$-conditional expectation. Further suppose that we know (e.g., via a standard size estimate) that $\om^n_i$ is of size $\Del^{\pi}$ for some $\pi>0$ (uniformly in $i$). Suppose for the moment that $\ell_n=1$. Then $U'_n(t)$ is a martingale sum, so by Doob's inequality and the fact that $\E[\om^n_i\om^n_j] = \E[\om^n_i \E[\om^n_j\mid \calf^n_{j-1}]]= 0$ for all $i<j$, 
\[ \E[((U'_n)^\ast_T)^2] \lec \E[U'_n(T)^2] = \Del \sum_{i,j=1}^{T(n)} \E[\om^n_i\om^n_j] =\Del \sum_{i=1}^{T(n)}\E[(\om^n_i)^2]\lec \Del T(n)\Del^{2\pi} \lec \Del^{2\pi},  \]
which implies by Jensen's inequality that $\E[(U'_n)^\ast_T] \leq \E[((U'_n)^\ast_T)^2]^{1/2} \lec \Del^\pi$. This is better than the standard size estimate $\E[(U'_n)^\ast_T]\lec \Del^{-1/2+\pi}$! Of course, this is due to the martingale structure of $U'_n$ (note that only the terms with $i=j$ remain in the third step of the previous display). 

For general $\ell_n$, we will get something between $\Del^{-1/2+\pi}$ and $\Del^{\pi}$. Indeed, let us rearrange the terms in $U'_n$ in the following way:
\[ U'_n(t)=\sum_{k=1}^{\ell_n} V^n_k(t),\qquad V^n_k(t):=\Delh\sum_{j=1}^{[t(n)/\ell_n]} \om^n_{k+(j-1)\ell_n}.  \]
The point is now that for fixed $k$, the term $V^n_k(t)$ is a martingale sum relative to the filtration $(\calf^n_{k+(j-1)\la_n})_{j=1,\dots,[t(n)/\ell_n]}$. In fact, $\om^n_{k+(j-1)\ell_n}$ is $\calf^n_{k+(j-1)\la_n}$-measurable with vanishing $\calf^n_{k+(j-2)\la_n}$-conditional expectation. Thus, in analogy to the simple case discussed above, we have
\[ \E[((V^n_k)^\ast_T)^2]\lec \E[V^n_k(T)^2] = \Del \sum_{j=1}^{[T(n)/\ell_n]} \E[(\om^n_{k+(j-1)\ell_n})^2] \lec \Del[T(n)/\ell_n]\Del^{2\pi}\lec \ell_n^{-1}\Del^{2\pi}. \]
As a result, again by Jensen's inequality,
\[ \E[(U'_n)^\ast_T] \leq \sum_{k=1}^{\ell_n} \E[(V^n_k)^\ast_T]\lec \sum_{k=1}^{\ell_n} \ell_n^{-\frac12}\Del^{\pi} = \ell_n^{\frac12}\Del^{\pi} \lec \Del^{-\frac\ell2+\pi}, \]
which is, as promised, between $\Del^{-1/2+\pi}$ and $\Del^{\pi}$ (recall that $0<\ell<1$). In summary, such a \emph{martingale size estimate} is applicable whenever we have a sum of terms that are conditionally independent of each other to a certain degree. The main parameter is the length of overlap (this was $\ell_n$ above), that is, the distance at which the terms start to become conditionally independent. The conclusion is then that the sum $\Delh\sum_{i=1}^{t(n)}$, which by a standard size estimate would be of magnitude $\Del^{-1/2}$, is effectively only of order $\sqrt{\ell_n}$. In what follows, very often $\om^n_i$ is actually only $\calf^n_{i+L-1}$-measurable. But we can easily convince ourselves that a shift like this by finitely many units does not affect the order of the estimates at all.

\bpr[Proof of Lemma~\ref{lem:remove-lan}] What remains from the left-hand side of \eqref{eq:remove-lan} is $\Delh$ times a sum of $\la_n$ many $f$-terms, which are all of size $1$ [by \eqref{eq:mom}, \eqref{eq:Pin1}, and the polynomial growth of $f$]. So the total size is $\sqrt{\Del}\la_n \leq \Del^{1/2-a}$, which tends to $0$ as $n\to\infty$ if $a$ is close enough to $\frac1{2\Theta}<\frac12$.
\epr

\bpr[Proof of Lemma~\ref{lem:trunc}] Introduce a sequence of intermediate truncation levels $\la_n^r:=[\Del^{-a_r}]$, where $r=0,\dots,R$, by choosing the numbers 
$a_0>\dots > a_R$ in such a way that $1>a_0>\frac1\Theta$, $a_R=a$, and $a_r>\frac{a_{r-1}}{\Theta}$ for all $r=1,\dots,R$ (this is possible because $\Theta=1+\frac\al2>1$). We can then write the left-hand side of \eqref{eq:approx-ga} as $A^n_1(t)+A^n_2(t)+ \sum_{r=1}^R (A^{n,r}_3(t)+A^{n,r}_4(t))$ where
\begin{align*}
A^n_1(t)&:=\Delh \sumtla \Biggl\{ f\biggl( \frac{\Delta^n_i u}{\tau_n},\dots, \frac{\Delta^n_{i+L-1} u}{\tau_n}\biggr) - f\Biggl( \frac{\Delta^n_i u^{(i+L-1)\Del}_0}{\tau_n}\Biggr) \Biggr\},\\
A^n_2(t)&:=\Delh \sumtla \Biggl\{ f\biggl( \frac{\Delta^n_i u^{(i+L-1)\Del}_0}{\tau_n}\biggr) - f\Biggl( \frac{\Delta^n_i u^{(i+L-1)\Del}_{(i-\la_n^0)\Del}}{\tau_n}\Biggr) \Biggr\},\\
A^{n,r}_3(t)&:=\Delh\sumtla \ov\al^{n,r}_i, \qquad A^{n,r}_4(t):=\Delh\sumtla \bbe\Bigl[\al^{n,r}_i\mathrel{\big|} \calf^n_{i-\la_n^{r-1}}\Bigr],
\end{align*}
and, for $i=\la_n+1,\dots,[t/\Del]$,
\beq\label{eq:al}  \al^{n,r}_i := f\Biggl( \frac{\Delta^n_i u^{(i+L-1)\Del}_{(i-\la_n^{r-1})\Del}}{\tau_n}\Biggr) - f\Biggl(\frac{\Delta^n_i u^{(i+L-1)\Del}_{(i-\la_n^r)\Del}}{\tau_n}\Biggr),\qquad \ov\al^{n,r}_i := \al^{n,r}_i - \bbe\Bigl[\al^{n,r}_i \mathrel{\big|} \calf^n_{i-\la_n^{r-1}}\Bigr]. \eeq

It was shown in \cite[Lemmas 3.6 and 3.7]{Chong19} that $A^n_2$ and $A^{n,r}_3$ (for every $r$) converge to $0$ in $L^1$, uniformly on compacts. But let us include the short proofs for the reader's convenience (and to get used to standard and martingale size estimates). By the mean value theorem, we have
\beq\label{eq:help3}
A^n_2(t)=\Delh\sumtla \sum_{\nu\in\cali} \partial_{\nu}f(\wh\al^{n,i}_1)\Biggl( \frac{\Delta^n_i u^{(i+L-1)\Del}_0-\Delta^n_i u^{(i+L-1)\Del}_{(i-\la_n^0)\Del}}{\tau_n}\Biggr)_\nu
\eeq
for some $\wh\al^{n,i}_1$ between $\Delta^n_i u^{(i+L-1)\Del}_0/\tau_n$ and $\Delta^n_i u^{(i+L-1)\Del}_{(i-\la_n^0)\Del}/\tau_n$.
By a standard size estimate and \eqref{eq:Theta}, the difference in parentheses is of order $\Del^{\Theta a_0/2}$, while $\partial_\nu f(\wh \al^n_i)$ is of order $1$. Thus,
$\E[(A^n_2)^\ast_T]\lec \Del^{-1/2+\Theta /2  a_0}$, which converges to $0$ for every $T>0$ because $a_0>\frac1\Theta$ by assumption.

By a similar argument, each term $\al^{n,r}_i$ is of size $\Del^{\Theta a_r/2}$, hence also $\ov\al^{n,r}_i$ by the contraction property of the conditional expectation on $L^p$-spaces. Moreover, $\ov\al^{n,r}_i$ is $\calf^n_{i+L-1}$-measurable with zero conditional expectation given $\calf^n_{i-\la_n^{r-1}}$, so by a martingale size argument and the hypothesis that $a_r>\frac{a_{r-1}}{\Theta}$, 
$$\E[(A^{n,r}_3)^\ast_T]\lec\sqrt{\la^{r-1}_n}\Del^{\frac\Theta 2 a_r} = \Del^{-\frac{a_{r-1}}{2}+\frac{\Theta}{2}a_r} \to 0.$$

Concerning $A^n_1$, notice that the only difference between $\Delta^n_i u$ and $\Delta^n_i u^{(i+L-1)\Del}_0$ is that the former takes into account the increment of $u^{(0)}$ from \eqref{eq:initial}, while the latter does not. As explained in \cite[Remark~2.7]{Chong19}, hypothesis H4 implies that $|u^{(0)}(i\Del,x)-u^{(0)}((i-1)\Del,x)|/\tau_n=o(\Delh)$, uniformly in $i=\la_n+1,\dots,T(n)$. Hence, using the mean value theorem and a standard size estimate [similar to \eqref{eq:help3}], we deduce $\E[(A^n_1)^\ast_T]\to0$. 

For the last term $A^{n,r}_4$, the proof of \cite[Lemma 3.8]{Chong19} breaks down because $\si(u(t,x))$ is not $\frac12$-Hölder continuous in time (in general), so we need another argument in the multiplicative case. 
The idea is to apply Taylor's theorem twice and write
\beq\label{deltanri} \begin{split} \al^{n,r}_i &= \sum_{\nu\in\cali} \partial_\nu f \Biggl(\frac{\Delta^n_i u^{(i+L-1)\Del}_{(i-\la_n^r)\Del}}{\tau_n} \Biggr) \frac{\Bigl(\Delta^n_i u_{(i-\la_n^{r-1})\Del}^{(i-\la_n^r)\Del}\Bigr)_\nu}{\tau_n} +\al^{n,r,4}_i=\sum_{j=1}^4
	\al^{n,r,j}_i , \end{split}\eeq
where
\begin{align*}
\al^{n,r,1}_i &:=\sum_{\nu\in\cali} \partial_\nu f \Biggl(\frac{\Delta^n_i u^{(i+L-1)\Del}_{(i-\la_n^r)\Del,(i-\la_n^{r-1})\Del}}{\tau_n} \Biggr) \frac{\Bigl(\Delta^n_i u_{(i-\la_n^{r-1})\Del}^{(i-\la_n^r)\Del}\Bigr)_\nu}{\tau_n}, \\
\al^{n,r,2}_i &:=  \sum_{\nu_1,\nu_2\in\cali} \partial^2_{\nu_1\nu_2} f \Biggl(\frac{\Delta^n_i u^{(i+L-1)\Del}_{(i-\la_n^r)\Del,(i-\la_n^{r-1})\Del}}{\tau_n} \Biggr)\frac{\Bigl(\Delta^n_i u_{(i-\la_n^{r-1})\Del}^{(i-\la_n^r)\Del}\Bigr)_{\nu_1}}{\tau_n}\\
&\mathrel{\hphantom{:=} \mathrel{\times}}\frac{\Bigl(\Delta^n_i u_{(i-\la_n^r)\Del}^{(i+L-1)\Del}-\Delta^n_i u^{(i+L-1)\Del}_{(i-\la_n^r)\Del,(i-\la_n^{r-1})\Del}\Bigr)_{\nu_2}}{\tau_n},\\
\al^{n,r,3}_i &:= \frac12\sum_{\nu_1,\nu_2,\nu_3\in\cali} \partial^3_{\nu_1\nu_2\nu_3} f(\wh\al^{n,i}_2)\frac{\Bigl(\Delta^n_i u_{(i-\la_n^{r-1})\Del}^{(i-\la_n^r)\Del}\Bigr)_{\nu_1}}{\tau_n}\frac{\Bigl(\Delta^n_i u_{(i-\la_n^r)\Del}^{(i+L-1)\Del}-\Delta^n_i u^{(i+L-1)\Del}_{(i-\la_n^r)\Del,(i-\la_n^{r-1})\Del}\Bigr)_{\nu_2}}{\tau_n}\\
&\mathrel{\hphantom{:=} \mathrel{\times}}\frac{\Bigl(\Delta^n_i u_{(i-\la_n^r)\Del}^{(i+L-1)\Del}-\Delta^n_i u^{(i+L-1)\Del}_{(i-\la_n^r)\Del,(i-\la_n^{r-1})\Del}\Bigr)_{\nu_3}}{\tau_n},\\
\al^{n,r,4}_i&:=\frac{1}{2} \sum_{\nu_1,\nu_2\in\cali} \partial^2_{\nu_1\nu_2} f (\wh\al^{n,i}_3) \frac{\Bigl(\Delta^n_i u_{(i-\la_n^{r-1})\Del}^{(i-\la_n^r)\Del}\Bigr)_{\nu_1}\Bigl(\Delta^n_i u_{(i-\la_n^{r-1})\Del}^{(i-\la_n^r)\Del}\Bigr)_{\nu_2}}{\tau_n^2},
\end{align*}
and $\wh\al^{n,i}_2$ and $\wh\al^{n,i}_3$ are some intermediate points. We obtain $A^{n,r}_3=\sum_{j=1}^4 A^{n,r}_{3,j}$ by setting 
$$A^{n,r}_{3,j}(t):=\Delh\sumtla \E[\al^{n,r,j}_i\mid \calf^n_{i-\la_n^{r-1}}],\qquad j=1,\dots,4.$$ 

The last two quotients in $\al^{n,r,3}_i$ are both of order $(\la_n^{r-1}\Del)^{1/2-\al/4}$ by \eqref{eq:siest}. Moreover, the term $\Delta^n_i u_{(i-\la_n^{r-1})\Del}^{(i-\la_n^r)\Del}/\tau_n$ is of size $\Del^{\Theta a_r/2}$ by \eqref{eq:Theta}. Because $\al\leq1$ and $a_r>\frac{a_{r-1}}{\Theta}$, it follows that
\beq\label{eq:A33} \E[(A^{n,r}_{3,3})^\ast_T]\lec \Del^{-\frac12+\frac\Theta 2 a_r +2(1-a_{r-1})(\frac12-\frac\al4)} \leq\Del^{-\frac12+\frac\Theta 2 a_r+\frac12-\frac{a_{r-1}}{2}} =\Del^{\frac{\Theta a_r-a_{r-1}}{2}}\to0.  \eeq
Similarly, we have $\E[(A^{n,r}_{3,4})^\ast_T]\lec \Del^{-1/2+\Theta a_r}\to0$ because $a_r\geq a_R=a>\frac1{2\Theta}$.

Concerning $A^{n,r}_{3,1}$, we use the tower property of the conditional expectation, which leads to 
\begin{align*}
&\bbe\Bigl[\al^{n,r,1}_i \mathrel{\big|} \calf^n_{i-\la_n^{r-1}}\Bigr]\\ 
&\quad= 
\sum_{\nu\in\cali} \E\Biggl[ \E\Biggl[\partial_\nu f \Biggl(\frac{\Delta^n_i u^{(i+L-1)\Del}_{(i-\la_n^r)\Del,(i-\la_n^{r-1})\Del}}{\tau_n} \Biggr)\mathrel{\bigg|} \calf^n_{i-\la_n^r}\Biggr] \frac{\Bigl(\Delta^n_i u_{(i-\la_n^{r-1})\Del}^{(i-\la_n^r)\Del}\Bigr)_\nu}{\tau_n}\mathrel{\Bigg|}\calf^n_{i-\la_n^{r-1}}\Biggr].
\end{align*}
Notice now that the inner conditional expectation, and hence $A^{n,r}_{3,1}$, is $0$ because $\partial_\nu f$ is an odd function and its argument follows a centered normal law, conditionally on $\calf^n_{i-\la_n^r}$.

For $A^{n,r}_{3,2}$, the situation is more complicated. Indeed, if we used a standard size argument, the $\partial^2_{\nu_1\nu_2}f$-term in the definition of $\al^{n,r,2}_i$ would be of order $1$, and the next two factors would be of order $\Del^{\Theta a_r/2}$ and $(\la_n^{r-1}\Del)^{1/2-\al/4}$. So in the borderline case, if $\al=1$ (so $\Theta=\frac32$), we  only get 
\beq\label{eq:help11} \E[(A^{n,r}_{3,2})^\ast_T]\lec \Del^{-\frac12+\frac\Theta 2 a_r+(1-a_{r-1})(\frac12-\frac\al4)} = \Del^{-\frac14+\frac34 a_r-\frac{a_{r-1}}{4}}.\eeq 
Since $a_{r-1}>a_r$, the exponent is clearly less than $-\frac14+\frac12 a_r$, so if $a_r$ is close to $\frac{1}{2\Theta}=\frac13$, it becomes negative. This shows that for $\al$ close to $1$, a standard size estimate is not sufficient.

Instead, in order to show that $A^{n,r}_{3,2}$ is negligible, we have to use, for the first time, that $\sigma(u(t,x))$ is a function of the solution process itself. As an auxiliary step, we replace $\al^{n,r,2}_i$ by $\al^{\prime n,r,2}_i$, which is defined in the same way as the former but with $(\Delta^n_i u_{(i-\la_n^{r-1})\Del}^{(i-\la_n^r)\Del})_{\nu_1}$ substituted by $(\Delta^n_i u_{(i-\la_n^{r-1})\Del,(i-\la_n^{r-1})\Del}^{(i-\la_n^r)\Del})_{\nu_1}$. Let $A^{\prime n,r}_{3,2}$ denote the corresponding modified version of $A^{n,r}_{3,2}$. A standard size estimate shows that the incurred error is of order $\Del^{-1/2+\Theta a_r/2}(\la_n^{r-1}\Del)^{1/2-\al/4}$ [as in \eqref{eq:help11}] times an additional factor $(\la_n^{r-1}\Del)^{1/2-\al/4}$ (due to the substitution). In total, by \eqref{eq:A33}, 
\beq\label{eq:A32prime} \E[(A^{n,r}_{3,2}-A^{\prime n,r}_{3,2})^\ast_T]\lec\Del^{-\frac12+\frac\Theta 2 a_r+(1-a_{r-1})(1-\frac\al2)}\to0.\eeq 

The last fraction in the definition of $\al^{n,r,2}_i$ (which also appears  in $\al^{\prime n,r,2}_i$) can be rewritten as
\begin{align*}
&\iint_{(i-\la_n^r)\Del}^{(i+l_2-1)\Del} \frac{\Delta^n_i G_y(s)_{\nu_2}}{\tau_n}\Bigl(\si(u(s,y))-\si(u(s,y)_0^{(i-\la_n^{r-1})\Del})\Bigr)\,W(\dd s,\dd y)\\
&\quad=\iint_{(i-\la_n^r)\Del}^{(i+l_2-1)\Del} \frac{\Delta^n_i G_y(s)_{\nu_2}}{\tau_n}\si'(u(s,y)_0^{(i-\la_n^{r-1})\Del})\Bigl(u(s,y)-u(s,y)_0^{(i-\la_n^{r-1})\Del}\Bigr)\,W(\dd s,\dd y)\\
&\qquad+\iint_{(i-\la_n^r)\Del}^{(i+l_2-1)\Del} \frac{\Delta^n_i G_y(s)_{\nu_2}}{\tau_n}\frac{\si''(\wh \al^{n,i}_3(s,y))}{2}\Bigl(u(s,y)-u(s,y)_0^{(i-\la_n^{r-1})\Del}\Bigr)^2\,W(\dd s,\dd y)
\end{align*}
for some $\wh\al^{n,i}_3(s,y)$ between $u(s,y)$ and $u(s,y)_0^{(i-\la_n^{r-1})\Del}$. Because of the square, the contribution of the second term can be bounded as in \eqref{eq:A32prime} and is therefore negligible.
The first term equals
\begin{align*}
&\iint_{(i-\la_n^r)\Del}^{(i+l_2-1)\Del} \frac{\Delta^n_i G_y(s)_{\nu_2}}{\tau_n}\si'(u(s,y)_0^{(i-\la_n^{r-1})\Del})\\
&\quad\quad\times\iint_{(i-\la_n^{r-1})\Del}^s G_{y-z}(s-r)\si(u(r,z)_0^{(i-\la_n^{r-1})\Del})\,W(\dd r,\dd z)\,W(\dd s,\dd y)\\
&\quad+\iint_{(i-\la_n^r)\Del}^{(i+l_2-1)\Del} \frac{\Delta^n_i G_y(s)_{\nu_2}}{\tau_n}\si'(u(s,y)_0^{(i-\la_n^{r-1})\Del})\iint_{(i-\la_n^{r-1})\Del}^s G_{y-z}(s-r)\\
&\quad\quad\times\Bigl(\si(u(r,z))-\si(u(r,z)_0^{(i-\la_n^{r-1})\Del})\Bigr)\,W(\dd r,\dd z)\,W(\dd s,\dd y).
\end{align*}
Again, the second summand is negligible as it has an additional factor $\la_n^{r-1}\Del$ from the $\si$-difference. So we are left to consider $A^{\prime\prime n,r}_{3,2}$, which is defined in the same way as $A^{\prime n,r}_{3,2}$ but with the last factor in $\al^{n,r,2}_i$ replaced by the first summand of the previous display. Written out explicitly, 
\begin{align*}
A^{\prime\prime n,r}_{3,2}(t)&=\sqrt{\Del}\sumtla\sum_{\nu_1,\nu_2\in\cali}\E\Biggl[\partial^2_{\nu_1\nu_2} f \Biggl(\frac{\Delta^n_i u^{(i+L-1)\Del}_{(i-\la_n^r)\Del,(i-\la_n^{r-1})\Del}}{\tau_n} \Biggr)\frac{\Bigl(\Delta^n_i u_{(i-\la_n^{r-1})\Del,(i-\la_n^{r-1})\Del}^{(i-\la_n^r)\Del}\Bigr)_{\nu_1}}{\tau_n}\\
&\quad\times\iint_{(i-\la_n^r)\Del}^{(i+l_2-1)\Del} \frac{\Delta^n_i G_y(s)_{\nu_2}}{\tau_n}\si'(u(s,y)_0^{(i-\la_n^{r-1})\Del})\iint_{(i-\la_n^{r-1})\Del}^s G_{y-z}(s-r)\\
&\quad\quad\times\si(u(r,z)_0^{(i-\la_n^{r-1})\Del})\,W(\dd r,\dd z)\,W(\dd s,\dd y)\mathrel{\Bigg|}\vphantom{\Biggl(\frac{\Delta^n_i u^{(i+L-1)\Del}_{(i-\la_n^r)\Del,(i-\la_n^{r-1})\Del}}{\tau_n} \Biggr)}\calf^n_{i-\la_n^{r-1}}\Biggr].
\end{align*}
The crucial observation is now that under the $\calf^n_{i-\la_n^{r-1}}$-conditional probability measure, the double $W$-integral (resp., the preceding fraction) is an element of the second (resp., first) Wiener chaos, while the $\partial^2_{\nu_1\nu_2}f$-term, because $f$ is even, belongs to the direct sum of all Wiener chaoses of even order. Consequently, the product of the three terms belongs to the direct sum of all odd-order Wiener chaoses and therefore has zero $\calf^n_{i-\la_n^{r-1}}$-conditional expectation. The aforementioned statements concerning Wiener chaoses can be found in \cite[Propositions~1.1.3 and 1.1.4]{Nualart06}.
\epr

\bpr[Proof of Lemma~\ref{lem:CLT-core}] Let
\begin{align*}
C^{n,m}_3(t):=\Delh\sum_{i=J^{n,m}(t)((m+1)\la_n+L-1)+1}^{t(n)} \ga^n_{i,k^{n,m}_i},\qquad
J^{n,m}(t):=\biggl[\frac{t(n)}{(m+1)\la_n+L-1}\biggr].
\end{align*}
The sum defining $C^{n,m}_3(t)$ contains at most $(m+1)\la_n+L-1$ terms. As $\E[|\ga^n_{i,k}|]$ is uniformly bounded, we have $\E[(C^{n,m}_3)^\ast_T] \lec \Delh \la_n = \Del^{1/2-a}$, which converges to $0$ when $a$ is close enough to $\frac{1}{2\Theta}$. Moreover, a close inspection of the respective definitions reveals that $C^{n,m}-C^{n,m}_3$ is equal to $\wh V^{n,m,1}+ \wh V^{n,m,2}$ (modulo the first $\la_n$ terms, which we have removed in Lemma~\ref{lem:remove-lan} and which are negligible by the same reason), where $\wh V^{n,m,1}$ and $\wh V^{n,m,2}$ are the terms defined and subsequently analyzed in \cite[(3.11)]{Chong19} and \cite[(D.8)]{Chong19a}, respectively. So with identical proofs as in the mentioned references [they do not rely on the Hölder properties of $\si(u(s,y))$!], the statement of the lemma follows upon setting $C^{n,m}_1:=\wh V^{n,m,1}$ and $C^{n,m}_2:=\wh V^{n,m,2}$ (with the first $\la_n$ terms removed).
\epr

\bpr[Proof of Lemma~\ref{lem:B13}] The argument is the same as for $A^{n,r}_3$ in the proof of Lemma~\ref{lem:trunc}: as $\beta^n_i-\ov\beta^n_i$ and $\ov\beta^n_{i,k}-\beta^n_{i,k}$ are both of order $\Del^{\Theta \ov a/2}$ [and thus, also $f(\beta^n_i)-f(\ov\beta^n_i)$ and $f(\ov\beta^n_{i,k})-f(\beta^n_{i,k})$ by the mean value theorem], a martingale estimate shows $\E[(B^n_1)^\ast_T]+\E[(B^{n,m}_3)^\ast_T]\lec \la_n^{1/2}\Del^{\Theta \ov a/2} \leq \Del^{(\Theta\ov a-a)/2}$, which vanishes as $n\to\infty$ because $\ov a>\frac{a}{\Theta}$ by hypothesis. 
\epr

\bpr[Proof of Lemma~\ref{lem:D45}] A standard size estimate gives 
\beq\label{eq:help13}\E[(D^n_4)^\ast_T]\lec \Del^{-\frac12}((\la_n\Del)^{\frac12-\frac\al4})^3 = \Del^{-\frac12+3(1-a)(\frac12-\frac\al4)}.\eeq
This converges to $0$ as $n\to\infty$ if $a$ is sufficiently close to $\frac{1}{2\Theta}=\frac1{2+\al}$, because the last exponent in \eqref{eq:help13} gets arbitrarily close to $-\frac12+3(1-\frac{1}{2+\al})(\frac12-\frac\al4)$, which is strictly positive for $\al\in(0,1)$. 
\epr

\bpr[Proof of Lemma~\ref{lem:remove-s}] 
For  $\nu_1$ and $\nu_2$ fixed, as a consequence of the identity 
\beq\label{eq:id2}xy-x_0y_0=y_0(x-x_0)+x(y-y_0),\eeq 
we have $v^{n,i}_{\nu_1,\nu_2} - v^{\prime n,i}_{\nu_1,\nu_2}=E^{n,i}_1+E^{n,i}_2+E^{n,i}_3$, where
\begin{align*}
E^{n,i}_1&:= \iiint_0^{\la_n\Del} \si(u(i\Del,x_{k_{12}}-y')_0^{(i-\la_n)\Del})\Bigl[\si(u(t^{n,i}_{l_{12}-1}-s,x_{k_{12}}-y)_0^{(i-\la_n)\Del})\\
&\mathrel{\hphantom{:=} \mathrel{-}} \si(u(i\Del,x_{k_{12}}-y)_0^{(i-\la_n)\Del})\Bigr]\,\Pi^n_{|l_1-l_2|,x_{k_1}-x_{k_2}}(\dd s,\dd y,\dd y'),\\
E^{n,i}_2&:=\iiint_0^{\la_n\Del} \si(u(t^{n,i}_{l_{12}-1}-s,x_{k_{12}}-y)_0^{(i-\la_n)\Del})\Bigl[\si(u(t^{n,i}_{l_{12}-1}-s,x_{k_{12}}-y')_0^{(i-\la_n)\Del})\\
&\mathrel{\hphantom{:=} \mathrel{-}} \si(u(i\Del,x_{k_{12}}-y')_{0}^{(i-\la_n)\Del})\Bigr]\,\Pi^n_{|l_1-l_2|,x_{k_1}-x_{k_2}}(\dd s,\dd y,\dd y'),\\
E^{n,i}_3&:=\iiint_{\la_n\Del}^{(\la_n+l_{12}-1)\Del} \si(u(t^{n,i}_{l_{12}-1}-s,x_{k_{12}}-y)_0^{(i-\la_n)\Del})\\
&\mathrel{\hphantom{:=} \mathrel{\times}}\si(u(t^{n,i}_{l_{12}-1}-s,x_{k_{12}}-y')_0^{(i-\la_n)\Del})\,\Pi^n_{|l_1-l_2|,x_{k_1}-x_{k_2}}(\dd s,\dd y,\dd y').
\end{align*}
If we manage to prove that $E^{n,i}_1$, $E^{n,i}_2$, and $E^{n,i}_3$ are of size $o(\Delh)$ uniformly in $i$, then, because $\un\mu_f$ is differentiable by \cite[(D.46)]{Chong19a}, \eqref{eq:help12} follows from the mean value theorem. A standard size estimate immediately gives $\E[|E^{n,i}_3|]\lec \Del^{\Theta a} = o(\Delh)$ because $a>\frac{1}{2\Theta}$. Concerning the other two terms, we  only consider $E^{n,i}_1$ further as $E^{n,i}_2$ can be treated in a completely analogous manner.

Using the mean value theorem on $\si$, we can find $\eps^{n,i}_2(s,y)$ between   $u(t^{n,i}_{l_{12}-1}-s,x_{k_{12}}-y)_0^{(i-\la_n)\Del}$ and $u(i\Del,x_{k_{12}}-y)_0^{(i-\la_n)\Del}$ such that 
\begin{align*} E^{n,i}_1 &= \iiint_0^{\la_n\Del} \si(u(i\Del,x_{k_{12}}-y')_0^{(i-\la_n)\Del})\si'(\eps^{n,i}_2(s,y))\\ 
&\quad\times \Bigl[u(t^{n,i}_{l_{12}-1}-s,x_{k_{12}}-y)_0^{(i-\la_n)\Del}-u(i\Del,x_{k_{12}}-y)_0^{(i-\la_n)\Del}\Bigr]\,\Pi^n_{|l_1-l_2|,x_{k_1}-x_{k_2}}(\dd s,\dd y,\dd y')\\
&= \iiint_0^{\la_n\Del} \si(u(i\Del,x_{k_{12}}-y')_0^{(i-\la_n)\Del})\si'(\eps^{n,i}_2(s,y))\Biggl(\iint_0^{(i-\la_n)\Del} \Bigl( G_{x_{k_{12}}-y-z}(t^{n,i}_{l_{12}-1}-s-r)\\
&\quad-  G_{x_{k_{12}}-y-z}(i\Del-r) \Bigr)\si(u(r,z))\,W(\dd r,\dd z)\Biggr)\,\Pi^n_{|l_1-l_2|,x_{k_1}-x_{k_2}}(\dd s,\dd y,\dd y').
\end{align*}

Next, in the spirit of \cite[Lemma~3]{BN11} and \cite[Lemma~3.14]{Chong19}, let $a^{(0)},\dots,a^{(Q)}$ be numbers such that $a=:a^{(0)}>\dots>a^{(Q-1)}>0$ and $a^{(Q)}:=-\infty$ and define $\la_n^{(q)}:=[\Del^{-a^{(q)}}]$ for $q=0,\dots,Q$. Then 
\allowdisplaybreaks[4]
\begin{align*} E^{n,i}_1&= \sum_{q=1}^Q  \iiint_{\la_n^{(q)}\Del}^{\la_n^{(q-1)}\Del} \si(u(i\Del,x_{k_{12}}-y')_0^{(i-\la_n)\Del})\si'(\eps^{n,i}_2(s,y))\\
&\quad\times\Biggl(\iint_0^{(i-\la_n)\Del} \Bigl( G_{x_{k_{12}}-y-z}(t^{n,i}_{l_{12}-1}-s-r) \\
&\qquad-  G_{x_{k_{12}}-y-z}(i\Del-r) \Bigr)\si(u(r,z))\,W(\dd r,\dd z)\Biggr)\,\Pi^n_{|l_1-l_2|,x_{k_1}-x_{k_2}}(\dd s,\dd y,\dd y'). \end{align*}  
If $s\in(\la_n^{(q)}\Del,\la_n^{(q-1)}\Del]$, then ${|t^{n,i}_{l_{12}-1}-s-i\Del|}\lec \la_n^{(q-1)}\Del$. Moreover, as $\la_n\Del$ is of the same order as $(\la_n^{(q-1)}\Del)^{1- a_0^{(q-1)}}$, where $a_0^{(q-1)}:=(a-a^{(q-1)})/(1-a^{(q-1)})$, it follows from \eqref{eq:Theta} [with $\la_n^{(q-1)}\Del$ and $a_0^{(q-1)}$ in the roles of $\Del$ and $a$, respectively] that the $W$-integral above has size 
\beq\label{eq:sizeW} (\la_n^{(q-1)}\Del)^{\frac12-\frac\al4}(\la_n^{(q-1)}\Del)^{\frac\Theta2 a_0^{(q-1)}}\leq \Del^{(1-a^{(q-1)})(1-\frac\Theta2)+\frac\Theta2(a-a^{(q-1)})}
=\Del^{1-a^{(q-1)}-\frac\Theta2(1-a)}. 
\eeq
At the same time, for $q\neq Q$, the $\Pi^n$-integral on $(\la_n^{(q)}\Del,\la_n^{(q-1)}\Del]\times\R^d\times\R^d$ contributes a factor of $\Del^{\Theta a^{(q)}}$, while it is of order $1$ for $q=Q$. Thus, by a standard size estimate, 
\[ \E[|E^{n,i}_1|] \lec \Del^{1-a^{(Q-1)}-\frac\Theta2(1-a)}+\sum_{q=1}^{Q-1} \Del^{1-a^{(q-1)}-\frac\Theta2(1-a)+\Theta a^{(q)}}.  \]
This is $o(\Delh)$ if  $a^{(Q-1)}<\frac12-\frac\Theta2(1-a)$ and
\beq\label{eq:rec1} a^{(q-1)}<\frac12-\frac\Theta2(1-a)+\Theta a^{(q)} \eeq
for all $q=1,\dots,Q-1$. As $a>\frac1{2\Theta}$ and $\Theta\leq \frac32$, the interval $(0,\frac12-\frac\Theta2(1-a))$ is nonempty, allowing us to choose $a^{(Q-1)}$ herein. Next, consider the linear recurrence relation $b_{n+1}=g(b_n)$ with $g(x):=\frac12-\frac\Theta2(1-a)+\Theta x$ and $b_1=a^{(Q-1)}$. The only fixed point of $g$ is $\frac{\Theta(1-a)-1}{2\Theta-2}$, which is negative because $\Theta\in(1,\frac32]$ and $a>\frac1{2\Theta}$. 
Moreover, $g$ has positive slope and satisfies $g(0)>0$. As we start with a strictly positive $b_1$, we conclude that $b_n$ increases to infinity as $n\to\infty$. Thus, letting $Q$ be the smallest $n$ such that $b_n>\frac1{2\Theta}$, we can meet condition \eqref{eq:rec1} by choosing $a^{(Q-n)}$  slightly smaller than $b_n$ for $n=2,\dots,Q-1$.
\epr

\bpr[Proof of Lemma~\ref{lem:remove-y}] Applying Itô's formula to the product of the two $\si$-terms in the last integral of \eqref{eq:v} (with $i\Del$ instead of $t^{n,i}_{l_{12}-1}-s$) and noting that $\Pi^n$ is symmetric with respect to swapping $\dd y$ and $\dd y'$, we have, for every $\nu_1$ and $\nu_2$, that $v^{\prime n,i}_{\nu_1,\nu_2}- v^{\prime\prime n,i}_{\nu_1,\nu_2} = F^{n,i}_1+F^{n,i}_2+F^{n,i}_3$, where
\beq\label{eq:delta-split}\begin{split}
	F^{n,i}_1&:=2\iiint_0^{\la_n\Del} \iint_{0}^{(i-\la_n)\Del} \Bigl[G_{x_{k_{12}}-y-z}(i\Del-r)\si'(u(i\Del,x_{k_{12}}-y)_{0}^r)\si(u(i\Del,x_{k_{12}}-y')_{0}^r)\\
	&\mathrel{\hphantom{:=}}\quad{\mathrel{-} G_{x_{k_{12}}-z}(i\Del-r)\si'(u(i\Del,x_{k_{12}})_{0}^r)\si(u(i\Del,x_{k_{12}})_{0}^r)\Bigr]}\\
	&\mathrel{\hphantom{:=} \mathrel{\times}} \si(u(r,z))\,W(\dd r,\dd z)\,\Pi^n_{|l_1-l_2|,x_{k_1}-x_{k_2}}(\dd s,\dd y,\dd y'),\\
	F^{n,i}_2	&:=\iiint_0^{\la_n\Del}\iiint_{0}^{(i-\la_n)\Del} \Bigl[G_{x_{k_{12}}-y-z}(i\Del-r)G_{x_{k_{12}}-y-z'}(i\Del-r)\\
	&\mathrel{\hphantom{:=}}\quad{\mathrel{\times}\si''(u(i\Del,x_{k_{12}}-y)_{0}^r)\si(u(i\Del,x_{k_{12}}-y')_{0}^r)-G_{x_{k_{12}}-z}(i\Del-r)}\\
	&\mathrel{\hphantom{:=}}\quad{\mathrel{\times} G_{x_{k_{12}}-z'}(i\Del-r)\si''(u(i\Del,x_{k_{12}})_{0}^r)\si(u(i\Del,x_{k_{12}})_{0}^r)\Bigr]}\\
	&\mathrel{\hphantom{:=} \mathrel{\times}}\si(u(r,z))\si(u(r,z'))\,\La(\dd z,\dd z')\,\dd r\,\Pi^n_{|l_1-l_2|,x_{k_1}-x_{k_2}}(\dd s,\dd y,\dd y'),\\
	F^{n,i}_3	&:=\iiint_0^{\la_n\Del}\iiint_{0}^{(i-\la_n)\Del} \Bigl[ G_{x_{k_{12}}-y-z}(i\Del-r)G_{x_{k_{12}}-y'-z'}(i\Del-r)\si'(u(i\Del,x_{k_{12}}-y)_{0}^r)\\
	&\mathrel{\hphantom{:=}}\quad{\mathrel{\times}\si'(u(i\Del,x_{k_{12}}-y')_{0}^r)-G_{x_{k_{12}}-z}(i\Del-r)G_{x_{k_{12}}-z'}(i\Del-r)(\si')^2(u(i\Del,x_{k_{12}})_{0}^r)\Bigr]}\\
	&\mathrel{\hphantom{:=} \mathrel{\times}} \si(u(r,z))\si(u(r,z'))\,\La(\dd z,\dd z')\,\dd r\,\Pi^n_{|l_1-l_2|,x_{k_1}-x_{k_2}}(\dd s,\dd y,\dd y').
\end{split}\raisetag{9\baselineskip}\eeq

Each $F^{n,i}_j$, $j=1,2,3$, involves in brackets a difference of products of $G$- and $\si$-terms. So we can use the identity
\beq\label{eq:id3} xy-x_0y_0=y_0(x-x_0)+x_0(y-y_0)+(x-x_0)(y-y_0) \eeq
to further decompose each $F^{n,i}_j$ into three integrals $F^{n,i}_{j,1}+F^{n,i}_{j,2}+F^{n,i}_{j,3}$ such that the first (resp., second, third) term contains a difference of $G$-terms (resp., $\si$-terms, both).
For example, 
\begin{align*}
F^{n,i}_{1,1}&=2\iiint_0^{\la_n\Del} \iint_{0}^{(i-\la_n)\Del} \Bigl[G_{x_{k_{12}}-y-z}(i\Del-r)-G_{x_{k_{12}}-z}(i\Del-r)\Bigr]\si'(u(i\Del,x_{k_{12}})_{0}^r)\\
&\quad\times\si(u(i\Del,x_{k_{12}})_{0}^r)  \si(u(r,z))\,W(\dd r,\dd z)\,\Pi^n_{|l_1-l_2|,x_{k_1}-x_{k_2}}(\dd s,\dd y,\dd y'),\\
F^{n,i}_{1,2} &= 2\iiint_0^{\la_n\Del} \iint_{0}^{(i-\la_n)\Del} G_{x_{k_{12}}-z}(i\Del-r) \Bigl[\si'(u(i\Del,x_{k_{12}}-y)_{0}^r)\si(u(i\Del,x_{k_{12}}-y')_{0}^r)\\
&\quad  -\si'(u(i\Del,x_{k_{12}})_{0}^r)\si(u(i\Del,x_{k_{12}})_{0}^r)\Bigr] \si(u(r,z))\,W(\dd r,\dd z)\,\Pi^n_{|l_1-l_2|,x_{k_1}-x_{k_2}}(\dd s,\dd y,\dd y'),\\
F^{n,i}_{1,3} &=2\iiint_0^{\la_n\Del} \iint_{0}^{(i-\la_n)\Del} \Bigl(G_{x_{k_{12}}-y-z}(i\Del-r)-G_{x_{k_{12}}-z}(i\Del-r)\Bigr)\\
&\quad\times\Bigl(\si'(u(i\Del,x_{k_{12}}-y)_{0}^r)\si(u(i\Del,x_{k_{12}}-y')_{0}^r) -\si'(u(i\Del,x_{k_{12}})_{0}^r)\si(u(i\Del,x_{k_{12}})_{0}^r)\Bigr)\\
&\quad \times \si(u(r,z))\,W(\dd r,\dd z)\,\Pi^n_{|l_1-l_2|,x_{k_1}-x_{k_2}}(\dd s,\dd y,\dd y').
\end{align*}
Similarly to Lemma~\ref{lem:remove-s}, the lemma is proved  once we can show that each term $F^{n,i}_{j,j'}$, $j,j'=1,2,3$, is of size $o(\Delh)$ uniformly in $i$.

For $F^{n,i}_{1,1}$, we apply Taylor's theorem with integral remainder to the $G$-difference and obtain
\begin{align*}
F^{n,i}_{1,1}&=2\sum_{j=1}^d \iiint_0^{\la_n \Del}  -y_j  \iint_0^{(i-\la_n)\Del} \partial_{x_j} G (i\Del-r,x_{k_{12}}-z)\si'(u(i\Del,x_{k_{12}})_0^r)\\
&\quad\quad\times\si(u(i\Del,x_{k_{12}})_0^r)\si(u(r,z))\,W(\dd r,\dd z)\,\Pi^n_{|l_1-l_2|,x_{k_1}-x_{k_2}}(\dd s,\dd y,\dd y')\\
&\quad+2\sum_{j,j'=1}^d \iiint_0^{\la_n\Del} y_jy_{j'} \iint_0^{(i-\la_n)\Del}  \int_0^1 (1-u) \partial^2_{x_jx_{j'}} G(i\Del-r,x_{k_{12}}-uy-z)\,\dd u\\
&\qquad\times \si'(u(i\Del,x_{k_{12}})_0^r)\si(u(i\Del,x_{k_{12}})_0^r)\si(u(r,z))\,W(\dd r,\dd z)\,\Pi^n_{|l_1-l_2|,x_{k_1}-x_{k_2}}(\dd s,\dd y,\dd y').
\end{align*}
The first part is identically zero because the integral of $y_j$ with respect to $\Pi^n$ vanishes for symmetry reasons. Let us now determine the size of the second part, which we denote by $\wh F^{n,i}_{1,1}$ in the following. To this end, we fix $j$ and $j'$ and permute the $W$-integral with the $\dd u$-integral. Leaving aside terms that do not depend on $r$ and $z$ for a moment and using the fact that all $\si$-terms are of order $1$, the size of the $W$-integral is bounded by a constant times
\begin{align*}
&\iiint_0^{(i-\la_n)\Del} \Bigl|\partial^2_{x_jx_{j'}} G (i\Del-r,x_{k_{12}}-uy-z)\partial^2_{x_jx_{j'}} G(i\Del-r,x_{k_{12}}-uy-z')\Bigr|\,\La(\dd z,\dd z')\,\dd r\\
&\quad\leq\iiint_{\la_n\Del}^T \Bigl|\partial^2_{x_jx_{j'}} G  (r,z)\partial^2_{x_jx_{j'}} G (r,z')\Bigr|\,\La(\dd z,\dd z')\,\dd r.
\end{align*}
Since $\ee^{-\frac{w}{4}}w \leq \frac{4}{\ee}$, we have
\beq\label{eq:sp-der-G}\begin{split} \Bigl|\partial^2_{x_jx_{j'}} G(t,x)\Bigr| 
	&=\begin{cases} G(t,x)\Bigl|\frac{x^2_j}{t^2}-\frac1{t}\Bigr| \leq \frac{G(t,\frac{x}{\sqrt 2})}{t}\Bigl(\ee^{-\frac{|x|^2}{4 t}}\frac{|x|^2}{t}+\ee^{-\frac{|x|^2}{4 t}}\Bigr) &\text{if } j=j',\\ G(t,x)\frac{|x_jx_{j'}|}{t^2}\leq \frac{G(t,\frac{x}{\sqrt 2})}{ t}\ee^{-\frac{|x|^2}{4 t}}\frac{|x|^2}{t}&\text{if } j\neq j'\end{cases} \\
	&\leq\frac{\frac4 \ee +1}{t}G(t,\textstyle\frac{x}{\sqrt 2}).\end{split}\eeq
Thus, if we switch to the Fourier domain and use the identity
\beq\label{eq:FG} \calf G(t,\cdot)(\xi) = \ee^{-2\pi^2 t |\xi|^2},\eeq 
we derive, omitting multiplicative factors in the exponential by a scaling argument,
\beq\label{eq:2der-prod}\begin{split}
	&\iiint_{\la_n\Del}^T \Bigl|\partial^2_{x_jx_{j'}} G (r,z)\partial^2_{x_jx_{j'}} G (r,z')\Bigr|\,\La(\dd z,\dd z')\,\dd r\lec \int_{\la_n\Del}^T \frac{1}{r^2}\int_{\R^d} \ee^{-r |\xi|^2}\,\mu(\dd \xi)\,\dd r\\
	&\quad\lec \int_{\la_n\Del}^T \frac{1}{r^2}\int_0^\infty \ee^{-rw^2}w^{\al-1}\,\dd w\,\dd r \lec \int_{\la_n\Del}^T \frac{1}{r^2}r^{-\frac12-\frac{\al-1}{2}}\,\dd r \lec (\la_n\Del)^{-1-\frac{\al}{2}},
\end{split}\eeq
where $\mu$ is the spectral measure from \eqref{eq:isometry}.
So by a size estimate, we get from \eqref{eq:y2},
\beq\label{eq:whF-1}  \E[|\wh F^{n,i}_{1,1}|] \lec (\la_n\Del)^{-\frac12-\frac\al4} \iiint ({|y|^2}+{|y'|}^2)\,|\Pi^n_{|l_1-l_2|,x_{k_1}-x_{k_2}}|  (\dd s,\dd y,\dd y')\lec \Del^{-\frac\Theta 2(1-a)+1},\eeq
which is $o(\Delh)$ because $\Theta\leq\frac32$ and $a>\frac1{2\Theta}$ and thus, $-\frac\Theta 2(1-a)+1>1-\frac\Theta 2(1-\frac1{2\Theta}) =  \frac54 -\frac\Theta2\geq\frac12$.

For $F^{n,i}_{1,2}$, we have to develop the $\si$-difference further using Taylor's theorem. This results in
\beq\label{eq:Fn2}\begin{split}
	F^{n,i}_{1,2}&=2\iiint_0^{\la_n\Del}\iint_0^{(i-\la_n)\Del} G_{x_{k_{12}}-z}(i\Del-r) \biggl[ \si''(u(i\Del,x_{k_{12}})_0^r)\si(u(i\Del,x_{k_{12}})_0^r)\\
	&\quad\quad\times\Bigl(u(i\Del,x_{k_{12}}-y)_0^r -u(i\Del,x_{k_{12}})_0^r\Bigr) +(\si')^2(u(i\Del,x_{k_{12}})_0^r)\\
	&\quad\quad\times\Bigl(u(i\Del,x_{k_{12}}-y')_0^r -u(i\Del,x_{k_{12}})_0^r\Bigr)\biggr]\si(u(r,z))\,W(\dd r,\dd z)\,\Pi^n_{|l_1-l_2|,x_{k_1}-x_{k_2}}(\dd s,\dd y,\dd y')\\
	&\quad+\iiint_0^{\la_n\Del}\iint_0^{(i-\la_n)\Del} G_{x_{k_{12}}-z}(i\Del-r) \biggl[ \si'''(\zeta^{n,i}_1(r,y))\si(\zeta^{n,i}_2(r,y'))\\
	&\quad\quad\quad\times\Bigl(u(i\Del,x_{k_{12}}-y)_0^r -u(i\Del,x_{k_{12}})_0^r\Bigr)^2+ \si'(\zeta^{n,i}_1(r,y))\si''(\zeta^{n,i}_2(r,y'))\\
	&\quad\quad\quad\times\Bigl(u(i\Del,x_{k_{12}}-y')_0^r -u(i\Del,x_{k_{12}})_0^r\Bigr)^2+ 2\si''(\zeta^{n,i}_1(r,y))\si'(\zeta^{n,i}_2(r,y'))\\
	&\quad\quad\quad\times\Bigl(u(i\Del,x_{k_{12}}-y)_0^r -u(i\Del,x_{k_{12}})_0^r\Bigr)\Bigl(u(i\Del,x_{k_{12}}-y')_0^r -u(i\Del,x_{k_{12}})_0^r\Bigr)\biggr]\\
	&\quad\quad\times\si(u(r,z))\,W(\dd r,\dd z)\,\Pi^n_{|l_1-l_2|,x_{k_1}-x_{k_2}}(\dd s,\dd y,\dd y')
\end{split}\raisetag{6\baselineskip}\eeq
for some intermediate values $\zeta^{n,i}_1(r,y)$ and $\zeta^{n,i}_2(r,y')$. The integrand of the first $\Pi^n$-integral depends on $y^{(\prime)}$ only through the two $u$-differences, which can be written as
\beq\label{eq:udiff}\begin{split} &u(i\Del,x_{k_{12}}-y^{(\prime)})_0^r -u(i\Del,x_{k_{12}})_0^r\\
	&\quad= \iint_0^r \Bigl(G_{x_{k_{12}}-y^{(\prime)}-w}(i\Del-v)-G_{x_{k_{12}}-w}(i\Del-v)\Bigl)\si(u(v,w))\,W(\dd v,\dd w). \end{split}\eeq
Therefore, the first $\Pi^n$-integral in \eqref{eq:Fn2} has a very similar structure to $F^{n,i}_{1,1}$. In fact, the arguments from above can be applied analogously to show that this part is asymptotically negligible. 

For the second $\Pi^n$-integral in \eqref{eq:Fn2}, a size estimate suffices. Indeed, the order of magnitude of the integral in \eqref{eq:udiff} is 
\beq\label{eq:help}
|y|\sum_{j=1}^d\Biggl(\iiint_0^r \Bigl|\partial_{x_j}G(i\Del-v,w)\partial_{x_j}G(i\Del-v,w')\Bigr|\,\La(\dd w,\dd w')\,\dd v\Biggr)^{\frac12}.
\eeq
Similarly to \eqref{eq:sp-der-G} and  \eqref{eq:2der-prod}, we have
\beq\label{eq:g-firstder} \Bigl|\partial_{x_j}G(t,x)\Bigr| 
= G(t,x)\frac{|x_j|}{t}\leq \frac{G(t,\frac{x}{\sqrt 2})}{\sqrt{t}}\ee^{-\frac{|x|^2}{4 t}}\frac{|x|}{\sqrt{t}} \leq \sqrt{\frac{2}{\ee t}}G(t,\textstyle\frac x{\sqrt 2}) \eeq
and, because $r\leq (i-\la_n)\Del$,
\begin{align*}
&\iiint_0^r \Bigl|\partial_{x_j}G(i\Del-v,w)\partial_{x_j}G(i\Del-v,w')\Bigr|\,\La(\dd w,\dd w')\,\dd v\\
&\quad\lec \int_{\la_n\Del}^T v^{-1}\int_0^\infty \ee^{-vw^2}w^{\al-1}\,\dd w\,\dd v\lec \int_{\la_n\Del}^T   v^{-\frac\al2-1}\lec (\la_n\Del)^{-\frac\al2}.
\end{align*}
As each of the three summands within the second pair of brackets in  \eqref{eq:Fn2} contains two $u$-differences
and integration of $|y|^2$ against $|\Pi^n|$ yields a factor $\Del$ by \eqref{eq:y2}, the second $\Pi^n$-integral in \eqref{eq:Fn2} is of size $\Del^{-(1-a)\al/2+1}$, which is $o(\Delh)$ because $1-(1-a)\frac\al 2\geq1-\frac12(1-a) = \frac {a+1} 2>\frac12$.

For $F^{n,i}_{1,3}$, the reasoning is similar since the $G$-difference is of size $|y|$ times the first spatial derivatives of $G$ [cf.\ \eqref{eq:help}] and the difference of the $\si$-terms, after linearizing, becomes differences of $u$-terms, which are also of size $|y^{(\prime)}|$ times the first spatial derivatives of $G$ [cf.\ \eqref{eq:Fn2} and \eqref{eq:udiff}].

A careful inspection of the remaining terms reveals that  both $F^{n,i}_{2,j'}$ and $F^{n,i}_{3,j'}$ can be analyzed similarly to $F^{n,i}_{1,j'}$, for all $j'=1,2,3$. The only notable differences are first, that instead of a $W(\dd r,\dd z)$-integral, they have integrals with respect to $\La(\dd z,\dd z')\,\dd r$, and second, that there are products of two $G$-terms instead of single ones. The former does not affect the reasoning at all. Concerning the latter, we can simply use \eqref{eq:id3} to turn differences of products of $G$-terms into differences of single $G$-terms. Notice that we need up to four derivatives of $\si$ for this argument.
\epr

\bpr[Proof of Lemma~\ref{lem:remove-Pi}] Recall that we may assume $M=K=1$ according to Remark~\ref{rem:2}. By the mean value theorem (recall that $\un\mu_f$ is differentiable by \cite[(D.46)]{Chong19a}), a standard size estimate, and the relation \eqref{eq:muid}, we are therefore left to show that [with $\nu_i=(1,l_i)$ for $i=1,2$]
\beq\label{eq:diag} \E\Bigl[\Bigl|v^{\prime\prime n,i}_{\nu_1,\nu_2}-\Ga_{|l_1-l_2|}\si^2(u(i\Del,x_1)_0^{(i-\la_n)\Del})\Bigr|\Bigr] = o(\Delh),  \eeq
uniformly in $i$. The left-hand side of the previous line is bounded by a constant times 
\begin{align*}
&\E[|\si^2(u(i\Del,x_1)_0^{(i-\la_n)\Del})|]|\Pi^n_{|l_1-l_2|,0}|((\la_n\Del,\infty)\times\R^d\times\R^d)\\
&\quad+\E[|\si^2(u(i\Del,x_1)_0^{(i-\la_n)\Del})|]\bigl|\Pi^n_{|l_1-l_2|,0}([0,\infty)\times\R^d\times\R^d)-\Ga_{|l_1-l_2|}\bigr|, 
\end{align*}
which is $\lec\Del^{\Theta a}+0$ by \eqref{eq:Theta} and \eqref{eq:PinGa}. So \eqref{eq:diag} follows from the hypothesis that $a>\frac1{2\Theta}$.
\epr

\bpr[Proof of Lemma~\ref{lem:iDels}]
The left-hand side of \eqref{eq:help2} is equal to
\begin{align*} &\frac{1}{\sqrt{\Del}}\sumtla \int_{(i-1)\Del}^{i\Del} \Biggl(  \mu_f\Bigl(\si^2(u(i\Del,\un x)_0^{(i-\la_n)\Del})\Bigr)-  \mu_f\Bigl(\si^2(u(s,\un x)_0^{s-\la_n\Del})\Bigr)  \Biggr)\,\dd s\\
&\quad 
-\frac{1}{\sqrt{\Del}}\int_{[t/\Del]\Del}^t \mu_f\Bigl(\si^2(u(s,\un x)_0^{s-\la_n\Del})\Bigr)\,\dd s.  \end{align*}
The last line is of order $\Del^{-1/2}\Del=\Del^{1/2}$ and therefore negligible. So only the term in the first line needs to be considered further, which  will be denoted by $G^n(t)$ from now on. Letting $\Phi\colon \R^K\to\R^M$ be the function that maps $x=(x_1,\dots,x_K)\in\R^K$ to $\mu_f((\si^2(x_k))_{k=1,\dots,K})$,  we have  by Taylor's theorem ($\mu_f$ is twice differentiable by \cite[(D.45)]{Chong19a}), 
\beq\label{eq:Gn}\begin{split}
	G^n(t)&=\frac{1}{\sqrt{\Del}}\sumtla \sum_{k=1}^K \int_{(i-1)\Del}^{i\Del} \partial_k\Phi(u(s,\un x)_0^{s-\la_n\Del})\Bigl(u(i\Del,x_k)_0^{(i-\la_n)\Del}-u(s,x_k)_0^{s-\la_n\Del}\Bigr)\,\dd s\\
	&\quad+\frac{1}{2\sqrt{\Del}}\sumtla  \sum_{k_1,k_2=1}^K \int_{(i-1)\Del}^{i\Del} \partial^2_{k_1k_2}\Phi(\eta^n_i(s))\Bigl(u(i\Del,x_{k_1})_0^{(i-\la_n)\Del}-u(s,x_{k_1})_0^{s-\la_n\Del}\Bigr)\\
	&\qquad\times\Bigl(u(i\Del,x_{k_2})_0^{(i-\la_n)\Del}-u(s,x_{k_2})_0^{s-\la_n\Del}\Bigr) \,\dd s\\
	& = G^n_1(t)+G^n_2(t)+G^n_3(t),
\end{split}\raisetag{2\baselineskip}\eeq
where $\eta^n_i(s)$ is some intermediate value and 
\begin{align*}
G^n_1(t)&:=\frac{1}{\sqrt{\Del}}\sumtla \sum_{k=1}^K\int_{(i-1)\Del}^{i\Del} \partial_k\Phi(u(s,\un x)_0^{s-\la_n\Del})\\
&\mathrel{\hphantom{:=} \mathrel{\times}}\iint_0^{s-\la_n\Del} (G_{x_k-z}(i\Del-r) - G_{x_k-z}(s-r))\si(u(r,z))\,W(\dd r,\dd z)\,\dd s,\\
G^n_2(t)&:=\frac{1}{\sqrt{\Del}}\sumtla\sum_{k=1}^K \int_{(i-1)\Del}^{i\Del} \partial_k\Phi(u(s,\un x)_0^{s-\la_n\Del})\\
&\mathrel{\hphantom{:=} \mathrel{\times}} \iint_{s-\la_n\Del}^{(i-\la_n)\Del} G_{x_k-z}(i\Del-r)  \si(u(r,z))\,W(\dd r,\dd z)\,\dd s,
\end{align*}
and $G^n_3(t)$ is the expression spanning over the second and third line of \eqref{eq:Gn}.

Observing that the $W$-integral in the definition of $G^n_1$ is of order 
\beq\label{eq:sizeG} \Biggl(\iint_{\la_n\Del}^T |G_{z}(r+i\Del-s)-G_{z}(r)||G_{z'}(r+i\Del-s)-G_{z'}(r)|\,\La(\dd z,\dd z')\,\dd r\Biggr)^{\frac12}\lec \Bigl(\Del^{1-\frac\al2+\Theta a}\Bigr)^{\frac12}\eeq
by \eqref{eq:Theta} and the fact that ${|i\Del-s|}\leq \Del$, we derive
\[ \E[(G^n_1)^\ast_T]\lec \Del^{-\frac12+\frac12-\frac\al4+\frac\Theta2 a} = \Del^{\frac\Theta2 a-\frac\al4} \to0,\]
because $\frac\Theta2 a-\frac\al4>\frac\Theta 2 \times \frac1{2\Theta}-\frac14=0$ by hypothesis.

For $G^n_2$, let us first interchange the $W$- and the $\dd s$-integrals, which leads to
\beq\label{eq:Gn2}\begin{split}
	G^n_2(t)&=\frac{1}{\sqrt{\Del}}\sumtla \sum_{k=1}^K \iint_{(i-1-\la_n)\Del}^{(i-\la_n)\Del}  \Biggl(\int_{(i-1)\Del}^{r+\la_n\Del} \partial_k\Phi(u(s,\un x)_0^{s-\la_n\Del})\,\dd s\Biggr)\\
	&\quad\times G_{x_k-z}(i\Del-r)  \si(u(r,z))\,W(\dd r,\dd z).
\end{split}\eeq
Now the $\dd s$-integral is of size $\Del$ at most, while the $W$-integral (without the former) is of size
\beq\label{eq:help4}\begin{split} &\Biggl(\iiint_{(i-1-\la_n)\Del}^{(i-\la_n)\Del} G_{x_k-z}(i\Del-r)G_{x_k-z'}(i\Del-r)\,\La(\dd z,\dd z')\,\dd r\Biggr)^{\frac12} \\
	&\quad = \Biggl(\iiint_{\la_n\Del}^{(\la_n+1)\Del} G_{z}(r)G_{z'}(r)\,\La(\dd z,\dd z')\,\dd r\Biggr)^{\frac12} \lec \Bigl((\la_n+1)\Del)^{1-\frac\al2}-(\la_n\Del)^{1-\frac\al2}\Bigr)^{\frac12}\\
	&\quad= (\la_n\Del)^{\frac12-\frac\al4}((1+\la_n^{-1})^{1-\frac\al2}-1)^{\frac12} \lec \la_n^{-\frac12}(\la_n\Del)^{\frac12-\frac\al4}.
\end{split}\eeq
Moreover, the $i$th summand in \eqref{eq:Gn2} is $\calf^n_{i-\la_n}$-measurable with vanishing  $\calf^n_{i-1-\la_n}$-conditional  expectation. Thus, a martingale size estimate gives us
\[ \E[(G^n_2)^\ast_T]\lec \Del^{-\frac12}\Del^{-\frac12}\Del\Del^{\frac a 2+(1-a)(\frac12-\frac\al4)}\to0. \]

Finally, since $u(i\Del,x_k)_0^{(i-\la_n)\Del}-u(s,x_k)_0^{s-\la_n\Del}$ is of size  $(\Del^{1-\frac\al2+\Theta a})^{\frac12}+\la_n^{-\frac12}(\la_n\Del)^{\frac12-\frac\al4}$ by \eqref{eq:sizeG} and \eqref{eq:help4}, a standard size estimate and the fact that $1-\frac\al2\geq \frac12$ imply
\[ \E[(G^n_3)^\ast_T] \lec \Del^{-\frac12}\Bigl(\Del^{1-\frac\al2+\Theta a}+\Del^{a+(1-a)(1-\frac\al2)}\Bigr) \leq \Del^{\Theta a}+\Del^{\frac{a}{2}}\to0, \]
which completes the proof of the lemma.
\epr

\bpr[Proof of Lemma~\ref{lem:Bn1}] By Taylor's theorem, we have 
\beq\label{eq:Bn1-Bn12} B^{n,m}_2(t)-B^{n,m,1}_2(t)=\Delh\sumtla \sum_{\nu\in\cali} (\kappa^{n,m,i}_\nu-\E[\kappa^{n,m,i}_\nu \mid \calf^n_{i-\la_n}]),\eeq 
where
\[ \kappa^{n,m,i}_\nu:=\frac12\sum_{\nu'\in\cali}\partial^2_{\nu\nu'} f(\wt\kappa^{n,m}_i)  (\ov\beta^n_i-\ov\beta^n_{i,k^{n,m}_i})_\nu(\ov\beta^n_i-\ov\beta^n_{i,k^{n,m}_i})_{\nu'} \]
for some intermediate point $\wt\kappa^{n,m}_i$. By \eqref{eq:increments}, the definition of the $\beta$-terms in \eqref{eq:betaga}, and the fact that $k^{n,m}_i=O(\la_n)$ for fixed $m$, we deduce that $\kappa^{n,m,i}_\nu$ is of size $((\la_n\Del)^{1/2-\al/4})^2$. Moreover, the $i$th term in \eqref{eq:Bn1-Bn12} has zero $\calf^n_{i-\la_n}$-conditional expectation. Thus, a martingale size estimate gives 
\beq\label{eq:comparison}\E[(B^{n,m}_2-B^{n,m,1}_2)^\ast_T]\lec \la_n^{\frac12}(\la_n\Del)^{1-\frac\al2} \leq \la_n^{\frac12}(\la_n\Del)^{\frac12} \leq \Del^{\frac12-a},\eeq
which goes to $0$ if $a$ is close to $\frac{1}{2\Theta}$.
\epr

\bpr[Proof of Lemma~\ref{lem:Bn2}] The easiest term is $B^{n,m,1,4}_2$: the $W(\dd r,\dd z)$-integral is of size $(\la_n\Del)^{1/2-\al/4}$, so a martingale size estimate immediately gives $\E[(B^{n,m,1,4}_2)^\ast_T]\lec \la_n^{1/2}(\la_n\Del)^{1-\al/2}\to0$ by \eqref{eq:comparison}.

Next, observe that the term $\theta^{n,m,i}_{1,3}$ is of size $(\ov\la_n\Del)^{1/2-\al/4}$ [due to the $W(\dd r,\dd z)$-integral]. So if we replace $\si(u(r,z))$ by $\si(u((i-\la_n)\Del,z))$ in $\theta^{n,m,i}_{1,3}$, the difference resulting from this modification will come with  an additional factor of $(\la_n\Del)^{1/2-\al/4}$. Hence, if we denote the so-obtained term by $\wt\theta^{n,m,i}_{1,3}$ and the corresponding approximation of $B^{n,m,1,3}_2$ by $\wt B^{n,m,1,3}_2$, a martingale size estimate yields $$\E[(B^{n,m,1,3}_2-\wt B^{n,m,1,3}_2)^\ast_T]\lec \la_n^{\frac12}(\ov\la_n\Del)^{\frac12-\frac\al4}(\la_n\Del)^{\frac12-\frac\al4},$$ which converges to $0$ as $n\to\infty$ by comparison with \eqref{eq:comparison}.

The crucial observation is now that the $\calf^n_{i-\la_n}$-conditional law of $\theta^{n,m,i}_{1,1}$ (resp., $\wt\theta^{n,m,i}_{1,3}$) agrees with its $\calf^n_{i-\ov\la_n}$-conditional law. As a consequence, subtracting its $\calf^n_{i-\la_n}$-conditional expectation in $B^{n,m,1,1}_2$ (resp., $\wt B^{n,m,1,3}_2$) amounts to subtracting its $\calf^n_{i-\ov\la_n}$-conditional expectation. This insight allows us to improve our martingale size estimate and derive 
\beq\label{eq:ov} \E[(B^{n,m,1,1}_2)^\ast_T]+\E[(\wt B^{n,m,1,3}_2)^\ast_T]\lec \ov\la_n^{\frac12}(\la_n\Del)^{\frac12-\frac\al4}= \Del^{-\frac{\ov a}{2}+(1-a)(\frac12-\frac\al4)}.\eeq 
By choosing $a$ and $\ov a$ sufficiently small, the exponent gets arbitrarily close to $$-\frac{1}{(2\Theta)^2}+\biggl(1-\frac1{2\Theta}\biggr)\biggl(\frac12 -\frac\al 4\biggr)=-\frac1{(2+\al)^2}+\frac{(1+\al)(2-\al)}{4(2+\al)} = \frac{\al(4-\al-\al^2)}{4(2+\al)^2},$$
which is strictly positive for all $0<\al\leq1$.

It remains to prove that $B^{n,m,1,2}_2-B^{n,2}_2\limL0$ for every $m\in\N$. To this end, notice that $k^{n,m}_i$ appears twice in the definition of $(\theta^{n,m,i}_{1,2})_\nu$: once in the $\partial_\nu f$-term and once in the $\si'$-term. Recalling that $k^{n,m}_i\leq (m+1)\la_n+L-1 = O(\la_n)$ for fixed $m\in\N$, we see that replacing $k^{n,m}_i$ by $0$ leads to an error of size $(\la_n\Del)^{1/2-\al/4}$ for both terms and therefore, by \eqref{eq:id2}, also for the difference $(\theta^{n,m,i}_{1,2})_\nu-(\theta^{n,i}_2)_\nu$. As the $W(\dd r,\dd z)$-integral yields another factor of $(\la_n\Del)^{1/2-\al/4}$, a martingale size estimate shows that $\E[(B^{n,m,1,2}_2-B^{n,2}_2)^\ast_T]\lec \la_n^{1/2}(\la_n\Del)^{1-\al/2}\to0$ by \eqref{eq:comparison}.
\epr

\bpr[Proof of Lemma~\ref{lem:Bn3}] As a first step, notice that $(\wt \theta^{n,i}_3)_\nu$ arises from $(\theta^{n,i}_2)_\nu$ upon replacing $\ov\beta^n_{i,0}=\Delta^n_i \wt u^{(i+L-1)\Del}_{(i-\ov\la_n)\Del,(i-\la_n)\Del}/\tau_n$  [resp., $\si'(u((i-\la_n)\Del, y))$, $\si(u(r,z))$] in the definition of $(\theta^{n,i}_2)_\nu$ by $\Delta^n_i  u^{(i+L-1)\Del}_{(i-\ov\la_n)\Del,(i-\la_n)\Del}/\tau_n$ [resp., $\si'(u(s,y)_0^{(i-\la_n)\Del})$, $\si(u(r,z)_0^{(i-\la_n)\Del})$]. Each change leads to a factor of order $(\la_n\Del)^{1/2-\al/4}$, in addition to the $(\la_n\Del)^{1/2-\al/4}$-factor coming from the $W(\dd r,\dd z)$-integral in \eqref{eq:wttheta3}. Thus, $B^{n,2}_2-\wt B^{n,3}_2 \limL 0$ by \eqref{eq:comparison}, where 
$$\wt B^{n,3}_2(t):=\Delh\sumtla \sum_{\nu\in\cali} ((\wt \theta^{n,i}_3)_\nu-\E[(\wt \theta^{n,i}_3)_\nu\mid \calf^n_{i-\la_n}]).$$

Observe now that we actually have $\E[(\wt\theta^{n,i}_3)_\nu \mid \calf^n_{i-\la_n}] = 0$. This is because conditionally on $\calf^n_{i-\la_n}$, the iterated $W$-integral in \eqref{eq:wttheta3} is an element of the second Wiener chaos, while the $\partial_\nu f$-term belongs to the direct sum of odd-order Wiener chaoses.
As a consequence, we have $\wt B^{n,3}_2(t)=\Delh\sumtla \sum_{\nu\in\cali} (\wt \theta^{n,i}_3)_\nu$. So if we subtract  from each $(\wt \theta^{n,i}_3)_\nu$ the term $(\theta^{n,i}_3)_\nu$ (which, by definition, is the $\calf^n_{i-\ov\la_n}$-conditional expectation  of the former), then, as we have seen in \eqref{eq:ov}, a martingale size estimate yields $\E[(\wt B^{n,3}_2-B^{n,3}_2)^\ast_T]\lec \ov\la_n^{1/2}(\la_n\Del)^{1/2-\al/4}\to0$.
\epr

\bpr[Proof of Lemma~\ref{lem:Bn345}] Because ${|s-i\Del|}\leq \ov\la_n\Del$, we know from \eqref{eq:Pin1} (with $|s-i\Del|$ for the value of $\Del$) that
\beq\label{eq:Gest} \iint_0^{i\Del} (G_{y-z}(i\Del-r)-G_{y-z}(s-r))\si(u(r,z)_0^{(i-\la_n)\Del})\,W(\dd r,\dd z) \eeq
is of size at most $(\ov\la_n\Del)^{1/2-\al/4}$. Together with \eqref{eq:siest} and \eqref{eq:id2}, we have that $\wt \theta^{n,i}_3-\wt\theta^{n,i}_4$ is of the same size $(\ov\la_n\Del)^{1/2-\al/4}$. Moreover, by the same reason as in the proof of Lemma~\ref{lem:Bn3}, we have \beq\label{eq:condexp0} \E[(\theta^{n,i}_{3\mid 4\mid 5})_{\nu}\mid \calf^n_{i-\la_n}] = \E[(\wt \theta^{n,i}_{3\mid 4\mid 5})_{\nu}\mid \calf^n_{i-\la_n}] = 0.\eeq
Hence, when $a$ and $\ov a$ are sufficiently small, a martingale size estimate yields 
\beq\label{eq:help6} \E[(B^{n,3}_2-B^{n,4}_2)^\ast_T]\lec \la_n^{\frac12}(\ov\la_n\Del)^{\frac12-\frac\al4} \leq \Del^{-\frac a2+(1-\ov a)(\frac12-\frac\al4)} \to 0. \eeq
Indeed, recalling that $\Theta=1+\frac\al2$, the exponent $-\frac a2 +(1-\ov a)(\frac12-\frac\al4)$ is then close to 
\[ -\frac{1}{4\Theta}+\biggl(1-\frac1{2\Theta^2}\biggr)\frac{2-\al}{4} = \frac{-2\Theta+(2-\al)(2\Theta^2-1)}{8\Theta^2}= \frac{\al(4-\al^2-2\al)}{16\Theta^2},\]
which is strictly positive for $\al\in(0,1]$.

For $B^{n,4}_2-B^{n,5}_2$, observe that $\ov v^{\prime n,i}_{\nu_1,\nu_2}-\ov v^{\prime\prime n,i}_{\nu_1,\nu_2}$ is exactly the same as $v^{\prime n,i}_{\nu_1,\nu_2}-v^{\prime\prime n,i}_{\nu_1,\nu_2}$  (defined before Lemma~\ref{lem:remove-s}) except that the domain of integration of the $\Pi^n$-integral goes up only to $(\ov\la_n+l_{12}-1)\Del$ instead of $\la_n\Del$. This difference has no impact on the estimates that have been carried out in \eqref{eq:delta-split} and thereafter. We conclude that $\ov v^{\prime n,i} -\ov v^{\prime\prime n,i}$ is of size $o(\Delh)$, uniformly in $i$. 

Next, using \eqref{eq:id3} on the $\si'$-terms and the $W$-integrals in $\ov w^{\prime n,i}_{\nu_1,\nu_2}$ and $\ov w^{\prime\prime n,i}_{\nu_1,\nu_2}$ (hereby grouping a $W$-integral and a $\si'$-term together as one factor) and recalling that  $\Pi^n$ is symmetric in interchanging $\dd y$ and $\dd y'$, we obtain
\begin{align*}
&\ov w^{\prime n,i}_{\nu_1,\nu_2}-\ov w^{\prime\prime n,i}_{\nu_1,\nu_2}\\
&\quad=2\iiint \Biggl(\iint_{(i-\la_n)\Del}^{(i-\ov\la_n)\Del}  \Bigl[ G_{x_{k_{12}}-y-z}(i\Del-r)\si'(u(i\Del,x_{k_{12}}-y)_0^{(i-\la_n)\Del})\\
&\qquad\qquad-G_{x_{k_{12}}-z}(i\Del-r)\si'(u(i\Del,x_{k_{12}})_0^{(i-\la_n)\Del})\Bigr]  \si(u(r,z)_0^{(i-\la_n)\Del})\,W(\dd r,\dd z)\Biggr)\\
&\quad\qquad\times \si'(u(i\Del,x_{k_{12}})_0^{(i-\la_n)\Del})\iint_{(i-\la_n)\Del}^{(i-\ov\la_n)\Del}  G_{x_{k_{12}}-z}(i\Del-r)  \si(u(r,z)_0^{(i-\la_n)\Del})\,W(\dd r,\dd z)\\
&\quad\qquad\times\bone_{s\in(0,(\ov\la_n+l_{12}-1)\Del]}\,\Pi^n_{|l_1-l_2|,x_{k_1}-x_{k_2}}(\dd s,\dd y,\dd y')\\
&\quad\quad+\iiint \Biggl(\iint_{(i-\la_n)\Del}^{(i-\ov\la_n)\Del}  \Bigl[ G_{x_{k_{12}}-y-z}(i\Del-r)\si'(u(i\Del,x_{k_{12}}-y)_0^{(i-\la_n)\Del})\\
&\qquad\qquad-G_{x_{k_{12}}-z}(i\Del-r)\si'(u(i\Del,x_{k_{12}})_0^{(i-\la_n)\Del})\Bigr]  \si(u(r,z)_0^{(i-\la_n)\Del})\,W(\dd r,\dd z)\Biggr)\\
&\quad\qquad\times \Biggl(\iint_{(i-\la_n)\Del}^{(i-\ov\la_n)\Del}  \Bigl[ G_{x_{k_{12}}-y'-z}(i\Del-r)\si'(u(i\Del,x_{k_{12}}-y')_0^{(i-\la_n)\Del})\\
&\qquad\qquad-G_{x_{k_{12}}-z}(i\Del-r)\si'(u(i\Del,x_{k_{12}})_0^{(i-\la_n)\Del})\Bigr]  \si(u(r,z)_0^{(i-\la_n)\Del})\,W(\dd r,\dd z)\Biggr)\\
&\quad\qquad\times\bone_{s\in(0,(\ov\la_n+l_{12}-1)\Del]}\,\Pi^n_{|l_1-l_2|,x_{k_1}-x_{k_2}}(\dd s,\dd y,\dd y').
\end{align*}
Inside the first $\Pi^n$-integral, there is a $W$-integral involving, in brackets, a difference of $G$- and $\si'$-terms. There is also another expression (the third $\si'$-term and the second $W$-integral) that does not depend on the variables $s$, $y$, and $y'$. We can take this out of the $\Pi^n$-integral, which then has the same structure as the term $F^{n,i}_1$ in the proof of Lemma~\ref{lem:remove-y}. With virtually the same arguments that have led to \eqref{eq:whF-1}, we deduce that  the first $\Pi^n$-integral above is of magnitude $(\ov\la_n\Del)^{-1/2-\al/4}\Del$. [We have $\ov\la_n$ now because it appears in the upper limit of the first $W$-integral above. Moreover, the second one, which we have taken outside the $\Pi^n$-integral, is actually of size $(\la_n\Del)^{1/2-\al/4}$ but it suffices in the following to bound this simply by $1$.] 

Concerning the second $\Pi^n$-integral above, notice that the two $W$-integrals equal 
\begin{align*} &\si'(u(i\Del,x_{k_{12}}-y^{(\prime)})_0^{(i-\la_n)\Del}) u(i\Del,x_{k_{12}}-y^{(\prime)})_{(i-\la_n)\Del,(i-\la_n)\Del}^{(i-\ov \la_n)\Del}\\
&\quad-\si'(u(i\Del,x_{k_{12}})_0^{(i-\la_n)\Del}) u(i\Del,x_{k_{12}})_{(i-\la_n)\Del,(i-\la_n)\Del}^{(i-\ov \la_n)\Del}\end{align*}
up to negligible contributions from $u^{(0)}$. Using again \eqref{eq:id3}, we obtain terms that are similar to \eqref{eq:udiff}. Following the subsequent arguments, one can show that the difference in the previous display is of order $|y^{(\prime)}|(\ov\la_n\Del)^{-\al/4}$, which, together with \eqref{eq:y2}, means that the second $\Pi^n$-integral is of order $(\ov\la_n\Del)^{-\al/2}\Del$. Altogether, $\ov w^{\prime n,i}_{\nu_1,\nu_2}-\ov w^{\prime\prime n,i}_{\nu_1,\nu_2}$ is of magnitude $(\ov\la_n\Del)^{-1/2-\al/4}\Del$.

A similar reasoning shows that $\ov c^{\prime n,i}_{\nu_1,\nu_2}-\ov c^{\prime\prime n,i}_{\nu_1,\nu_2}$ is of the same size. Moreover, 
by assumption H1, we know from \cite[(D.45) and (D.46)]{Chong19a} that $\ov\mu_f$ is twice differentiable in its arguments with all derivatives up to order two having at most polynomial growth. In particular, the mean value theorem and the identities in \eqref{eq:thetamu} imply that $\theta^{n,i}_4-\theta^{n,i}_5$ is of size $\Delh+(\ov\la_n\Del)^{-1/2-\al/4}\Del$. Thus, recalling that $\theta^{n,i}_4-\theta^{n,i}_5$ has a vanishing $\calf^n_{i-\la_n}$-conditional expectation, a martingale size estimate finally shows $\E[(B^{n,4}_2-B^{n,5}_2)^\ast_T]\lec \la_n^{1/2}\Delh+\la_n^{1/2}(\ov\la_n\Del)^{-1/2-\al/4}\Del = (\la_n\Del)^{1/2} +\la_n^{1/2}(\ov\la_n\Del)^{1/2-\al/4}\ov\la_n^{-1}\to0$ by comparison with \eqref{eq:help6}.
\epr

\bpr[Proof of Lemma~\ref{lem:easy}] Recall from Remark~\ref{rem:2} that we may assume $M=K=1$. In this case, $|\ov v^{\prime\prime n,i}_{\nu_1,\nu_2}-\ov v''_n(i\Del)_{\nu_1,\nu_2}|$, $|\ov w^{\prime\prime n,i}_{\nu_1,\nu_2}-\ov w''_n(i\Del)_{\nu_1,\nu_2}|$, and $|\ov c^{\prime\prime n,i}_{\nu_1,\nu_2}-\ov c''_n(i\Del)_{\nu_1,\nu_2}|$ are of magnitude $\Del^{\Theta  \ov a}$ by \eqref{eq:Theta} and \eqref{eq:PinGa}. 
Furthermore, we have, similarly to \eqref{eq:condexp0}, for all values of $s$,
\beq\label{eq:condexp0-2} \E[\ov\mu_f(\ov v''_n(s),\ov w''_n(s), \ov c''_n(s)) \mid \calf_{s-\la_n\Del}] =\E[\ov\mu_f(\ov v_n(s),\ov w_n(s), \ov c_n(s)) \mid \calf_{s-\la_n\Del}] = 0.\eeq
In particular, $\E[(\theta^{n,i}_5)_\nu-(\theta^{n,i}_6)_\nu \mid \calf^n_{i-\la_n}] = 0$ for all values of $i$.
Thus, a martingale size estimate  yields $\E[(B^{n,5}_2-B^{n,6}_2)^\ast_T]\lec\la_n^{1/2}\Del^{\Theta \ov a} \leq \Del^{-a/2+\Theta \ov a} \to 0$ because $\ov a>\frac a \Theta$.

Concerning the second statement of the lemma, observe that $|\ov v''_n(i\Del)_{\nu_1,\nu_2}-\ov v_n(i\Del)_{\nu_1,\nu_2}|$, $|\ov w''_n(i\Del)_{\nu_1,\nu_2}-\ov w_n(i\Del)_{\nu_1,\nu_2}|$, and $|\ov c''_n(i\Del)_{\nu_1,\nu_2}-\ov c_n(i\Del)_{\nu_1,\nu_2}|$ are all of size at most $(\ov\la_n\Del)^{1/2-\al/4}$. Thus, combined with a martingale size estimate, which is possible thanks to \eqref{eq:condexp0-2}, we deduce $\E[(B^{n,6}_2-B^{n,7}_2)^\ast_T]\lec \la_n^{1/2}(\ov\la_n\Del)^{1/2-\al/4}$, which vanishes in the limit as $n\to\infty$ by \eqref{eq:help6}.
\epr

\bpr[Proof of Lemma~\ref{lem:B56}] Recalling \eqref{eq:condexp0-2}, we have $B^{n,7}_2(t)-B^{n,8}_2(t)=I^n_1(t)+I^n_2(t)$, where
\begin{align*} I^n_1(t)&:=\frac1{\sqrt{\Del}} \sumtla \int_{(i-1)\Del}^{i\Del} \sum_{\nu\in\cali}\bigg[ \ov\mu_f(\ov v_n(i\Del),\ov w_n(i\Del), \ov c_n(i\Del))_\nu-\ov\mu_f(\ov v_n(s),\ov w_n(s), \ov c_n(s))_\nu\bigg]\,\dd s,\\
I^n_2(t) &:=
-\frac1{\sqrt{\Del}}\int_{[t/\Del]\Del}^t \sum_{\nu\in\cali} \ov\mu_f(\ov v_n(s),\ov w_n(s), \ov c_n(s))_\nu\,\dd s. \end{align*}
The boundary term $I^n_2$ is negligible (cf.\ the proof of Lemma~\ref{lem:iDels}). For $I^n_1$, we have from \eqref{eq:siest}, \eqref{eq:Gest}, and the differentiability of $\ov \mu_f$ that the difference in brackets is of size $(i\Del-s)^{1/2-\al/4}\leq \Del^{1/2-\al/4}$. Moreover, by \eqref{eq:condexp0-2}, the $i$th term in $I^n_1(t)$ is $\calf^n_i$-measurable with vanishing $\calf^n_{i-1-\la_n}$-conditional expectation. Therefore, a martingale size estimate yields $\E[(I^n_1)^\ast_T]\lec \la_n^{1/2}\Del^{1/2-\al/4}\leq \Del^{-a/2+1/2- \al/ 4}\leq \Del^{1/4-a/2}\to0$, because the exponent is close to $\frac14-\frac1{4\Theta}>0$. 
\epr

\bpr[Proof of Lemma~\ref{lem:Dn21}] 
Notice that $D^{n,1}_2$ is the second-order term that arises from a Taylor expansion of the $\si$-difference in
\begin{align*} &\frac{(\Delta^n_i u_{(i-\la_n)\Del}^{(i+L-1)\Del}-\Delta^n_i  u_{(i-\la_n)\Del,(i-\la_n)\Del}^{(i+L-1)\Del})_\nu}{\tau_n} \\
&\quad=\iint_{(i-\la_n)\Del}^{(i+l-1)\Del} \frac{\Delta^n_i G_y(s)_\nu}{\tau_n}\Bigl(\si(u(s,y))-\si(u(s,y)_0^{(i-\la_n)\Del})\Bigr) \,W(\dd s,\dd y),
\end{align*}
which is part of $\delta^{n,i}_2$.
The third-order term involves $(u(s,y)-u(s,y)_0^{(i-\la_n)\Del})^3$, which is of size $((\la_n\Del)^{1/2-\al/4})^3= o(\Delh)$ by \eqref{eq:help13} if $a$ is sufficiently close to $\frac1{2\Theta}$.

So it remains to prove that the first-order term
\allowdisplaybreaks[4]
\begin{align*}
K^n_1(t)&:=\Delh\sumtla\sum_{\nu\in\cali} \E\Biggl[ \partial_\nu f\Biggl( \frac{\Delta^n_i  u_{(i-\la_n)\Del,(i-\la_n)\Del}^{(i+L-1)\Del}}{\tau_n}\Biggr)\iint_{(i-\la_n)\Del}^{(i+l-1)\Del} \frac{\Delta^n_i G_y(s)_\nu}{\tau_n} \\
&\qtimes\si'(u(s,y)_0^{(i-\la_n)\Del})\Bigl(u(s,y)-u(s,y)_0^{(i-\la_n)\Del}\Bigr) \,W(\dd s,\dd y)\mathrel{\bigg|}\calf^n_{i-\la_n}\Biggr]\\
&\hphantom{:}=\Delh\sumtla\sum_{\nu\in\cali} \E\Biggl[ \partial_\nu f\Biggl( \frac{\Delta^n_i  u_{(i-\la_n)\Del,(i-\la_n)\Del}^{(i+L-1)\Del}}{\tau_n}\Biggr)\iint_{(i-\la_n)\Del}^{(i+l-1)\Del} \frac{\Delta^n_i G_y(s)_\nu}{\tau_n} \\
&\qtimes\si'(u(s,y)_0^{(i-\la_n)\Del})\iint_{(i-\la_n)\Del}^s G_{y-z}(s-r)\si(u(r,z))\,W(\dd r,\dd z) \,W(\dd s,\dd y)\mathrel{\bigg|}\calf^n_{i-\la_n}\Biggr]
\end{align*}
is also negligible. If we had $\si(u(r,z)_0^{(i-\la_n)\Del})$ instead of $\si(u(r,z))$ in the last line, the conditional expectation would be identically zero because the double integral belongs to the second Wiener chaos (conditionally on $\calf^n_{i-\la_n}$), while the $\partial_\nu f$-term only has components from the odd-order Wiener chaoses. Thus, it suffices to show that
\beq\label{eq:Kn2}\begin{split}
	K^n_2(t)&:=\Delh\sumtla\sum_{\nu\in\cali} \E\Biggl[ \partial_\nu f\Biggl( \frac{\Delta^n_i  u_{(i-\la_n)\Del,(i-\la_n)\Del}^{(i+L-1)\Del}}{\tau_n}\Biggr)\iint_{(i-\la_n)\Del}^{(i+l-1)\Del} \frac{\Delta^n_i G_y(s)_\nu}{\tau_n} \\
	&\qtimes\si'(u(s,y)_0^{(i-\la_n)\Del})\iint_{(i-\la_n)\Del}^s G_{y-z}(s-r)\\
	&\qtimes\Bigl(\si(u(r,z))-\si(u(r,z)_0^{(i-\la_n)\Del})\Bigr)\,W(\dd r,\dd z) \,W(\dd s,\dd y)\mathrel{\bigg|}\calf^n_{i-\la_n}\Biggr]
\end{split}\eeq
vanishes as $n$ gets large. Using again Taylor's formula on the $\si$-difference, we only need to keep the first-order term given by
\beq\label{eq:Kn3}\begin{split}
K^n_3(t)&:=\Delh\sumtla\sum_{\nu\in\cali} \E\Biggl[ \partial_\nu f\Biggl( \frac{\Delta^n_i  u_{(i-\la_n)\Del,(i-\la_n)\Del}^{(i+L-1)\Del}}{\tau_n}\Biggr)\iint_{(i-\la_n)\Del}^{(i+l-1)\Del} \frac{\Delta^n_i G_y(s)_\nu}{\tau_n} \\
&\qtimes\si'(u(s,y)_0^{(i-\la_n)\Del})\iint_{(i-\la_n)\Del}^s G_{y-z}(s-r)\si'(u(r,z)_0^{(i-\la_n)\Del})\\
&\qtimes\iint_{(i-\la_n)\Del}^r G_{z-w}(r-v)\si(u(v,w))\,W(\dd v,\dd w)\,W(\dd r,\dd z) \,W(\dd s,\dd y)\mathrel{\bigg|}\calf^n_{i-\la_n}\Biggr],
\end{split}\eeq
because the second-order term contains two $(\la_n\Del)^{1/2-\al/4}$-factors, in addition to a third one from  the $W(\dd r,\dd z)$-integral of $G_{y-z}(s-r)$ in \eqref{eq:Kn2}, which makes it $o(\Delh)$ as seen above. 

Observe that $K^n_3$ already contains two terms of size $(\la_n\Del)^{1/2-\al/4}$: the  $W(\dd r,\dd z)$-integral of $G_{y-z}(s-r)$ and the $W(\dd v,\dd w)$-integral of $G_{z-w}(r-v)$. Thus, we may further replace $\si(u(v,w))$ by $\si(u(v,w)_0^{(i-\la_n)\Del})$, as this approximation would give a third $(\la_n\Del)^{1/2-\al/4}$-factor. Let us assume for the moment that we can, at the same time, replace the argument of the $\partial_\nu f$-term by $\Delta^n_i  u_{(i-\la^{\prime(q)}_n)\Del,(i-\la_n)\Del}^{(i+L-1)\Del}/\tau_n$, where $\la^{\prime (q)}_n:=[\Del^{-a^{\prime (q)}}]$ and $a^{\prime (q)}\leq a$ will be specified later. 
Then, instead of $K^n_3$, we only have to consider
\begin{align*}
K^n_{4,q}(t)&:=\Delh\sumtla\sum_{\nu\in\cali} \E\Biggl[ \partial_\nu f\Biggl( \frac{\Delta^n_i  u_{(i-\la^{\prime(q)}_n)\Del,(i-\la_n)\Del}^{(i+L-1)\Del}}{\tau_n}\Biggr)\iint_{(i-\la_n)\Del}^{(i+l-1)\Del} \frac{\Delta^n_i G_y(s)_\nu}{\tau_n} \\
&\qtimes\si'(u(s,y)_0^{(i-\la_n)\Del})\iint_{(i-\la_n)\Del}^s G_{y-z}(s-r)\si'(u(r,z)_0^{(i-\la_n)\Del})\\
&\qtimes\iint_{(i-\la_n)\Del}^r G_{z-w}(r-v)\si(u(v,w)_0^{(i-\la_n)\Del})\,W(\dd v,\dd w)\,W(\dd r,\dd z) \,W(\dd s,\dd y)\mathrel{\bigg|}\calf^n_{i-\la_n}\Biggr].
\end{align*}

Now the crucial observation is this: the $W(\dd r,\dd z)$-integral can be split into an integral on $((i-\la_n)\Del,(i-\la^{\prime(q)}_n)\Del]\times\R^d$ and a second one on $((i-\la^{\prime(q)}_n)\Del,s]\times\R^d$. The contribution of the former one is exactly zero by Lemma~\ref{lem:chaos}. Indeed, conditionally on $\calf^n_{i-\la_n}$, the argument of the $\partial_\nu f$-term is a Wiener integral starting at $(i-\la^{\prime(q)}_n)\Del$, while the triple $W$-integral has one part, namely the $W(\dd r,\dd z)$-integral, that ends at $(i-\la^{\prime(q)}_n)\Del$. So only the contribution of the second integral on $((i-\la^{\prime(q)}_n)\Del,s]\times\R^d$ remains. But for this part, the same argument shows that we can also reduce the $W(\dd v,\dd w)$-integral to $((i-\la^{\prime(q)}_n),r]$. Once this is done, we can further approximate the $\partial_\nu f$-argument by $\Delta^n_i  u_{(i-\la^{\prime(q+1)}_n)\Del,(i-\la_n)\Del}^{(i+L-1)\Del}/\tau_n$ (i.e., approximate $K^n_{4,q}$ by $K^n_{4,q+1}$) as long as $\Del^{-1/2}((\la^{\prime(q)}_n\Del)^{1/2-\al/4})^2\Del^{\Theta a^{\prime (q+1)}/2}\to0$, that is, as long as
\[ -\frac12 + (1-a^{\prime (q)})\biggl(1-\frac\al2\biggr)+\frac\Theta2 a^{\prime (q+1)}>0. \]

Starting with $\la_n^{\prime(1)}=\la_n$ [i.e., $a^{\prime (1)}:= a$] and iterating the argument above $Q'$ times, we end up with $K^n_{4,Q'}$, which, after shrinking the $W(\dd r,\dd z)$- and the $W(\dd v,\dd w)$-integral as described above, is of size
\[ \E[(K^n_{4,Q'})^\ast_T]\lec \Del^{-\frac12}((\la_n^{\prime(Q')}\Del)^{\frac12-\frac\al4})^2. \]
This goes to zero if $-\frac12+(1- a^{\prime (Q')})(1-\frac\al2)>0$, and because $\al<1$, this happens if $a^{\prime (Q')}$ is sufficiently small. As a consequence, similarly to the proof of Lemma~\ref{lem:remove-s}, we are left to show that the sequence defined via
\[ -\frac12 + (1-b'_n)\biggl(1-\frac\al2\biggr)+\frac\Theta2 b'_{n+1}=0,\qquad b'_1=a, \]
gets arbitrarily close to or smaller than $0$ as $n\to\infty$. To this end, note that $b'_{n+1}=h(b'_n)$ with $h(x) := \frac{2}{\Theta}(\frac12-(1-x)(1-\frac\al2))$.  Recalling that $\Theta = 1+\frac\al2$, it is easy to verify the function $h$ is a line that is parallel and below the identity for $\al=\frac23$, so in this case, $b'_{n+1}-b'_n = h(b'_n)-b'_n<0$ for all $n\in\N$. If $\al\in(0,1)\setminus\{\frac23\}$, then $h$ has positive slope, satisfies $h(0)<0$, and intersects the identity at $x=\frac{2(\al-1)}{3\al-2}$, which is either strictly negative (if $\al>\frac23$) or larger than $\frac{1}{2\Theta}$ (if $\al<\frac23$). Thus, if $a$ is close enough to $\frac{1}{2\Theta}$ (so that it is smaller than the point of intersection), we have $b'_{n+1}-b_n = h(b'_n)-b'_n<0$ for all $n\in\N$ as well. We conclude that $b'_n$ converges to a negative limit if $\al>\frac23$ and to $-\infty$ if $\al\leq \frac23$, which implies what we wanted to show.
\epr

\bpr[Proof of Lemma~\ref{lem:Dn22}] The $u$-difference in \eqref{eq:wtrhophi-2} equals $\iint_{(i-\la_n)\Del}^s G_{y-z}(s-r)\si(u(r,z))\,W(\dd r,\dd z)$. The reader can easily verify that $\si(u(r,z))$ can be changed to $\si(u(r,z)_0^{(i-\la_n)\Del})$ because this will incur an error of size $\Del^{-1/2}((\la_n\Del)^{1/2-\al/4})^3$, which is negligible by \eqref{eq:help13}. With this modification, if we expand the square using the integration by parts formula, then $\wt\rho^{n,i}_2$ corresponds exactly to the quadratic variation part. Therefore, we are left to show that 
\begin{align*}
&\Delh\sumtla\sum_{\nu\in\cali} 
\E\Biggl[ \partial_\nu f\Biggl( \frac{\Delta^n_i  u_{(i-\la_n)\Del,(i-\la_n)\Del}^{(i+L-1)\Del}}{\tau_n}\Biggr)\iint_{(i-\la_n)\Del}^{(i+l-1)\Del} \frac{\Delta^n_i G_y(s)_\nu}{\tau_n} \si''(u(s,y)_0^{(i-\la_n)\Del})\\
&\qtimes\iint_{(i-\la_n)\Del}^s G_{y-z}(s-r)\si(u(r,z)_0^{(i-\la_n)\Del})\\
&\qtimes\iint_{(i-\la_n)\Del}^r G_{y-w}(s-v)\si(u(v,w)_0^{(i-\la_n)\Del})\,W(\dd v,\dd w)\,W(\dd r,\dd z) \,W(\dd s,\dd y)\mathrel{\bigg|} \calf^n_{i-\la_n}\Biggr]\limL0.
\end{align*}
The left-hand side has exactly the same structure as the term $K^n_3$ in \eqref{eq:Kn3} [when $\si(u(v,w))$ in \eqref{eq:Kn3} is replaced by $\si(u(v,w)_0^{(i-\la_n)\Del})$]. So identical arguments can be employed to show the statement of the previous display.
\epr

\bpr[Proof of Lemma~\ref{lem:Dn234}]
For the first statement, we claim that changing the domain of integration of the $\La(\dd z,\dd z')\,\dd r$-integrals in $\wt w^{n,i}$ and $\wt c^{n,i}$ to $[0,\la_n\Del]\times\R^d\times\R^d$ only leads to an asymptotically negligible error. Indeed, the difference resulting from modifying $\wt c^{n,i}_{\nu_1,\nu_2}$ in this way is, by a size estimate, of order
\beq\label{eq:intG}\begin{split} &\iiint_{(\la_n+l_{12}-1)\Del-s}^{\la_n\Del} G_z(r)G_{z'}(r)\,\La(\dd z,\dd z')\,\dd r \leq  \iiint_{\la_n\Del-s}^{\la_n\Del} G_z(r)G_{z'}(r)\,\La(\dd z,\dd z')\,\dd r\\
&\quad\lec (\la_n\Del)^{1-\frac\al2}-(\la_n\Del-s)^{1-\frac\al2} \leq (\la_n\Del)^{1-\frac\al2}\biggl(1-\biggl(1-\frac{s}{\la_n\Del}\biggr)^{1-\frac\al2}\biggr)\lec (\la_n\Del)^{1-\frac\al2}\frac{s}{\la_n\Del}, \end{split}\raisetag{3\baselineskip}\eeq
where the last step is true by the mean value theorem if $\frac{s}{\la_n\Del}$ is bounded away from $1$. Let $\la'_n:=[\Del^{-a'}]$, where $a'$ is chosen within $(\frac{a}{2\Theta},\frac a2)$. Then the last $\Pi^n$-integral in \eqref{eq:wc}, restricted to $[0,\la'_n\Del]\times\R^d\times\R^d$ and with the $\La(\dd z,\dd z')\,\dd r$-integral taken on $((\la_n+l_{12}-1)\Del-s,\la_n\Del]\times\R^d\times\R^d$, is of size $(\la_n\Del)^{1-\frac\al2}\frac{\la'_n}{\la_n} \leq \Del^{(1-a)/2+a-a'} = \Del^{1/2+a/2-a'}=o(\Delh)$ because $a'<\frac a2$. For the  same integral on $[\la'_n\Del,(\la_n+l_{12}-1)\Del]\times\R^d\times\R^d$, we simply estimate the $\La(\dd z,\dd z')\,\dd r$-integral by $(\la_n\Del)^{1-\frac\al2}$ and use \eqref{eq:Theta} to obtain in total the bound $(\la_n\Del)^{1-\frac\al2}\Del^{\Theta a'}\leq \Del^{(1-a)/2+\Theta a'} = \Del^{1/2-a/2+\Theta a'}=o(\Delh)$ because $a'>\frac{a}{2\Theta}$. An analogous argument shows that the modifications described at the beginning of the proof only change $\wt w^{n,i}$ by an $o(\Delh)$-term. As $\ov\mu_f$ is differentiable (see the proof of Lemma~\ref{lem:Bn345}), our claim above follows from the mean value theorem. With the same reasoning as for the term $E^{n,i}_3$ in the proof of Lemma~\ref{lem:remove-s}, we may further restrict the domain of integration of the $\Pi^n$-integrals to $[0,\la_n\Del]\times\R^d\times\R^d$.

Let us now consider the remaining changes to go from $(\wt v^{n,i},\wt w^{n,i},\wt c^{n,i})$ to $(\wt v^{\prime n,i},\wt w^{\prime n,i},\wt c^{\prime n,i})$. To this end,
note that $\rho^{n,i}_3=\E[\wt \rho^{n,i}_3\mid \calf^n_{i-\la_n}]$, where $(\wt \rho^{n,i}_3)_\nu$ arises from replacing $u(s|r,y|z^{(\prime)})_0^{(i-\la_n)\Del}$ by $u(i\Del,y|z^{(\prime)})_0^{(i-\la_n)\Del}$ in
\beq\label{eq:rhoni3}\begin{split}
	  &\frac12\partial_\nu f\Biggl( \iint_{(i-\la_n)\Del}^{(i+L-1)\Del} \frac{\Delta^n_i G_y(s)}{\tau_n} \si(u(s,y)_0^{(i-\la_n)\Del})\,W(\dd s,\dd y)\Biggr)\\
	&\qtimes\iint_{(i+l-1-\la_n)\Del}^{(i+l-1)\Del} \frac{\Delta^n_i G_y(s)_\nu}{\tau_n} \si''(u(s,y)_0^{(i-\la_n)\Del})\iiint_{s-\la_n\Del}^s G_{y-z}(s-r)G_{y-z'}(s-r)\\
	&\cquad{\mathrel{\times}\si(u(r,z)_0^{(i-\la_n)\Del})\si(u(r,z')_0^{(i-\la_n)\Del})\,\La(\dd z,\dd z')\,\dd r  \,W(\dd s,\dd y).}
\end{split}\eeq
The $\La(\dd z,\dd z')\,\dd r$-integral is already of size $(\la_n\Del)^{1-\al/2}$. By \eqref{eq:siest} and \eqref{eq:id2}, fixing the time variable of $u$ at $i\Del$ yields another factor of $(\la_n\Del)^{1/2-\al/4}$, which together is $o(\Delh)$ by \eqref{eq:help13} and completes the proof of $D^{n,2}_2-D^{n,3}_2\limL0$.

As in the proof of Lemma~\ref{lem:Bn345}, $D^{n,3}_2-D^{n,4}_2\limL0$ follows if we can show that $\wt v^{\prime n,i}-\wt v^{\prime\prime n,i}$, $\wt w^{\prime n,i}-\wt w^{\prime\prime n,i}$, and $\wt c^{\prime n,i}-\wt c^{\prime\prime n,i}$ are of size $o(\Delh)$, uniformly in $i$. For $\wt v^{\prime n,i}-\wt v^{\prime\prime n,i}$, this has been done in the proof of Lemma~\ref{lem:remove-y} [because $\wt v^{\prime(\prime) n,i}=v^{\prime(\prime) n,i}$ by definition]. For $\wt w^{\prime n,i}-\wt w^{\prime\prime n,i}$, simply notice that for every $p>0$,
\beq\label{eq:sizew} \E[|\wt w^{\prime n,i}|^p]^{\frac1p}+\E[|\wt w^{\prime\prime n,i}|^p]^{\frac1p}\lec((\la_n\Del)^{1-\frac\al2})^2\eeq
 by a standard size estimate. Finally, notice that the $\La(\dd z,\dd z')\,\dd r$-integral in $\wt c^{\prime n,i}$ is already of size $(\la_n\Del)^{1-\al/2}$. By \eqref{eq:id2} and Lemma~\ref{lem:momincr}, taking the difference $\wt c^{\prime n,i}-\wt c^{\prime\prime n,i}$ gives us a factor of $(\lvert y\rvert^{1-\al/2-\eps}+\lvert y'\rvert^{1-\al/2-\eps})$ (for any $\eps>0$) within the $\Pi^n$-integral. By \eqref{eq:Pin1}, the total mass of all the measures $|\Pi^n_{r,h}|$ is uniformly bounded in $n$, $r$, and $h$, so Jensen's inequality and \eqref{eq:y2} show that
\beq\label{eq:yalpha} \iiint_0^{\la_n\Del} (\lvert y \rvert^{1-\frac\al2-\eps}+\lvert y'\rvert^{1-\frac\al2-\eps})\,|\Pi^n_{|l_1-l_2|,x_{k_1}-x_{k_2}}|(\dd s,\dd y,\dd y') \lec \Del^{\frac12-\frac\al4-\frac\eps2}. \eeq
Altogether, if $\eps>0$ is small enough, then $\E[|\wt c^{\prime n,i}-\wt c^{\prime\prime n,i}|]\lec (\la_n\Del)^{1-\al/2}\Del^{1/2-\al/4-\eps/2} = o(\Delh)$ by comparison with \eqref{eq:help13}.
\epr

\bpr[Proof of Lemma~\ref{lem:Dn2456}]
The proof of $D^{n,4}_2-D^{n,5}_2\limL0$ is completely analogous to that of Lemma~\ref{lem:remove-Pi}. For $D^{n,5}_2-D^{n,6}_2\limL0$, observe that this difference is given by
\begin{align*}
\frac{1}{\Delh} \sumtla\sum_{\nu\in\cali}\int_{(i-1)\Del}^{i\Del} \Bigl(\ov\mu_f(\wt v(i\Del),\wt w(i\Del), \wt c(i\Del))_\nu-\ov\mu_f(\wt v(s),\wt w(s), \wt c(s))_\nu\Bigr)\,\dd s
\end{align*}
plus some boundary terms that are $o(\Delh)$. Now develop the $\ov\mu_f$-difference up to second order using Taylor's theorem. Because $\wt w(i\Del)-\wt w(s)$ is of size $o(\Delh)$ by \eqref{eq:sizew} and the same holds true for $\wt c(i\Del)-\wt c(s)$ [the $\La(\dd z,\dd z')\,\dd r$-integral gives $(\la_n\Del)^{1-\al/2}$ and approximating $u(i\Del,x_k)_0^{(i-\la_n)\Del}$ by $u(s,x_k)_0^{(i-\la_n)\Del}$ gives another $\Del^{1/2-\al/4}$], only the terms corresponding to taking one derivative of $\ov\mu_f(v,w,c)$ in the $v$-variables remain [the reader may easily confirm that second-order terms are $o(\Delh)$ as well]. And for those, one can proceed in the same way as for the terms $G^n_1$ and $G^n_2$ in the proof of Lemma~\ref{lem:iDels}. We omit the details at this point.
\epr

\bpr[Proof of Lemma~\ref{lem:Dn3}]
The second and third fraction in $\delta^{n,i}_3$ are both of size $(\la_n\Del)^{1/2-\al/4}$. As in the proofs concerning $D^n_2$, we are allowed to make any modification that gives rise to an additional error term of size $(\la_n\Del)^{1/2-\al/4}$. Thus, if we write out the last two fractions in $\delta^{n,i}_3$ and Taylor expand $\si(u(s,y))-\si(u(s,y)_0^{(i-\la_n)\Del})$, we only need to keep the first-order term 
\begin{align*}
&\si'(u(s,y)_0^{(i-\la_n)\Del})(u(s,y)-u(s,y)_0^{(i-\la_n)\Del}) \\
&\quad= \si'(u(s,y)_0^{(i-\la_n)\Del})\iint_{(i-\la_n)\Del}^s G_{y-z}(s-r)\si(u(r,z))\,W(\dd r,\dd z).\end{align*} Clearly,  $\si(u(r,z))$ here can be replaced by $\si(u(r,z)_0^{(i-\la_n)\Del})$. Furthermore, in $u(s,y)_0^{(i-\la_n)\Del}$, we can fix the time coordinate at $i\Del$. As a result, we have $D^n_3-D^{n,1}_3\limL0$, where $D^{n,1}_3(t):=\frac{\Delh}{2}\sumtla\sum_{\nu_1,\nu_2\in\cali} (\phi^{n,i}_{1})_{\nu_1,\nu_2}$, $(\phi^{n,i}_{1})_{\nu_1,\nu_2}:=\E[(\wt \phi^{n,i}_{1})_{\nu_1,\nu_2}\mid\calf^n_{i-\la_n}]$, and
\begin{align*}
(\wt \phi^{n,i}_{1})_{\nu_1,\nu_2}&:=\partial^2_{\nu_1\nu_2} f(\chi^n_i) \iint_{(i-\la_n)\Del}^{(i+l_1-1)\Del} \frac{\Delta^n_i G_y(s)_{\nu_1}}{\tau_n}\si'(u(i\Del,y)_0^{(i-\la_n)\Del})\\
&\cquad{\mathrel{\times} \iint_{(i-\la_n)\Del}^s G_{y-z}(s-r)\si(u(r,z)_0^{(i-\la_n)\Del})\,W(\dd r,\dd z)\,W(\dd s,\dd y)}\\
&\qtimes\iint_{(i-\la_n)\Del}^{(i+l_2-1)\Del} \frac{\Delta^n_i G_y(s)_{\nu_2}}{\tau_n}\si'(u(i\Del,y)_0^{(i-\la_n)\Del})\\
&\cquad{\mathrel{\times} \iint_{(i-\la_n)\Del}^s G_{y-z}(s-r)\si(u(r,z)_0^{(i-\la_n)\Del})\,W(\dd r,\dd z)\,W(\dd s,\dd y)}
\end{align*}
and 
\beq\label{eq:chi} \chi^n_i:=\iint_{(i-\la_n)\Del}^{(i+L-1)\Del} \frac{\Delta^n_i G_y(s)}{\tau_n} \si(u(i\Del,y)_0^{(i-\la_n)\Del})\,W(\dd s,\dd y).\eeq

Next, we use integration by parts to compute the product of the two $W(\dd s,\dd y)$-integrals in $(\wt\phi^{n,i}_1)_{\nu_1,\nu_2}$, leading to $D^{n,1}_3=D^{n,1,1}_3+D^{n,1,2}_3$, where $D^{n,1,j}_3(t):=\frac{\Delh}{2}\sumtla\sum_{\nu_1,\nu_2\in\cali} (\phi^{n,i}_{1,j})_{\nu_1,\nu_2}$ for $j=1,2$, $(\phi^{n,i}_{1,1})_{\nu_1,\nu_2}:=\E[(\wt \phi^{n,i}_{1,1})_{\nu_1,\nu_2}+(\wt \phi^{n,i}_{1,1})_{\nu_2,\nu_1}\mid\calf^n_{i-\la_n}]$, $(\phi^{n,i}_{1,2})_{\nu_1,\nu_2}:=\E[(\wt \phi^{n,i}_{1,2})_{\nu_1,\nu_2}\mid\calf^n_{i-\la_n}]$,
\begin{align*}
(\wt \phi^{n,i}_{1,1})_{\nu_1,\nu_2}&:=\partial^2_{\nu_1\nu_2} f(\chi^n_i) \iint_{(i-\la_n)\Del}^{(i+l_1-1)\Del} \frac{\Delta^n_i G_y(s)_{\nu_1}}{\tau_n}\si'(u(i\Del,y)_0^{(i-\la_n)\Del})\\
&\qtimes\iint_{(i-\la_n)\Del}^s G_{y-z}(s-r)\si(u(r,z)_0^{(i-\la_n)\Del})\,W(\dd r,\dd z)\\
&\qtimes \iint_{(i-\la_n)\Del}^{s} \frac{\Delta^n_i G_z(r)_{\nu_2}}{\tau_n}\si'(u(i\Del,z)_0^{(i-\la_n)\Del})\\
&\cquad{\mathrel{\times}\iint_{(i-\la_n)\Del}^r G_{z-w}(r-v)\si(u(v,w)_0^{(i-\la_n)\Del})\,W(\dd v,\dd w)\,W(\dd r,\dd z)\,W(\dd s,\dd y)},
\end{align*}
and
\begin{align*}
(\wt \phi^{n,i}_{1,2})_{\nu_1,\nu_2}&:=\partial^2_{\nu_1\nu_2} f(\chi^n_i)\iiint_{(i-\la_n)\Del}^{(i+l_{12}-1)\Del} \frac{\Delta^n_i G_y(s)_{\nu_1}\Delta^n_i G_{y'}(s)_{\nu_2}}{\tau^2_n} \si'(u(i\Del,y)_0^{(i-\la_n)\Del})\\
&\qtimes\si'(u(i\Del,y')_0^{(i-\la_n)\Del})\iint_{(i-\la_n)\Del}^s G_{y-z}(s-r)\si(u(r,z)_0^{(i-\la_n)\Del})\,W(\dd r,\dd z)\\
&\qtimes  \iint_{(i-\la_n)\Del}^s G_{y'-z}(s-r)\si(u(r,z)_0^{(i-\la_n)\Del})\,W(\dd r,\dd z)\,\La(\dd y,\dd y')\,\dd s.
\end{align*}

Inside the $W(\dd s,\dd y)$-integral of $ (\wt \phi^{n,i}_{1,1})_{\nu_1,\nu_2}$, the product of the $W(\dd r,\dd z)$-integrals can again be evaluated using integration by parts, resulting in three terms: a triple $W$-integral [not counting the $W(\dd s,\dd y)$-integral!], a $W(\dd r,\dd z)$-integral of a product of two $W(\dd v,\dd w)$-integrals, and a $\La(\dd z,\dd z')\,\dd r$-integral of a single $W(\dd v,\dd w)$-integral. For the second term, integrate by parts once more to get two triple $W$-integrals and a $W(\dd r,\dd z)$-integral of a $\La(\dd w,\dd w')\,\dd v$-integral. This very last contribution is equal to $D^{n,2,1}_3(t):=\frac{\Delh}{2}\sumtla\sum_{\nu_1,\nu_2\in\cali} (\phi^{n,i}_{2,1})_{\nu_1,\nu_2}$, where  $(\phi^{n,i}_{2,1})_{\nu_1,\nu_2}:=\E[(\wt \phi^{n,i}_{2,1})_{\nu_1,\nu_2}+(\wt \phi^{n,i}_{2,1})_{\nu_2,\nu_1}\mid\calf^n_{i-\la_n}]$ and
\beq\label{eq:phi21}\begin{split}
	(\wt \phi^{n,i}_{2,1})_{\nu_1,\nu_2}&:=\partial^2_{\nu_1\nu_2} f(\chi^n_i)\iint_{(i-\la_n)\Del}^{(i+l_1-1)\Del} \frac{\Delta^n_i G_y(s)_{\nu_1}}{\tau_n}\si'(u(i\Del,y)_0^{(i-\la_n)\Del})\\
	&\qtimes \iint_{(i-\la_n)\Del}^{s} \frac{\Delta^n_i G_z(r)_{\nu_2}}{\tau_n}\si'(u(i\Del,z)_0^{(i-\la_n)\Del})\\ 
	&\cquad{\mathrel{\times}\iiint_{(i-\la_n)\Del}^r G_{y-w}(s-v)G_{z-w'}(r-v)\si(u(v,w)_0^{(i-\la_n)\Del})}\\
	&\cqquad{\mathrel{\times}	 	\si(u(v,w')_0^{(i-\la_n)\Del})\,\La(\dd w,\dd w')\,\dd v\,W(\dd r,\dd z)\,W(\dd s,\dd y).}
\end{split}\eeq

The crucial point is now the following: with the same trick as in the proof of Lemma~\ref{lem:Dn21} after \eqref{eq:Kn3}, all the other terms can be shown to be asymptotically negligible, that is, $D^{n,1,1}_3-D^{n,2,1}_3\limL0$ as $n\to\infty$. The trick does not apply to $D^{n,2,1}_3$ itself as it only allows us to reduce the domain of $W$-integrals, say, from $\iint_{(i-\la_n)\Del}^{\cdot}$ to $\iint_{(i-\la_n'')\Del}^{\cdot}$, where $\la_n''=[\Del^{-a''}]$ for some $0<a''<a$. Since $D^{n,2,1}_3$ is the only term where all $W$-integrals are taken on $\frac{\Delta^n_i G}{\tau_n}$-terms (and none of them is taken on a $G$-term), this procedure does not  decrease the asymptotic size of $(\wt \phi^{n,i}_{2,1})_{\nu_1,\nu_2}$. Nevertheless, it will be useful later and we denote this modified version of $(\wt \phi^{n,i}_{2,1})_{\nu_1,\nu_2}$ by $(\wt \phi^{n,i}_{3,1})_{\nu_1,\nu_2}$.

Similarly, if we multiply the $W(\dd r,\dd z)$-integrals in $(\wt \phi^{n,i}_{1,2})_{\nu_1,\nu_2}$, we only have to keep the quadratic covariation part, that is,  $D^{n,1,2}_3-D^{n,2,2}_3\limL0$, where $D^{n,2,2}_3(t):=\frac{\Delh}{2}\sumtla\sum_{\nu_1,\nu_2\in\cali} (\phi^{n,i}_{2,2})_{\nu_1,\nu_2}$, $(\phi^{n,i}_{2,2})_{\nu_1,\nu_2}:=\E[(\wt \phi^{n,i}_{2,2})_{\nu_1,\nu_2}\mid\calf^n_{i-\la_n}]$, and
\begin{align*}
(\wt \phi^{n,i}_{2,2})_{\nu_1,\nu_2}&:=\partial^2_{\nu_1\nu_2} f(\chi^n_i)\iiint_{(i-\la_n)\Del}^{(i+l_{12}-1)\Del} \frac{\Delta^n_i G_y(s)_{\nu_1}\Delta^n_i G_{y'}(s)_{\nu_2}}{\tau^2_n} \si'(u(i\Del,y)_0^{(i-\la_n)\Del})\\
&\qtimes\si'(u(i\Del,y')_0^{(i-\la_n)\Del})\iiint_{(i-\la_n)\Del}^s G_{y-z}(s-r)G_{y'-z'}(s-r)\\
&\cquad{\mathrel{\times}\si(u(r,z)_0^{(i-\la_n)\Del})\si(u(r,z')_0^{(i-\la_n)\Del})\,\La(\dd z,\dd z')\,\dd r \,\La(\dd y,\dd y')\,\dd s.}
\end{align*}

Consider now the $\La(\dd w,\dd w')\,\dd v$-integral in \eqref{eq:phi21}, which appears identically in $(\wt \phi^{n,i}_{3,1})_{\nu_1,\nu_2}$ as well. By changing variables $r-v\mapsto v$ and $(y-w,y-w')\mapsto (w,w')$, it can be rewritten as 
\beq\label{eq:help14}\begin{split} &\iiint_0^{r-(i-\la_n)\Del} G_w(v+s-r)G_{w'+z-y}(v)\\
	&\quad\times\si(u(r-v,y-w)_0^{(i-\la_n)\Del})\si(u(r-v,y-w')_0^{(i-\la_n)\Del})\,\La(\dd w,\dd w')\,\dd v. \end{split}\eeq
We want to replace $r-v$ by $i\Del-v$ and $s-r$ by $0$. The first operation is fine by the reasons outlined at the beginning of the proof. For the second operation, as ${|r-(i-\la_n)\Del|}\lec \la_n\Del$ and ${|s-r|}\lec \la''_n\Del$ in $(\wt \phi^{n,i}_{3,1})_{\nu_1,\nu_2}$, we can use \eqref{eq:der-G} to bound the size of the resulting error term by 
\begin{align*}
&	\iiint_0^{\la_n\Del} |G_w(v+s-r)-G_w(v)|G_{w'+z-y}(v)\,\La(\dd w,\dd w')\dd v \\ 
&\quad\leq\int_0^{\la''_n\Del} \iiint_0^{\la_n\Del} |\partial_t G(v+p,w)|G_{w'+z-y}(v)\,\La(\dd w,\dd w')\dd v\,\dd p \\
&\quad\lec \int_0^{\la''_n\Del} \int_0^{\la_n\Del} \frac{1}{v+p} \int_0^\infty \ee^{-(2v+p) q^2}q^{\al-1}\,\dd q\,\dd v\,\dd p \lec \int_0^{\la''_n\Del} \int_0^{\la_n\Del} \frac{1}{(v+p)(2v+p)^{\frac\al2}}\,\dd v\,\dd p\\
&\quad\leq \int_0^{\la''_n\Del} \int_0^{\la_n\Del} \frac{1}{(v+p)^{1+\frac\al2}}\,\dd v\,\dd p\lec \int_0^{\la''_n\Del} (p^{-\frac\al2}-(\la_n\Del+p)^{-\frac\al2})\,\dd p \\
&\quad\lec (\la''_n\Del)^{1-\frac\al2}-(\la_n\Del+\la''_n\Del)^{1-\frac\al2}+(\la_n\Del)^{1-\frac\al2} \lec (\la_n\Del)^{1-\frac\al2}((\textstyle\frac{\la''_n}{\la_n})^{1-\frac\al2}-(1+\frac{\la''_n}{\la_n})^{1-\frac\al2}+1)\\
&\quad\lec (\la_n\Del)^{1-\frac\al2}((\textstyle\frac{\la''_n}{\la_n})^{1-\frac\al2}+\frac{\la''_n}{\la_n})\lec (\la_n\Del)^{1-\frac\al2}(\textstyle\frac{\la''_n}{\la_n})^{1-\frac\al2} = (\la''_n\Del)^{1-\frac\al2}.
\end{align*}
So if we choose $a''\in(0,1-\frac{1}{2-\al})$, which is possible because $\al\in(0,1)$, the modification above only leads to an $o(\Delh)$-error.

Next, we also want to replace $z-y$ in \eqref{eq:help14} by $0$, which by \eqref{eq:g-firstder} only incurs an error of size 
\begin{align*}
&\iiint_0^{\la_n\Del} G_w(v)|G_{w'+z-y}(v)-G_{w'}(v)|\,\La(\dd w,\dd w')\dd v\\
&\quad \leq \sum_{j=1}^d \int_0^1 \iiint_0^{\la_n\Del} G_w(v) |\partial_{x_j} G(v,w'+p(z-y))(z_j-y_j)|\,\La(\dd w,\dd w')\,\dd v\,\dd p\\
&\quad\lec |z-y| \int_0^{\la_n\Del} v^{-\frac12} \int_0^\infty \ee^{-2vq^2}q^{\al-1}\,\dd q\,\dd v\lec |z-y|\int_0^{\la_n\Del} v^{-\frac{1+\al}{2}}\,\dd v \lec |z-y|(\la_n\Del)^{\frac{1-\al}{2}}.
\end{align*}
Since ${|z-y|}={|(x_{k_{12}}-z)-(x_{k_{12}}-y)|}\leq{|x_{k_{12}}-z|}+{|x_{k_{12}}-y|}$, either the $W(\dd s,\dd y)$- or the $W(\dd r,\dd z)$-integral in \eqref{eq:phi21} will give us a factor of $\Del^{1/2}$ by \eqref{eq:y2}. Since $(\la_n\Del)^{\frac{1-\al}{2}}\to0$ [because $\al\in(0,1)$], we have shown that $z-y$ in \eqref{eq:help14} can be dropped in the limit $n\to\infty$. Finally, attaching back the integrals from $(i-\la_n)\Del$ to $(i-\la''_n)\Del$ (using the trick from Lemma~\ref{lem:Dn21} ``backwards''), we conclude that $D^{n,2,1}_3-D^{n,4,1}_3\limL0$, where $D^{n,4,1}_3(t):=\frac{\Delh}{2}\sumtla\sum_{\nu_1,\nu_2\in\cali} (\phi^{n,i}_{4,1})_{\nu_1,\nu_2}$ with  $(\phi^{n,i}_{4,1})_{\nu_1,\nu_2}:=\E[(\wt \phi^{n,i}_{4,1})_{\nu_1,\nu_2}+(\wt \phi^{n,i}_{4,1})_{\nu_2,\nu_1}\mid\calf^n_{i-\la_n}]$ and
\begin{align*}
(\wt \phi^{n,i}_{4,1})_{\nu_1,\nu_2}&:=\partial^2_{\nu_1\nu_2} f(\chi^n_i)\iint_{(i-\la_n)\Del}^{(i+l_1-1)\Del} \frac{\Delta^n_i G_y(s)_{\nu_1}}{\tau_n}\si'(u(i\Del,y)_0^{(i-\la_n)\Del})\\
&\qtimes \iint_{(i-\la_n)\Del}^{s} \frac{\Delta^n_i G_z(r)_{\nu_2}}{\tau_n}\si'(u(i\Del,z)_0^{(i-\la_n)\Del})\\ 
&\cquad{\mathrel{\times}\iiint_0^{r-(i-\la_n)\Del} G_w(v)G_{w'}(v)\si(u(i\Del-v,y-w)_0^{(i-\la_n)\Del})}\\
&\cqquad{\mathrel{\times}	 	\si(u(i\Del-v,y-w')_0^{(i-\la_n)\Del})\,\La(\dd w,\dd w')\,\dd v\,W(\dd r,\dd z)\,W(\dd s,\dd y).}
\end{align*}

Analogous arguments also allow us to modify the $\La(\dd z,\dd z')\,\dd r$-integral in $(\wt \phi^{n,i}_{2,2})_{\nu_1,\nu_2}$ such that $D^{n,2,2}_3-D^{n,4,2}_3\limL0$, where $D^{n,4,2}_3(t):=\frac{\Delh}{2}\sumtla\sum_{\nu_1,\nu_2\in\cali} (\phi^{n,i}_{4,2})_{\nu_1,\nu_2}$ with  $(\phi^{n,i}_{4,2})_{\nu_1,\nu_2}:=\E[(\wt \phi^{n,i}_{4,2})_{\nu_1,\nu_2}\mid\calf^n_{i-\la_n}]$ and
\begin{align*}
(\wt \phi^{n,i}_{4,2})_{\nu_1,\nu_2}&:=\partial^2_{\nu_1\nu_2} f(\chi^n_i) \iiint_{(i-\la_n)\Del}^{(i+l_{12}-1)\Del} \frac{\Delta^n_i G_y(s)_{\nu_1}\Delta^n_i G_{y'}(s)_{\nu_2}}{\tau^2_n} \si'(u(i\Del,y)_0^{(i-\la_n)\Del})\\
&\qtimes\si'(u(i\Del,y')_0^{(i-\la_n)\Del})\iiint_0^{s-(i-\la_n)\Del} G_{z}(r)G_{z'}(r)\si(u(i\Del-r,y-z)_0^{(i-\la_n)\Del})\\
&\cqquad{\mathrel{\times}\si(u(i\Del-r,y-z')_0^{(i-\la_n)\Del})\,\La(\dd z,\dd z')\,\dd r \,\La(\dd y,\dd y')\,\dd s.}
\end{align*}

We continue with spatial approximations 
and replace all $y$- and $z$-variables (resp., all $y$- and $y'$-variables)
in the $u$-terms of $(\wt \phi^{n,i}_{4,1})_{\nu_1,\nu_2}$ [resp., $(\wt \phi^{n,i}_{4,2})_{\nu_1,\nu_2}$] by $x_{k_{12}}$. This is permitted because the $\La(\dd w,\dd w')\,\dd v$-integral [resp., the $\La(\dd z,\dd z')\,\dd r$-integral] is already of size $(\la_n\Del)^{1-\al/2}$ and the substitutions give us at least another factor of $\Del^{1/2-\al/4-\eps/2}$ [cf.\ \eqref{eq:yalpha}]. So the total error is $o(\Delh)$. Finally, by the same argument, we can change $\chi^n_i$ to $\wh\chi^n_i$, which is obtained by replacing $\si(u(i\Del,y)_0^{(i-\la_n)\Del})$  by $\si(u(i\Del,x_k)_0^{(i-\la_n)\Del})$ for the $(k,l)$th entry of \eqref{eq:chi}. 

We call the terms obtained from $(\wt \phi^{n,i}_{4,1})_{\nu_1,\nu_2}$ and $(\wt \phi^{n,i}_{4,2})_{\nu_1,\nu_2}$ in this way $(\wt \phi^{n,i}_{5,1})_{\nu_1,\nu_2}$ and $(\wt \phi^{n,i}_{5,2})_{\nu_1,\nu_2}$, respectively.
If we
apply the (stochastic) Fubini theorem and take the $\La(\dd w,\dd w')\,\dd v$-integral [resp., the $\La(\dd z,\dd z')\,\dd r$-integral] in $(\wt \phi^{n,i}_{5,1})_{\nu_1,\nu_2}$ [resp., $(\wt \phi^{n,i}_{5,2})_{\nu_1,\nu_2}$] to the outside, a careful inspection of the resulting expression and integration by parts reveal that
$(\wt \phi^{n,i}_5)_{\nu_1,\nu_2}:=(\wt \phi^{n,i}_{5,1})_{\nu_1,\nu_2}+(\wt \phi^{n,i}_{5,1})_{\nu_2,\nu_1}+(\wt \phi^{n,i}_{5,2})_{\nu_1,\nu_2}$ is nothing else but
\beq\label{eq:phi5}\begin{split}
	(\wt \phi^{n,i}_5)_{\nu_1,\nu_2}&=(\si')^2(u(i\Del,x_{k_{12}})_0^{(i-\la_n)\Del})
	\iiint_0^{(\la_n+l_{12}-1)\Del} G_z(r)G_{z'}(r) \partial^2_{\nu_1\nu_2} f(\wh\chi^n_i)\\
	&\quad\times \iint_{r+(i-\la_n)\Del}^{(i+l_1-1)\Del} \frac{\Delta^n_i G_y(s)_{\nu_1}}{\tau_n}\,W(\dd s,\dd y)\iint_{r+(i-\la_n)\Del}^{(i+l_2-1)\Del} \frac{\Delta^n_i G_y(s)_{\nu_2}}{\tau_n}\,W(\dd s,\dd y)\\
	&\quad\times\si(u(i\Del-r,x_{k_{12}}-z)_0^{(i-\la_n)\Del})\si(u(i\Del-r,x_{k_{12}}-z')_0^{(i-\la_n)\Del})\,\La(\dd z,\dd z')\,\dd r.
\end{split}\raisetag{3\baselineskip}\eeq
Thus, if $D^{n,5}_3(t):=\frac{\Delh}{2}\sumtla\sum_{\nu_1,\nu_2\in\cali} (\phi^{n,i}_5)_{\nu_1,\nu_2}$ and $(\phi^{n,i}_5)_{\nu_1,\nu_2}:=\E[(\wt \phi^{n,i}_5)_{\nu_1,\nu_2}\mid\calf^n_{i-\la_n}]$, then $D^{n,2}_3-D^{n,5}_3\limL0$ by what we have shown so far.

Obviously, we can replace the upper limit of the $\La(\dd z,\dd z')\,\dd r$-integral in \eqref{eq:phi5} by $\la_n\Del$.
Moreover, introducing intermediate points as in the proof of Lemma~\ref{lem:Dn21} [cf.\ the paragraphs after \eqref{eq:Kn3}], we can further show that the lower limit of both $W(\dd s,\dd y)$-integrals in \eqref{eq:phi5} can be substituted by $(i-\la_n)\Del$. The $\calf^n_{i-\la_n}$-conditional expectation of the resulting term equals
\begin{align*}
&
(\si')^2(u(i\Del,x_{k_{12}})_0^{(i-\la_n)\Del}) \wh \mu_f(\wh v^{n,i},\wh w^{n,i},\wh c^{n,i})_{\nu_1,\nu_2} \iiint_0^{\la_n\Del} G_z(r)G_{z'}(r)\\
&\quad\times\si(u(i\Del-r,x_{k_{12}}-z)_0^{(i-\la_n)\Del})\si(u(i\Del-r,x_{k_{12}}-z')_0^{(i-\la_n)\Del})\,\La(\dd z,\dd z')\,\dd r,
\end{align*}
where 
\begin{align*}
\wh v^{n,i}_{\nu_1,\nu_2}&:=\si^2(u(i\Del,x_{k_{12}})_0^{i\Del-\la_n\Del})\Pi^n_{|l_1-l_2|,x_{k_1}-x_{k_2}}([0,(\la_n+l_{12}-1)\Del]\times\R^d\times\R^d),\\
\wh w^{n,i}_{\nu_1,\nu_2}&:=\Pi^n_{|l_1-l_2|,x_{k_1}-x_{k_2}}([0,(\la_n+l_{12}-1)\Del]\times\R^d\times\R^d),\\
\wh c^{n,i}_{\nu_1,\nu_2}&:=\si(u(i\Del,x_{k_{12}})_0^{i\Del-\la_n\Del})\Pi^n_{|l_1-l_2|,x_{k_1}-x_{k_2}}([0,(\la_n+l_{12}-1)\Del]\times\R^d\times\R^d).
\end{align*}

In the limit $n\to\infty$, we may change the $\Pi^n$-terms to $\Ga_{|l_1-l_2}\bone_{k_1=k_2=k}$ and, upon replacing $i\Del$ by $s$, approximate the discrete sum defining $D^{n,5}_3$ by a continuous integral. The resulting expression is exactly $D^{n,6}_3$ from \eqref{eq:Dn6}, that is, we have $D^{n,5}_3- D^{n,6}_3\limL0$.
\epr

\bpr[Proof of Proposition~\ref{prop:cancel}]
Let $Z$ be the Gaussian variable introduced after \eqref{eq:LLN} and $Z'$ be independent of $Z$ and with the same law as $Z$ except that $w_k$ in \eqref{eq:cov-Z} is replaced by $w'_k$ (where $w,w'\in[0,\infty)^K$). 
Assuming that $Z$ and $Z'$ are defined on an auxiliary probability space $(\wt \Om,\wt \calf,\wt \P)$ (with expectation $\wt \E$), we deduce from Taylor's formula that
\beq\label{eq:Taylormu}\begin{split} \mu_f(w)-\mu_f(w') &= \wt\E[f(Z)-f(Z')] \\
	&= \sum_{\nu\in\cali} \wt\E[\partial_\nu f(Z')(Z-Z')_\nu]+\frac12\sum_{\nu_1,\nu_2\in\cali}\wt\E[\partial^2_{\nu_1\nu_2}f(Z')(Z-Z')_{\nu_1}(Z-Z')_{\nu_2}]\\
	&\quad+\frac16\sum_{\nu_1\nu_2\nu_3\in\cali}\wt\E[\partial^3_{\nu_1\nu_2\nu_3}f(Z'')(Z-Z')_{\nu_1}(Z-Z')_{\nu_2}(Z-Z')_{\nu_3}] \end{split}\eeq
for some $Z''$ between $Z$ and $Z'$. We apply this to the difference in \eqref{eq:last} by taking $Z:=Z(s)$ and $Z':=Z_n(s)$, where $Z(s)_{kl}:=\si(u(s,x_k))(Z_0)_{kl}$, $Z_n(s)_{kl}:=\si(u(s,x_k)_0^{s-\la_n\Del})(Z_0)_{kl}$, and $Z_0$ has the same law as $Z$ with $w=(1,\dots,1)$. With this choice, $Z-Z'$ is of size $(\la_n\Del)^{1/2-\al/4}$, so the third-order term in \eqref{eq:Taylormu} is negligible by \eqref{eq:help13}. Similarly, if we develop $\si(u(s,x_k))-\si(u(s,x_k)_0^{s-\la_n\Del})$ (which is part of $Z-Z'$) via Taylor's formula, then only the first-order term (resp., the first- and the second-order term) has to be considered further in the second (resp., the first) summand on the right-hand side of \eqref{eq:Taylormu}.
Thus, 
$H^n - (H^n_1+H^n_2+H^n_3+H^n_4)\limL0$, where 
\begin{align*}
H^n_1(t)&:=-\frac{1}{\Delh} \int_{\la_n\Del}^t \sum_{\nu\in\cali}  \si'(u(s,x_k)_0^{s-\la_n\Del})(u(s,x_k)-u(s,x_k)_0^{s-\la_n\Del}) \wt\E[\partial_\nu f(Z_n(s))(Z_0)_\nu] \,\dd s,\\
H^n_2(t)&:=-\frac{1}{2\Delh} \int_{\la_n\Del}^t \sum_{\nu\in\cali} \si''(u(s,x_k)_0^{s-\la_n\Del}) (u(s,x_k)-u(s,x_k)_0^{s-\la_n\Del})^2\wt\E[\partial_\nu f(Z_n(s))(Z_0)_\nu] \,\dd s,\\
H^n_3(t)&:=-\frac{1}{2\Delh} \int_{\la_n\Del}^t  \sum_{\nu_1,\nu_2\in\cali\colon k_1=k_2=k}(\si'(u(s,x_k)_0^{s-\la_n\Del}))^2(u(s,x_k)-u(s,x_k)_0^{s-\la_n\Del})^2 \\
&\qtimes\wt\E[\partial^2_{\nu_1\nu_2}f(Z_n(s))(Z_0)_{\nu_1}(Z_0)_{\nu_2}]\,\dd s,\\
H^n_4(t)&:=-\frac{1}{\Delh}\int_0^{\la_n\Del} \mu_f(\si^2(u(s,\un x)))\,\dd s. 
\end{align*}

We have $\E[(H^n_4)^\ast_T]\lec \Del^{1/2}\la_n \leq \Del^{1/2-a}\to0$ whenever $a$ is close to $\frac1{2\Theta}$. For the other terms, 
we claim that $u(s,x_k)-u(s,x_k)_0^{s-\la_n\Del}$ in $H^n_1(t)$, $H^n_2(t)$, and $H^n_3(t)$ may be replaced by 
\beq\label{eq:tildeW}  \iint_{s-\la_n\Del}^s G_{x_k-z}(s-r)\si(u(r,z)_0^{s-\la_n\Del})\,W(\dd r,\dd z), \eeq
leading to, say, $\wt H^n_1(t)$, $\wt H^n_2(t)$, and $\wt H^n_3(t)$. Indeed, to justify these replacements,  a size argument as in \eqref{eq:help13} suffices for $H^n_2$ and $H^n_3$. For $H^n_1$, notice that $H^n_1(t)-\wt H^n_1(t)=L^n_1(t)+L^n_2(t)$, where
\begin{align*}
L^n_1(t)&:=-\frac{1}{\Delh} \int_{\la_n\Del}^t \sum_{\nu\in\cali}  \wt\E[\partial_\nu f(Z_n(s))(Z_0)_\nu]\si'(u(s,x_k)_0^{s-\la_n\Del})\iint_{s-\la_n\Del}^s G_{x_k-z}(s-r)\\
&\qtimes\si'(u(r,z)_0^{s-\la_n\Del})\iint_{s-\la_n\Del}^r G_{z-w}(r-v)\si(u(v,w))\,W(\dd v,\dd w)\,W(\dd r,\dd z)  \,\dd s,\\
L^n_2(t)&:=-\frac{1}{2\Delh} \int_{\la_n\Del}^t \sum_{\nu\in\cali}  \wt\E[\partial_\nu f(Z_n(s)) (Z_0)_\nu]\si'(u(s,x_k)_0^{s-\la_n\Del})\\
&\qtimes\iint_{s-\la_n\Del}^s G_{x_k-z}(s-r)\si''(\psi_n(s,r,z))(u(r,z)-u(r,z)_0^{s-\la_n\Del})^2\,W(\dd r,\dd z) \,\dd s,
\end{align*}
and $\psi_n(s,r,z)$ is some intermediate value. The term $L^n_2$ is of size $\Del^{-1/2}((\la_n\Del)^{1/2-\al/4})^3$ and vanishes as $n\to\infty$. Next,  up to negligible boundary terms, $L^n_1$ can be written as 
\begin{align*} &-\frac{1}{\Delh}\sum_{i=\la_n+1}^{[t/\Del]} \int_{(i-1)\Del}^{i\Del} \sum_{\nu\in\cali}  \wt\E[\partial_\nu f(Z_n(s))(Z_0)_\nu]\si'(u(s,x_k)_0^{s-\la_n\Del})\iint_{s-\la_n\Del}^s G_{x_k-z}(s-r)\\
&\quad\times\si'(u(r,z)_0^{s-\la_n\Del})\iint_{s-\la_n\Del}^r G_{z-w}(r-v)\si(u(v,w))\,W(\dd v,\dd w)\,W(\dd r,\dd z) \,\dd s.
\end{align*}
Clearly, we may change $Z_n(s)$ [resp., $u(s,x_k)_0^{s-\la_n\Del}$] to $Z_n((i-1)\Del)$ [resp., $u(s,x_k)_0^{(i-1)\Del-\la_n\Del}$] as this gives a third $(\la_n\Del)^{1/2-\al/4}$-factor. Having done that,
we realize
that the $i$th summand has a vanishing $\calf^n_{i-1-\la_n}$-conditional expectation, so by a martingale size estimate, the previous display is of size $\la_n^{1/2}((\la_n\Del)^{1/2-\al/4})^2$, which goes to $0$ by \eqref{eq:comparison} and proves our claim above.

Returning to the main line of our proof and recalling the definition of $\ov \mu_f$ from \eqref{eq:ovmu} and that of $\ov v_n(s)$, $\ov w_n(s)$, and $\ov c_n(s)$ before Lemma~\ref{lem:easy}, we realize that $B^{n,8}_2$ from \eqref{eq:Bn8} is exactly equal to $-\wt H^n_1$, that is, $B^{n,8}_2+\wt H^n_1\equiv 0$. Moreover, if we expand the square of \eqref{eq:tildeW} in $\wt H^n_2(t)$ [resp., $\wt H^n_3(t)$] using the integration by parts formula, it is not hard to see that the quadratic variation part is equal to $-D^{n,6}_2(t)$ [resp., $-D^{n,6}_3(t)$] from Lemma~\ref{lem:Dn2456} (resp., Lemma~\ref{lem:Dn3}). That is, we have
\begin{align*}
\wt H^n_2(t)+D^{n,6}_2(t) &= -\frac{1}{2\Delh} \int_{\la_n\Del}^t \sum_{\nu\in\cali} \si''(u(s,x_k)_0^{s-\la_n\Del}) \wt\E[\partial_\nu f(Z_n(s))(Z_0)_\nu]\\
&\quad\times \iint_{s-\la_n\Del}^s G_{x_k-z}(s-r)\si(u(r,z)_0^{s-\la_n\Del})\\
&\qquad\times\iint_{s-\la_n\Del}^r G_{x_k-w}(s-v)\si(u(v,w)_0^{s-\la_n\Del})\,W(\dd v,\dd w) \,W(\dd r,\dd z) \,\dd s
\end{align*}
and
\begin{align*}
\wt H^n_3(t)+D^{n,6}_3(t) &= -\frac{1}{2\Delh} \int_{\la_n\Del}^t  \sum_{\nu_1,\nu_2\in\cali\colon k_1=k_2=k}(\si'(u(s,x_k)_0^{s-\la_n\Del}))^2\wt\E[\partial^2_{\nu_1\nu_2}f(Z_n(s))(Z_0)_{\nu_1}(Z_0)_{\nu_2}] \\
&\quad\times\iint_{s-\la_n\Del}^s G_{x_k-z}(s-r)\si(u(r,z)_0^{s-\la_n\Del})\\
&\qquad\times\iint_{s-\la_n\Del}^r G_{x_k-w}(s-v)\si(u(v,w)_0^{s-\la_n\Del})\,W(\dd v,\dd w)\,W(\dd r,\dd z)\,\dd s,
\end{align*}
respectively.
In both cases, we have a double $W$-integral inside a $\dd s$-integral, which is the same structure as encountered in the analysis of $L^n_1$ above. Thus, by analogous arguments, we derive $\wt H^n_2+D^{n,6}_2\limL0$ and $\wt H^n_3+D^{n,6}_3\limL0$. 
\epr

\appendix
\section{Auxiliary results}
We first state some basis properties of the measures $\Pi^n_{r,h}$ and $|\Pi^n_{r,h}|$ introduced in Section~\ref{sect:step3}. Some (weaker) variants thereof were proved in \cite[Section~B]{Chong19a}. For the reader's convenience, and because the proof below is much simpler than in the mentioned reference, we give full details here.
\blem\label{lem:Pin} 
Let $\al\in(0,1)$, $\Theta=1+\frac\al2$, and $a\in(0,1)$ and recall the definition of $\Ga_r$ from \eqref{eq:Ga-formula}. Then the following estimates hold uniformly in $n$, $r$, and $h$:
\begin{align}
|\Pi^n_{r,h}|([0,\infty)\times\R^d\times\R^d)&\lec 1, \label{eq:Pin1}\\
|\Pi^n_{r,h}|((\Del^{1-a},\infty)\times\R^d\times\R^d)&\lec \Del^{a\Theta},\label{eq:Theta}\\
\iiint_0^\infty (\lvert y\rvert^2+\lvert y'\rvert^2)\,|\Pi^n_{r,h}|(\dd s,\dd y,\dd y')&\lec \Del.\label{eq:y2}
\end{align}
Moreover, we have for every $n\in\N$ and $r\in\N_0$,
\beq \Pi^n_{r,0}([0,\infty)\times\R^d\times\R^d)=\Ga_r \label{eq:PinGa}\eeq
\elem
\bpr Let us consider \eqref{eq:y2} first, where it suffices
to analyze the integral against $|y|^2$ (the other is completely analogous). It can be split into $Z^n_1+Z^n_2$, where
\begin{align*}
Z^n_1&:=\frac{1}{\tau_n^2} \iiint_0^{\Del} G_y(s)|G_{y'+h}(s+r\Del)-G_{y'+h}(s+(r-1)\Del)||y|^2\,\La(\dd y,\dd y')\,\dd s,\\
Z^n_2&:=\frac 1 {\tau_n^2} \iiint_0^\infty |G_y(s+\Del)-G_y(s)||G_{y'+h}(s+(r+1)\Del)-G_{y'+h}(s+r\Del)||y|^2\,\La(\dd y,\dd y')\,\dd s.
\end{align*}
Similarly to \eqref{eq:sp-der-G},
\beq\label{eq:der-G} {|\partial_t G(t,x)|} = \frac{G(t,x)}{2 t}\biggl|\frac{|x|^2}{ t}-d\biggr|\leq\frac{G(t,\textstyle\frac{x}{\sqrt 2})}{2 t}\biggl(\ee^{-\frac{|x|^2}{4 t}}\frac{|x|^2}{ t}+d\ee^{-\frac{|x|^2}{4 t}}\biggr) \leq\frac{\frac4 \ee +d}{2 t}G(t,\textstyle\frac{x}{\sqrt 2}). \eeq
Therefore, using the fundamental theorem of calculus and eliminating the $\frac1{\sqrt{2}}$-factor in \eqref{eq:der-G} by a change of variables, we obtain
\beq\label{eq:help15}\begin{split}
Z^n_2 &= \frac 1 {\tau_n^2} \iiint_0^\infty \Biggl|\int_0^\Del \partial_t G(s+u,y)\,\dd u\int_0^\Del \partial_t G(s+r\Del+v,y'+h)\,\dd v \Biggr||y|^2\,\La(\dd y,\dd y')\,\dd s\\
&\leq\frac 1 {\tau_n^2}\int_0^\Del\int_0^\Del \iiint_0^\infty | \partial_t G(s+u,y) \partial_t G(s+r\Del+v,y'+h) ||y|^2\,\La(\dd y,\dd y')\,\dd s\,\dd u\,\dd v\\
&\lec\frac 1 {\tau_n^2}\int_0^\Del\int_0^\Del \int_0^\infty (s+u)^{-1}(s+r\Del+v)^{-1}\\
&\quad\times\int_{\R^d\times\R^d} G(s+u,y)G(s+r\Del+v,y'+h) |y|^2\,\La(\dd y,\dd y')\,\dd s\,\dd u\,\dd v.
\end{split}\eeq

Now recall from \cite[Chapter~I, Theorem~1.8]{Stein71} that multiplication with a monomial in the time domain corresponds to differentiation in the frequency domain. Therefore, together with \eqref{eq:isometry} and the rule $\calf \phi(\cdot+h)(\xi)= \ee^{2\pi \ii h\cdot\xi}$, we obtain
\begin{align*}
Z^n_2&\lec \frac1 {\tau_n^2}\int_0^\Del\int_0^\Del \int_0^\infty (s+u)^{-1}(s+v)^{-1}\\
&\quad\times\int_{\R^d}  \Biggl(\sum_{j=1}^d \partial^2_{\xi_j\xi_j}\calf G(s+u,\cdot)(\xi)\Biggr)\calf G(s+r\Del+v,\cdot)(\xi)\ee^{-2\pi \ii h\cdot\xi}\,\mu(\dd \xi)\,\dd s\,\dd u\,\dd v.
\end{align*}
By \eqref{eq:FG}, we have $\partial^2_{\xi_j\xi_j}\calf G(s,\cdot)(\xi) = 4\pi^2 s\ee^{-2\pi^2 s |\xi|^2}(4\pi^2 s\xi_j^2-1)$, which implies
that
\begin{align*}
Z^n_2&\lec \frac1 {\tau_n^2}\int_0^\Del\int_0^\Del \int_0^\infty (s+v)^{-1} \int_{\R^d} \ee^{-2\pi^2 (2s+v+u+r\Del) |\xi|^2}  \Bigl|(s+u)|\xi|^2-\textstyle\frac d {4\pi^2}\Bigr|\,\mu(\dd \xi)\,\dd s\,\dd u\,\dd v\\
&\lec\frac1 {\tau_n^2}\int_0^\Del\int_0^\Del \int_0^\infty (s+v)^{-1} \int_0^\infty \ee^{-2\pi^2 (2s+v+u) z^2}  \Bigl|(s+u)z^2-\textstyle\frac d {4\pi^2}\Bigr| z^{\al-1}\,\dd z\,\dd s\,\dd u\,\dd v\\
&\lec\frac1 {\tau_n^2}\int_0^\Del\int_0^\Del \int_0^\infty (s+v)^{-1} \Biggl(\int_0^{\sqrt{d/(4\pi^2(s+u))}} \ee^{-2\pi^2 (2s+v+u) z^2}  \textstyle\frac d {4\pi^2}z^{\al-1}\,\dd z\\
&\quad+\int_{\sqrt{d/(4\pi^2(s+u))}}^\infty \ee^{-2\pi^2 (2s+v+u) z^2}  (s+u)z^{\al+1}\,\dd z\Biggr)\,\dd s\,\dd u\,\dd v.
\end{align*}
Bounding the first exponential simply by $1$ and the second one by $\ee^{-2\pi^2(s+u)z^2}$ and substituting $2\pi^2(s+u)z^2\mapsto z^2$ for the second integral, we further derive
\begin{align*}
Z^n_2&\lec\frac1 {\tau_n^2}\int_0^\Del\int_0^\Del \int_0^\infty (s+v)^{-1} \Biggl(\int_0^{\sqrt{d/(4\pi^2(s+u))}}z^{\al-1}\,\dd z\\
&\quad+\int_{\sqrt{d/(4\pi^2(s+u))}}^\infty\ee^{-2\pi^2(s+u)z^2}  (s+u)z^{\al+1}\,\dd z\Biggr)\,\dd s\,\dd u\,\dd v\\
&\lec \frac1 {\tau_n^2}\int_0^\Del\int_0^\Del \int_0^\infty (s+v)^{-1} \Biggl((s+u)^{-\frac\al2}+\frac{(s+u)}{(2\pi^2(s+u))^{\frac\al2+1}} \int_{\sqrt{d/2}}^\infty \ee^{-z^2} z^{\al+1}\,\dd z\Biggr)\,\dd s\,\dd u\,\dd v\\
&\lec \frac1 {\tau_n^2}\int_0^\Del\int_0^\Del \int_0^\infty (s+v)^{-1} (s+u)^{-\frac\al2}\,\dd s\,\dd u\,\dd v.
\end{align*}

If $u\leq v$ (resp., $v<u$) in the previous line, bound $\int_0^\Del (s+v)^{-1}\,\dd v$ [resp., $\int_0^\Del (s+u)^{-\al/2}\,\dd v$] by $\Del(s+u)^{-1}$ [resp., $\Del (s+v)^{-\al/2}$]. This has the consequence that 
\begin{align*}
Z^n_2&\lec \frac\Del {\tau_n^2}\int_0^\Del \int_0^\infty  (s+u)^{-1-\frac\al2}\,\dd s\,\dd u \lec \frac\Del {\tau_n^2} \int_0^\infty \Bigl(s^{-\frac\al2}-(s+\Del)^{-\frac\al2}\Bigr)\,\dd s\\
&=\frac{\Del^{2-\frac\al2}} {\tau_n^2} \int_0^\infty \Bigl(t^{-\frac\al2}-(t+1)^{-\frac\al2}\Bigr)\,\dd t \lec \Del.
\end{align*}
Notice that the last integral is finite because the integrand behaves like $t^{-\al/2}$ at zero and like $t^{-1-\al/2}$ at infinity.

Concerning $Z^n_1$, let us consider the case where $r=0$ first. From \cite[(B.17)]{Chong19a}, we know that $\int_{\R^d} G_{y'+h}(s+r\Del)|y'-y|^{-\al}\,\dd y'\lec (s+r\Del)^{-\al/2}\leq s^{-\al/2}$ (uniformly in $r$, $h$, $y$, and $n$), and because $G(s,\cdot)$ is a Gaussian density, we have $\int_{\R^d} G_y(s)|y|^2\,\dd y \lec s$. Thus, using \eqref{eq:taun} for the last step,
$$Z^n_1=\frac1{\tau_n^2} \iiint_0^\Del G_y(s)G_{y'+h}(s)|y|^2|y-y'|^{-\al}\,\dd y\,\dd y'\,\dd s\lec \frac1{\tau_n^2}\int_0^\Del s^{1-\frac\al2}\,\dd s \lec \Del,$$ 
as desired. For $r\geq1$, we have, combining elements of the previous calculations, 
\begin{align*}
Z^n_1 &\leq \frac1{\tau_n^2}\int_0^\Del\iiint_0^\Del G(s,y)|\partial_t G(s+(r-1)\Del+v,y'+h)||y|^2\,\La(\dd y,\dd y')\,\dd s\,\dd v\\
&\lec \frac1{\tau_n^2}\int_0^\Del\iiint_0^\Del G(s,y) (s+(r-1)\Del+v)^{-1}G(s+(r-1)\Del+v,\textstyle\frac{y'+h}{\sqrt{2}})|y|^2\,\La(\dd y,\dd y')\,\dd s\,\dd v\\
&\lec\frac1{\tau_n^2}\int_0^\Del\int_0^\Del s\times s^{-1} \times s^{-\frac\al2}\,\dd s\,\dd v= \frac{\Del}{\tau_n^2}\int_0^\Del s^{-\frac\al2}\,\dd s \lec \Del,
\end{align*}
which completes the proof of \eqref{eq:y2}.

Property \eqref{eq:Pin1} can be proved in a similar fashion. Indeed, the left-hand side of \eqref{eq:Pin1} equals $Z^n_3+Z^n_4$, which are defined in the same way as $Z^n_1$ and $Z^n_2$, respectively, but with $|y|^2$ replaced by $1$. Then, from \eqref{eq:help15},
\begin{align*}
Z^n_4 &\lec \frac 1 {\tau_n^2}\int_0^\Del\int_0^\Del \int_0^\infty (s+u)^{-1}(s+r\Del+v)^{-1}\\
&\quad\times\int_{\R^d\times\R^d} G(s+u,y)G(s+r\Del+v,y'+h) \,\La(\dd y,\dd y')\,\dd s\,\dd u\,\dd v\\
&\lec \frac 1 {\tau_n^2}\int_0^\Del\int_0^\Del \int_0^\infty (s+u)^{-1}(s+r\Del+v)^{-1} \int_0^\infty \ee^{-2\pi^2(2s+u+v+r\Del)z^2}z^{\al-1}\,\dd z\,\dd s\,\dd u\,\dd v\\
&\lec \frac 1 {\tau_n^2}\int_0^\Del\int_0^\Del \int_0^\infty (s+u)^{-1}(s+v)^{-1}(2s+u+v)^{-\frac\al2}\,\dd s\,\dd u\,\dd v\\
&\lec \frac{1}{\tau_n^2}\int_0^\infty \int_0^\Del\int_0^\Del (s+v)^{-1-\frac\al4}(s+u)^{-1-\frac\al4} \,\dd u\,\dd v\,\dd s\\
&\lec\frac{1}{\tau_n^2} \int_0^\infty (s^{-\frac\al4}-(s+\Del)^{-\frac\al4})^2\,\dd s\lec \frac{\Del^{1-\frac\al2}}{\tau_n^2} \int_0^\infty (t^{-\frac\al4}-(t+1)^{-\frac\al4})^2\,\dd t\lec 1.
\end{align*}
We leave it to the reader to prove $Z^n_3\lec 1$, which finishes the proof of \eqref{eq:Pin1}.

Regarding \eqref{eq:Theta}, we deduce from similar considerations that the left-hand side is
\begin{align*}
&\lec \frac 1 {\tau_n^2}\int_0^\Del\int_0^\Del \int_{\Del^{1-a}-\Del}^\infty (s+u)^{-1}(s+r\Del+v)^{-1}\\
&\quad\times\int_{\R^d\times\R^d} G(s+u,y)G(s+r\Del+v,y'+h) \,\La(\dd y,\dd y')\,\dd s\,\dd u\,\dd v\\
&\lec \frac 1 {\tau_n^2}\int_{\frac12\Del^{1-a}}^\infty (s^{-\frac\al4}-(s+\Del)^{-\frac\al4})^2\,\dd s \lec  \frac{\Del^{1-\frac\al2}}{\tau_n^2} \int_{\frac12 \Del^{-a}}^\infty (t^{-\frac\al4}-(t+1)^{-\frac\al4})^2\,\dd t \lec \Del^{\Theta a},
\end{align*}
because $t\mapsto (t^{-\al/4}-(t+1)^{-\al/4})^2$ behaves like $t^{-2-\al/2} = t^{-1-\Theta}$ at infinity.

Finally, for \eqref{eq:PinGa}, observe that $\Pi^n_{r,0}([0,\infty)\times\R^d\times\R^d)=Z^n_5+Z^n_6$, where \begin{align*}
Z^n_5&:=\frac{1}{\tau_n^2} \iiint_0^{\Del} G_y(s)(G_{y'}(s+r\Del)-G_{y'}(s+(r-1)\Del))\,\La(\dd y,\dd y')\,\dd s,\\
Z^n_6&:=\frac 1 {\tau_n^2} \iiint_0^\infty (G_y(s+\Del)-G_y(s))(G_{y'}(s+(r+1)\Del)-G_{y'}(s+r\Del))\,\La(\dd y,\dd y')\,\dd s.
\end{align*}
Using \eqref{eq:isometry} and \eqref{eq:FG} for the first step, changing to polar coordinates in the second step, and employing the identity $\int_0^\infty \ee^{-uz^2} z^{\al-1}\,\dd z=\frac12\Ga(\frac\al2)u^{-\al/2}$ in the third step, we have
\begin{align*}
Z^n_6 &= \frac{1}{\tau_n^2} \iint_0^\infty (\ee^{-2\pi^2(s+\Del)|\xi|^2}-\ee^{-2\pi^2s|\xi|^2})(\ee^{-2\pi^2(s+(r+1)\Del)|\xi|^2}-\ee^{-2\pi^2(s+r\Del)|\xi|^2}) \,\mu(\dd \xi)\,\dd s\\
&=\frac{2\pi^{\frac d2}}{\Ga(\frac d2)\tau_n^2} \int_0^\infty \int_0^\infty (\ee^{-2\pi^2(s+\Del)z^2}-\ee^{-2\pi^2sz^2})(\ee^{-2\pi^2(s+(r+1)\Del)z^2}-\ee^{-2\pi^2(s+r\Del)z^2})z^{\al-1}\,\dd z\,\dd s\\
&=\frac{\pi^{\frac d2-\al}\Ga(\frac\al2)}{2^{\frac\al2}\Ga(\frac d2)\tau_n^2}\int_0^\infty ((2s+(r+2)\Del)^{-\frac\al2}  - 2(2s+(r+1)\Del)^{-\frac\al2} + (2s+r\Del)^{-\frac\al2})\,\dd s.
\end{align*}
By \eqref{eq:taun} and a change of variables $s\mapsto t\Del$,
\beq\label{eq:help16}\begin{split}
 Z^n_6 &= (1-\textstyle\frac\al2)\displaystyle\int_0^\infty ((2t+(r+2))^{-\frac\al2}  - 2(2t+(r+1))^{-\frac\al2} + (2t+r)^{-\frac\al2})\,\dd t \\
 &= -\frac12((r+2)^{1-\frac\al2}-2(r+1)^{1-\frac\al2}+r^{1-\frac\al2}).
\end{split}\eeq
A similar calculation shows that for $r\geq1$, 
\begin{align*}
Z^n_5&= \frac{2\pi^{\frac d2}}{\Ga(\frac d2)\tau_n^2} \int_0^{\Del} \int_0^\infty \ee^{-2\pi^2sz^2}(\ee^{-2\pi^2(s+r\Del)z^2}-\ee^{-2\pi^2(s+(r-1)\Del)z^2})z^{\al-1}\,\dd z\,\dd s\\
&=\frac{\pi^{\frac d2-\al}\Ga(\frac\al2)}{2^{\frac\al2}\Ga(\frac d2)\tau_n^2}\int_0^\Del ((2s+r\Del)^{-\frac\al2}-(2s+(r-1)\Del)^{-\frac\al2})\,\dd s\\
&=(1-\textstyle\frac\al2)\displaystyle\int_0^1 ((2t+r)^{-\frac\al2}-(2t+(r-1))^{-\frac\al2})\,\dd t \\
&= \frac12((r+2)^{1-\frac\al2}-r^{1-\frac\al2}-(r+1)^{1-\frac\al2}+(r-1)^{1-\frac\al2}),
\end{align*}
 which, together with \eqref{eq:help16}, shows \eqref{eq:PinGa}  for $r\geq1$. For $r=0$, we have
 \begin{align*}
 Z^n_5&=\frac{2\pi^{\frac d2}}{\Ga(\frac d2)\tau_n^2} \int_0^{\Del} \int_0^\infty \ee^{-2\pi^2sz^2}\ee^{-2\pi^2(s+\Del)z^2}z^{\al-1}\,\dd z\,\dd s=\frac{\pi^{\frac d2-\al}\Ga(\frac\al2)}{2^{\frac\al2}\Ga(\frac d2)\tau_n^2}\int_0^\Del (2s)^{-\frac\al2}\,\dd s\\
 &=(1-\textstyle\frac\al2)\displaystyle\int_0^1 (2t)^{-\frac\al2}\,\dd t = 2^{-\frac\al2}.
 \end{align*}
 This shows $\Pi^n_{0,0}([0,\infty)\times\R^d\times\R^d) = 1$ and completes the proof of the lemma.
\epr

Next, we gather some moment bounds on the solution to \eqref{eq:SHE-si}.
\blem\label{lem:momincr} Let $\si\colon \R\to\R$ be a globally Lipschitz function and $u_0$ be bounded and $(1-\frac\al2)$-Hölder continuous. For any $p\in(0,\infty)$, $T\in(0,\infty)$, and $\eps\in(0,1-\frac\al2)$,  the solution $u$ to \eqref{eq:SHE-si} satisfies
\beq\label{eq:increments} \E[|u(t,x)-u(t+\tau,x)|^p]^{\frac1p} \lec \tau^{\frac12-\frac\al4},\qquad \E[|u(t,x)-u(t,x+h)|^p]^{\frac1p} \lec |h|^{1-\frac\al2-\eps}, \eeq
uniformly in $t,\tau\in[0,T]$ and $x,h\in\R^d$ with $\lvert h\rvert \leq 1$.
\elem
\bpr
The first bound in \eqref{eq:increments} follows easily from \cite[Theorem~5.1.2 (ii)]{Lunardi95} (for the deterministic part) and \eqref{eq:Pin1} together with a standard size estimate (for the stochastic part).  The second bound is proved in \cite[Equation (19)]{Sanzsole02}.
\epr

The last result concerns orthogonality of multiple Wiener integrals. It is not a new result, yet we were not able to find the specific statement we need in the literature.
\blem\label{lem:chaos}
Assume that $g\colon \R^N\to\R$ satisfies $\E[g(cX)]<\infty$ for all $c>0$, where $X$ is a standard Gaussian in $\R^N$. Further suppose that $h_1\colon [0,\infty)\times\R^d\to\R$ and $h_2\colon ([0,\infty)\times\R^d)^N \to\R$ are such that 
\[(s,y,y')\mapsto |h_1(s,y)h_1(s,y')|\]
and
\[ (s_1,\dots,s_N,y_1,y'_1,\dots,y_N,y'_N)\mapsto |h_2(s_1,\dots,s_N,y_1,\dots,y_N)h_2(s_1,\dots,s_N,y'_1,\dots,y'_N)|\] 
are integrable with respect to $\La(\dd y,\dd y')\,\dd s$ and $\La(\dd y_1,\dd y'_1)\,\dd s_1\dotsm\,\La(\dd y_N,\dd y'_N)\,\dd s_N$, respectively, where $\La$ is the measure from \eqref{eq:cov}. For given numbers $0\leq a \leq b$ and $0\leq a_i\leq b_i$, $i=1,\dots,N$, if there exists $j\in\{1,\dots,N\}$ with $b_j\leq a$, then
\beq\label{eq:exp}\begin{split} &\E\Biggl[g\biggl(\iint_a^b h_1(s,y)\,W(\dd s,\dd y)\biggr) \\
&\quad\times\iint_{a_1}^{b_1} \dotsi \iint_{a_N}^{b_N\wedge s_{N-1}} h_2(s_1,\dots,s_N,y_1,\dots,y_N)\,W(\dd s_N,\dd y_N)\dotsm\,W(\dd s_1,\dd y_1)\Biggr]=0. \end{split}\eeq
\elem
\bpr
By Wiener chaos expansion (see \cite[Theorem~1.1.2]{Nualart06}), one can write $g(\iint_a^b h_1(s,y)\,W(\dd s,\dd y))$ as a series $\sum_{j=0}^\infty I_j(g_j)$ of multiple integrals. As multiple integrals of different orders are orthogonal to each other (see \cite[Theorem~1.1.1 and Proposition~1.1.4]{Nualart06}), only the $N$th term remains in \eqref{eq:exp}, that is, the expectation on the left-hand side equals
\beq\label{eq:exp-2} \E\Biggl[I_N(g_N) \iint_{a_1}^{b_1} \dotsi \iint_{a_N}^{b_N\wedge s_{N-1}} h_2(s_1,\dots,s_N,y_1,\dots,y_N)\,W(\dd s_N,\dd y_N)\dotsm\,W(\dd s_1,\dd y_1)\Biggr]. \eeq
In order to evaluate this, we can first condition on $\calf_{b_{j_0}}$, where $b_{j_0}$ is the smallest of all $b_j$'s satisfying $b_j\leq a$ (the existence of such a $b_j$ is guaranteed by hypothesis). Clearly, all the functions $g_j$ (in particular, $g_N$) are supported on $([a,b]\times\R^d)^j$, so the law of $I_N(g_N)$ is not affected by conditioning on $\calf_{b_{j_0}}$ (because $b_{j_0}\leq a$ and $W$ has independent increments in time). At the same time, the $N$-fold integral of $h_2$ is only an $(N-j_0)$-fold iterated integral conditionally on $\calf_{b_{j_0}}$, because the $j_0$ innermost integrals are $\calf_{b_{j_0}}$-measurable. So the orthogonality property of multiple integrals implies that the $\calf_{b_{j_0}}$-conditional expectation and hence the whole expectation in \eqref{eq:exp-2} is zero.
\epr

	\bibliographystyle{plainnat}
	\bibliography{bib-VolaEstimation}

\begin{thebibliography}{28}
\providecommand{\natexlab}[1]{#1}
\providecommand{\url}[1]{\texttt{#1}}
\expandafter\ifx\csname urlstyle\endcsname\relax
  \providecommand{\doi}[1]{doi: #1}\else
  \providecommand{\doi}{doi: \begingroup \urlstyle{rm}\Url}\fi

\bibitem[A{\"i}t-Sahalia and Jacod(2014)]{AitSahalia14}
Y.~A{\"i}t-Sahalia and J.~Jacod.
\newblock \emph{High-Frequency Financial Econometrics}.
\newblock Princeton University Press, Princeton, 2014.

\bibitem[Barndorff-Nielsen et~al.(2011)Barndorff-Nielsen, Corcuera, and
  Podolskij]{BN11}
O.E. Barndorff-Nielsen, J.M. Corcuera, and M.~Podolskij.
\newblock Multipower variation for {B}rownian semistationary processes.
\newblock \emph{Bernoulli}, 17\penalty0 (4):\penalty0 1159--1194, 2011.

\bibitem[Bibinger and Trabs(2018)]{Bibinger17}
M.~Bibinger and M.~Trabs.
\newblock Volatility estimation for stochastic {PDE}s using high-frequency
  observations.
\newblock arXiv:1710.03519, 2018.

\bibitem[Bibinger and Trabs(2019)]{Bibinger19}
M.~Bibinger and M.~Trabs.
\newblock On central limit theorems for power variations of the solution to the
  stochastic heat equation.
\newblock arXiv:1901.01026, 2019.

\bibitem[Chong(2019{\natexlab{a}})]{Chong19}
C.~Chong.
\newblock High-frequency analysis of parabolic stochastic {PDE}s.
\newblock \emph{Ann. Statist.}, 2019{\natexlab{a}}.
\newblock Forthcoming.

\bibitem[Chong(2019{\natexlab{b}})]{Chong19a}
C.~Chong.
\newblock Supplement to ``{H}igh-frequency analysis of parabolic stochastic
  {PDE}s''.
\newblock arXiv:1806.06959, 2019{\natexlab{b}}.

\bibitem[Chong(2019{\natexlab{c}})]{Chong19c}
C.~Chong.
\newblock High-frequency analysis of parabolic stochastic {PDE}s with
  multiplicative noise: {P}art {II}.
\newblock In preparation, 2019{\natexlab{c}}.

\bibitem[Cialenco(2018)]{Cialenco18}
I.~Cialenco.
\newblock Statistical inference for {SPDE}s: an overview.
\newblock \emph{Stat. Inference Stoch. Process.}, 21\penalty0 (2):\penalty0
  309--329, 2018.

\bibitem[Cialenco and Huang(2019)]{Cialenco17}
I.~Cialenco and Y.~Huang.
\newblock A note on parameter estimation for discretely sampled {SPDE}s.
\newblock \emph{Stoch. Dyn.}, 2019.
\newblock Forthcoming.

\bibitem[Corcuera et~al.(2006)Corcuera, Nualart, and Woerner]{Corcuera06}
J.M. Corcuera, D.~Nualart, and J.H.C. Woerner.
\newblock Power variation of some integral fractional processes.
\newblock \emph{Bernoulli}, 12\penalty0 (4):\penalty0 713--735, 2006.

\bibitem[Corcuera et~al.(2013)Corcuera, Hedevang, Pakkanen, and
  Podolskij]{Corcuera13}
J.M. Corcuera, E.~Hedevang, M.S. Pakkanen, and M.~Podolskij.
\newblock Asymptotic theory for {B}rownian semi-stationary processes with
  application to turbulence.
\newblock \emph{Stochastic Process. Appl.}, 123\penalty0 (7):\penalty0
  2552--2574, 2013.

\bibitem[Dalang(1999)]{Dalang99}
R.C. Dalang.
\newblock Extending martingale measures stochastic integral with applications
  to spatially homogeneous {S.P.D.E's}.
\newblock \emph{Electron. J. Probab.}, 4, 1999.
\newblock 29 pages.

\bibitem[Foondun et~al.(2015)Foondun, Khoshnevisan, and Mahboubi]{Foondun15}
M.~Foondun, D.~Khoshnevisan, and P.~Mahboubi.
\newblock Analysis of the gradient of the solution to a stochastic heat
  equation via fractional {B}rownian motion.
\newblock \emph{Stoch. Partial Differ. Equ. Anal. Comput.}, 3\penalty0
  (2):\penalty0 133--158, 2015.

\bibitem[Huang et~al.(2018)Huang, Nualart, and Viitasaari]{Huang18}
J.~Huang, D.~Nualart, and L.~Viitasaari.
\newblock A central limit theorem for the stochastic heat equation.
\newblock arXiv:1810.09492, 2018.

\bibitem[Huang et~al.(2019)Huang, Nualart, Viitasaari, and Zheng]{Huang19}
J.~Huang, D.~Nualart, L.~Viitasaari, and G.~Zheng.
\newblock Gaussian fluctuations for the stochastic heat equation with colored
  noise.
\newblock \emph{Stoch. Partial Differ. Equ. Anal. Comput.}, 2019.
\newblock Forthcoming.

\bibitem[Istas and Lang(1997)]{Istas97}
J.~Istas and G.~Lang.
\newblock Quadratic variations and estimation of the local {H}{\"o}lder index
  of a {G}aussian process.
\newblock \emph{Ann. Inst. H. Poincar{\'e} Probab. Statist.}, 33\penalty0
  (4):\penalty0 407--436, 1997.

\bibitem[Jacod(1997)]{Jacod97}
J.~Jacod.
\newblock On continuous conditional {G}aussian martingales and stable
  convergence in law.
\newblock In J.~Az{\'e}ma, M.~Emery, and M.~Yor, editors, \emph{S{\'e}minaire
  de Probabilit{\'e}s XXXI}, pages 232--246. Springer, Berlin, 1997.

\bibitem[Jacod and Protter(2012)]{Jacod12}
J.~Jacod and P.~Protter.
\newblock \emph{Discretization of Processes}.
\newblock Springer, Berlin, 2012.

\bibitem[Lunardi(1995)]{Lunardi95}
A.~Lunardi.
\newblock \emph{Analytic Semigroups and Optimal Regularity in Parabolic
  Problems}.
\newblock Birkh{\"a}user, Basel, 1995.

\bibitem[Nualart(2006)]{Nualart06}
D.~Nualart.
\newblock \emph{The Malliavin Calculus and Related Topics}.
\newblock Springer, Berlin, 2nd edition, 2006.

\bibitem[Podolskij(2015)]{Podolskij15}
M.~Podolskij.
\newblock Ambit fields: Survey and new challenges.
\newblock In R.H. Mena, J.C. Pardo, V.~Rivero, and G.U. Bravo, editors,
  \emph{XI Symposium on Probability and Stochastic Processes}, pages 241--279.
  Springer, Cham, 2015.

\bibitem[Podolskij and Vetter(2010)]{Podolskij10}
M.~Podolskij and M.~Vetter.
\newblock Understanding limit theorems for semimartingales: a short survey.
\newblock \emph{Statist. Neerlandica}, 64\penalty0 (3):\penalty0 329--351,
  2010.

\bibitem[Posp{\'i}{\v s}il and Tribe(2007)]{Pospisil07}
J.~Posp{\'i}{\v s}il and R.~Tribe.
\newblock Parameter estimates and exact variations for stochastic heat
  equations driven by space-time white noise.
\newblock \emph{Stoch. Anal. Appl.}, 25\penalty0 (3):\penalty0 593--611, 2007.

\bibitem[Sanz-Sol{\'e} and Sarr{\`a}(2002)]{Sanzsole02}
M.~Sanz-Sol{\'e} and M.~Sarr{\`a}.
\newblock H{\"o}lder continuity for the stochastic heat equation with spatially
  correlated noise.
\newblock In R.C. Dalang, M.~Dozzi, and F.~Russo, editors, \emph{Seminar on
  Stochastic Analysis, Random Fields and Applications III}. Birkh{\"a}user,
  Basel, 2002.

\bibitem[Stein(1970)]{Stein70}
E.M. Stein.
\newblock \emph{Singular Integrals and Differentiability Properties of
  Functions}.
\newblock Princeton University Press, Princeton, 1970.

\bibitem[Stein and Weiss(1971)]{Stein71}
E.M. Stein and G.~Weiss.
\newblock \emph{Introduction to Fourier Analysis on Euclidean Spaces}.
\newblock Princeton University Press, Princeton, NJ, 1971.

\bibitem[Swanson(2007)]{Swanson07}
J.~Swanson.
\newblock Variations of the solution to a stochastic heat equation.
\newblock \emph{Ann. Probab.}, 35\penalty0 (6):\penalty0 2122--2159, 2007.

\bibitem[Walsh(1981)]{Walsh81}
J.B. Walsh.
\newblock A stochastic model of neural response.
\newblock \emph{Adv. Appl. Probab.}, 13\penalty0 (2):\penalty0 231--281, 1981.

\end{thebibliography}

\end{document}